\documentclass[12pt,a4paper]{article}
\usepackage{amssymb,amsmath,amsthm,euscript}
\usepackage{leftidx}
\usepackage[matrix,arrow,curve]{xy}

\usepackage[utf8]{inputenc}
\usepackage[english]{babel}

\newcommand{\chara}{\mathop{\mathrm{char}}\nolimits}
\newcommand{\End}{\mathop{\mathrm{End}}\nolimits}

\newcommand{\Hom}{\mathop{\mathrm{Hom}}\nolimits}
\newcommand{\bfhom}{\mathop{\mathrm{\bf hom}}\nolimits}
\newcommand{\bfend}{\mathop{\mathrm{\bf end}}\nolimits}
\newcommand{\bfcohom}{\mathop{\mathrm{\bf cohom}}\nolimits}
\newcommand{\bfcoend}{\mathop{\mathrm{\bf coend}}\nolimits}

\newcommand{\id}{\mathop{\mathrm{id}}\nolimits}

\newcommand{\isoright}{\xrightarrow{\smash{\raisebox{-0.65ex}{\ensuremath{\sim}}}}}

\numberwithin{equation}{section}

\newcommand{\op}{{\mathrm{op}}}

\renewcommand{\le}{\leqslant}
\renewcommand{\ge}{\geqslant}

\newcommand{\A}{{\mathcal A}}
\newcommand{\B}{{\mathcal B}}

\newcommand{\gR}{{\mathfrak R}}

\newcommand{\gA}{{\mathfrak A}}

\newcommand{\wt}{\widetilde}
\newcommand{\wh}{\widehat}

\newtheorem{Th}{Theorem}[section]
\newtheorem{Lem}[Th]{Lemma}
\newtheorem{Prop}[Th]{Proposition}
\newtheorem{Cor}[Th]{Corollary}

\newtheorem{Def}[Th]{Definition}
\newtheorem{Rem}[Th]{Remark}
\newtheorem{War}[Th]{Warning}
\newcommand{\BB}{{\mathbb{B}}}

\newcommand{\RR}{{\mathbb{R}}}
\newcommand{\CC}{{\mathbb{C}}}
\newcommand{\ZZ}{{\mathbb{Z}}}
\newcommand{\MM}{{\mathbb{M}}}
\newcommand{\OO}{{\mathbb{O}}}
\newcommand{\NN}{{\mathbb{N}}}
\newcommand{\KK}{{\mathbb{K}}}
\newcommand{\FF}{{\mathbb{F}}}

\newcommand{\SSS}{{\mathbb{S}}}

\newcommand{\Set}{{\mathbf{Set}}}
\newcommand{\ManInf}{{\mathbf{Man}^\infty}}
\newcommand{\SLManInf}{{\mathbf{SLMan}^\infty}}

\newcommand{\AlgSet}{{\mathbf{AlgSet}}}
\newcommand{\SLSet}{{\mathbf{SLSet}}}
\newcommand{\SLAlgSet}{{\mathbf{SLAlgSet}}}

\newcommand{\Vect}{{\mathbf{Vect}}}
\newcommand{\FVect}{{\mathbf{FVect}}}

\newcommand{\Alg}{{\mathbf{Alg}}}
\newcommand{\CommAlg}{{\mathbf{CommAlg}}}
\newcommand{\FAlg}{{\mathbf{FAlg}}}
\newcommand{\GrAlg}{{\mathbf{GrAlg}}}
\newcommand{\ZZGrAlg}{{\text{$\ZZ$-$\mathbf{GrAlg}$}}}

\newcommand{\QA}{{\mathbf{QA}}}
\newcommand{\QAsc}{{\mathbf{QA}}_{\mathrm sc}}
\newcommand{\FQA}{{\mathbf{FQA}}}
\newcommand{\FQAsc}{{\mathbf{FQA}}_{\mathrm sc}}

\newcommand{\bfB}{{\mathbf{B}}}
\newcommand{\bfC}{{\mathbf{C}}}
\newcommand{\bfD}{{\mathbf{D}}}
\newcommand{\bfP}{{\mathbf{P}}}
\newcommand{\bfQ}{{\mathbf{Q}}}

\newcommand{\Mon}{\mathop{\mathbf{Mon}}\nolimits}

\newcommand{\Comon}{\mathop{\mathbf{Comon}}\nolimits}

\newcommand{\Bimon}{\mathop{\mathbf{Bimon}}\nolimits}
\newcommand{\Lact}{\mathop{\mathbf{Lact}}\nolimits}
\newcommand{\Rep}{\mathop{\mathbf{Rep}}\nolimits}
\newcommand{\Corep}{\mathop{\mathbf{Corep}}\nolimits}
\newcommand{\Lcoact}{\mathop{\mathbf{Lcoact}}\nolimits}
\newcommand{\eval}{\mathop{\mathrm{ev}}\nolimits}
\newcommand{\coev}{\mathop{\mathrm{coev}}\nolimits}

\usepackage{hyperref}

\newcounter{bbcount}[subsection]
\renewcommand{\thebbcount}{{\thesection}.\arabic{subsection}.\arabic{bbcount}}
\newcommand{\bbo}[1]{\noindent\refstepcounter{bbcount}{\bf\thebbcount.}{ \bfseries #1}}
\newcommand{\bb}[1]{\vspace{2mm}\noindent\refstepcounter{bbcount}{\bf\thebbcount.}{ \bfseries #1}}

\title{General representation theory in relatively closed monoidal categories}
\author{Alexey Silantyev\thanks{aleksejsilantjev@gmail.com}}
\date{}

\begin{document}

\maketitle

\vspace{-5mm}
\begin{center}
{\it Joint Institute for Nuclear Research, 141980 Dubna, Moscow region, Russia} \\
{\it State University ``Dubna''{}, 141980 Dubna, Moscow region, Russia} \\
\end{center}
\vspace{5mm}

\begin{abstract}
 We apply the notion of relative adjoint functor to generalise closed monoidal categories. We define representations in such categories and give their relation with left actions of monoids. The translation of these representations under lax monoidal functors is investigated. We introduce tensor product of representations of bimonoids as a functorial binary operation and show how symmetric lax monoidal functors act on this product. Finally we apply the general theory to classical and quantum representations.
\end{abstract}



\tableofcontents

\section{Introduction}

The classical representation theory deals with representations of algebras, groups and other algebraic objects which can be represented by linear operators on vector spaces. An idea of quantum representations goes back to the works of Yuri Manin~\cite{Manin87,Manin88,ManinBook91}, where he studied the category of finitely generated quadratic algebras. He proposed to regard the objects of the opposite category as a `non-commutative' (or `quantum') generalisation of finite-dimensional vector spaces and called them {\it quantum linear spaces}.
 
In~\cite{Sqrt} we introduced a notion of quantum representation, which generalises finite-dimensional representations of finite-dimensional algebras and algebraic groups to the case of quantum linear representation space. To define and investigate quantum representations we used there a general approach to representations in monoidal categories. In the present paper we continue to develop this approach, which we call {\it general representation theory}.

A representation of an algebra $\gA$ on a vector space $V$ can be regarded in two ways: as an action $a\colon\gA\otimes V\to V$ or as an algebra homomorphism $\rho\colon\gA\to\bfend(V)$, where $\bfend(V)$ is the algebra of linear operators on $V$. We use the word `representation' for such morphisms $\rho$. This is a basic example of the general notion of representation in a monoidal category for the category of vector spaces, where the role of the monoidal product is played by the tensor product $\otimes$.

The notion of action is generalised for any (symmetric) monoidal category $(\bfC,\otimes)$ as an action of a monoid (we consider the case of symmetric monoidal categories only). This notion is equivalent to the notion of representation iff $(\bfC,\otimes)$ is closed. In this case we can use the internal hom-functor $\bfhom\colon\bfC^\op\times\bfC\to\bfC$ to define representations and relate them with the (left) actions. In~\cite[\S~3]{Sqrt} we considered more general situation, when the functor $-\otimes V\colon\bfC\to\bfC$ has a right adjoint for each $V\in\bfP$, where $\bfP$ is a full subcategory of $\bfC$. We used the theorem on adjunction with a parameter~\cite[\S~4.7, Th.~3]{Mcl} to obtain a generalised internal hom-functor $\bfhom\colon\bfP^\op\times\bfC\to\bfC$ with a natural isomorphism
 $\alpha_{X,V,Z}\colon\Hom(X\otimes V,Z)\isoright\Hom\big(X,\bfhom(V,Z)\big)$, where $X\in\bfC$, $V\in\bfP$, $Z\in\bfC$.

But it turned out that this condition is too strong. We need the isomorphisms $\alpha_{X,V,Z}$ for $Z\in\bfP$ only. This can not be formulated in terms of usual adjoint functors, so we use the notion of relative adjoints introduced by Friedrich Ulmer in~\cite{Ulm}.
We require that for each $V\in\bfP$ the functor $-\otimes V\colon\bfC\to\bfC$ has a right adjoint relative to $\bfP$ and then we say that the symmetric monoidal category $(\bfC,\otimes)$ is {\it closed relative} to the {\it parametrising subcategory} $\bfP\subset\bfC$. This approach allows us to generalise the results of~\cite{Sqrt} to the monoidal category $\bfC$ opposite to the category of all the $\ZZ$-graded algebras with the Manin product `$\circ$' and the parametrising subcategory $\bfP$ consisting of the quantum linear spaces.

In the situations such as Quantum Representation Theory it is more convenient to work in the opposite category to the category $\bfC$. By this reason we double all the important results for this case by using categorical duality: we consider monoidal categories coclosed relative to a full subcategory, comonoids, corepresentations etc.

The paper is organised as follows. In the section~\ref{sec2} we generalise the notion of (co)closed monoidal category by means of relative adjoints. Subsection~\ref{secNot} is preliminary. In Subsection~\ref{secYL} we formulate some corollaries of Yoneda Lemma. Subsection~\ref{secRelAdj} is devoted to the theory of relative adjunctions. In Subsection~\ref{secRelCl} we define relatively (co)closed monoidal category and describe internal (co)hom and (co)end for this case. Section~\ref{secGRT} is devoted to (co)representations in relatively (co)closed monoidal categories. In Subsection~\ref{secRep} we introduce (co)representations of (co)monoids, their morphisms and give their relation with (co)actions. In Subsection~\ref{secTP} we describe tensor product of (co)representations of bimonoids. The translation of the (co)representations and their tensor products under monoidal functors is investigated in Subsection~\ref{secTrans}. Then, in Section~\ref{secEx}, we give the basic examples. Classical representations of various types of groups, algebras and their mixed (semi-linear) versions is described in Subsection~\ref{secCR} in terms of the general representation theory. In Subsection~\ref{secQR} we apply the general theory to the quantum representations.

\vspace{3mm}
{\it Acknowledgements}. 
The author thanks V. Rubtsov for useful advice.

\section{Relative adjunctions and relatively closed monoidal categories}
\label{sec2}

Before construction of the general representation theory we need to introduce some types of monoidal categories which is appropriate for definition of (co)representations. The condition on such monoidal categories is given by means of the notion of relative adjoints (with a parameter). For completeness we formulate and prove some needed statements on the relative adjoints. This requires in turn some facts following from Yoneda Lemma.

\subsection{Notations and conventions}
\label{secNot}

Fist of all, let us fix some conventions. In this work we use the terms and notations introduced in~\cite[\S~2]{Sqrt}. Let us remind some of them. We also recall some definitions from~\cite{Mcl} and introduce some terminology used below. 

\bb{Morphisms and functors.} In order to unify the notations with the previous publications we denote the composition of morphisms $f,g$ in a category $\bfC$ by $f\cdot g$ instead of $f\circ g$. For a fixed object $X\in\bfC$ denote by $\Hom(X,-)\colon\bfC\to\Set$ the functor from $\bfC$ to the category of sets, it maps an object $Y\in\bfC$ to the set $\Hom(X,Y)=\Hom_\bfC(X,Y)\in\Set$ which consists of all the morphisms in $\bfC$ of the form $X\to Y$, a morphism $f\colon Y\to Z$ is translated by this functor to the map $\Hom(X,f)\colon\Hom(X,Y)\to\Hom(X,Z)$, $g\mapsto f\cdot g\in\Hom(X,Z)$, $\;\forall\,g\in\Hom(X,Y)$. Analogously we define the contravariant functor $\Hom(-,X)\colon\bfC\to\Set$, $\Hom(f,X)\colon\Hom(Z,X)\to\Hom(Y,X)$, $h\mapsto h\cdot f\in\Hom(Y,X)$, $h\in\Hom(Z,X)$. Together these formulae give the bifunctor $\Hom=\Hom_\bfC=\Hom(-,-)\colon\bfC^\op\times\bfC\to\Set$.

Let $F\colon\bfC\to\bfD$ be a functor from a category $\bfC$ to a category $\bfD$. If $\bfC'\subset\bfC$ is a subcategory, then we can restrict $F$ to $\bfC'$ and get a new functor $\bfC'\to\bfD$, which we usually denote by the same letter $F$; formally, this is a composition $F\cdot E_{\bfC'}$, where $E_{\bfC'}\colon\bfC'\hookrightarrow\bfC$ is the categorical embedding. If $\bfD'$ is a subcategory of $\bfD$ such that $Ff$ is a morphism in $\bfD'$ for any morphism $f$ in $\bfC'$, then the restriction gives a functor $F\colon\bfC'\to\bfD'$.

\bb{Natural transformations.} \label{bbNatTr}
Let $\alpha\colon F\to G$ be a natural transformation between functors $F\colon\bfC\to\bfD$ and $G\colon\bfC\to\bfD$. We sometimes omit the subscript $X$ of the components $\alpha_X\colon FX\to GX$. For instance, for a morphism $f\colon Y\to Z$ in $\bfC$ the formulae $f_*(g)=f\cdot g$ and $f^*(h)=h\cdot f$ define the components $f_*=(f_*)_X=\Hom(X,f)\colon\Hom(X,Y)\to\Hom(X,Z)$ and $f^*=(f^*)_X=\Hom(f,X)\colon\Hom(Z,X)\to\Hom(Y,X)$ of the natural transformations $f_*\colon\Hom(-,Y)\to\Hom(-,Z)$ and $f^*\colon\Hom(Z,-)\to\Hom(Y,-)$ respectively.

For natural transformations $\alpha\colon F\to G$ and $\beta\colon G\to H$ between functors $F,G,H\colon\bfC\to\bfD$ we denote their composition by $\beta\cdot\alpha\colon F\to H$. This is a natural transformation with the components $(\beta\cdot\alpha)_X=\beta_X\cdot\alpha_X\colon FX\to HX$.

We say that the morphisms $\alpha=\alpha_X\colon FX\to GX$ are natural in $X$ iff they define a natural transformation $\alpha\colon F\to G$.
Recall that the natural transformation $\alpha\colon F\to G$ is an isomorphism in the category of functors $\bfC\to\bfD$ iff all the components $\alpha=\alpha_X\colon FX\to GX$ are isomorphisms in $\bfD$; the inverse natural transformation $\alpha^{-1}\colon G\to F$ has components $\alpha^{-1}=(\alpha_X)^{-1}\colon GX\to FX$ inverse to $\alpha_X$. Let $\bfC'\subset\bfC$ be a subcategory, then we say that $\alpha=\alpha_X\colon FX\to GX$ are natural in $X\in\bfC'$ iff they define a natural transformation between the restricted functors $F$ and $G$, i.e. $\alpha\colon F\cdot E_{\bfC'}\to G\cdot E_{\bfC'}$. In this case it is enough to suppose that the morphisms $\alpha_X$ are defined for $X\in\bfC'$.

Let $\bfP$ be also a category and $F$ and $G$ be bifunctors $\bfC\times\bfP\to\bfD$. Let $\bfC'\subset\bfC$ and $\bfP'\subset\bfP$ be subcategories. Then the requirement that $\alpha=\alpha_{X,Y}\colon F(X,Y)\to G(X,Y)$ are natural in $X\in\bfC'$ and $Y\in\bfP'$ means that these morphisms give a natural transformation between the restrictions of $F$ and $G$ to the subcategory $\bfC'\times\bfP'\subset\bfC\times\bfP$. For a fixed object $Y\in\bfP$ we have functors $F_Y,G_Y\colon\bfC\to\bfD$ such that $F_Y(f)=F(f,Y)$ and $G_Y(f)=G(f,Y)$ for any morphism $f$ in $\bfC$. These functors can be interpreted as functors with a parameter. We say that $\alpha=\alpha_{X,Y}\colon F_Y(X)\to G_Y(X)$ are natural in $X\in\bfC'$ iff for any $Y\in\bfP$ they give a natural transformation between the restricted functors $F_Y,G_Y\colon\bfC'\to\bfD$ (in this case we do not suppose that $F_Y$ and $G_Y$ are functorial in $Y$, i.e. morphisms $F(f,X)$ and $G(f,X)$ may be undefined). Analogous conventions are supposed for the functors with three or more arguments. 

Note that any functor $F\colon\bfC\to\bfD$ can be considered as the natural transformation between the functors $\Hom_\bfC$ and $\Hom_\bfD\big(F(-),F(-)\big)\colon\bfC^\op\times\bfC\to\Set$ with the components $F\colon\Hom(X,Y)\to\Hom(FX,FY)$. By fixing one of the arguments we obtain the natural transformation $\Hom(X,-)\to\Hom\big(FX,F(-)\big)$ or $\Hom(-,Y)\to\Hom\big(F(-),FY\big)$.

\subsection{Corollaries of Yoneda Lemma}
\label{secYL}

Yoneda Lemma and its corollaries are the main tool in the theory of adjoints and their relative version.

\bb{Yoneda Lemma.}
Let $R$ be an object of a category $\bfD$ and $K\colon\bfD\to\Set$ be a functor.
Yoneda Lemma establishes a one-to-one correspondence between the natural transformations from the functor $\Hom(R,-)\colon\bfD\to\Set$ to the functor $K$ and the elements of the set $KR$ (see~\cite[\S~3.2]{Mcl}). It can be formulated as follows.

\begin{Lem} \label{LemYL}
 Any natural transformation $\alpha\colon\Hom(R,-)\to K$ can be written in the form $\alpha_Z(f)=K(f)(u)\in KZ$, $f\in\Hom(R,Z)$, for a unique element $u\in KR$.
\end{Lem}

We are interested in some particular cases of this Lemma and their corollaries.

\bb{The case $K=\Hom\big(X,G(-)\big)$.} \label{bbLYG}
Let $G\colon\bfD\to\bfC$ be a functor and $X\in\bfC$. Consider the functor $K=\Hom\big(X,G(-)\big)$, it maps an object $Z\in\bfD$ to the set $\Hom(X,GZ)$ and a morphism $f\colon Z\to Z'$ to the map $(Gf)_*\colon\Hom(X,GZ)\to\Hom(X,GZ')$. In this case Yoneda Lemma~\ref{LemYL} takes the following form.

\begin{Cor} \label{CorYLG}
 Any maps $\alpha_Z\colon\Hom(R,Z)\to\Hom(X,GZ)$ natural in $Z\in\bfD$ have the form $\alpha(f)=(Gf)\cdot\gamma$ for a unique morphism $\gamma\colon X\to GR$.
\end{Cor}

Denote the natural transformation $\alpha\colon\Hom(R,-)\to\Hom\big(X,G(-)\big)$ given by a morphism $\gamma\colon X\to GR$ as in Corollary~\ref{CorYLG} by $\alpha=\gamma^*\cdot G$, this is a composition of natural transformations in the sense of p.~\ref{bbNatTr}, where the functor $G$ is regarded as the natural transformation $\Hom(R,-)\to\Hom\big(GR,G(-)\big)$ (this is not `horizontal' composition defined in~\cite[\S~2.5]{Mcl}). If $\bfC=\bfD$ and $G=\id_{\bfC}$, then $\alpha=\gamma^*$.

Consider a category $\bfC'$ and functors $F_1\colon\bfC'\to\bfD$, $F_2\colon\bfC'\to\bfC$. 

\begin{Cor} \label{CorYLGF1F2}
 The morphisms $\gamma_X\colon F_2X\to GF_1X$ are natural in $X\in\bfC'$ iff the morphisms $\alpha_{X,Z}=\gamma_X^*\cdot G\colon\Hom(F_1X,Z)\to\Hom(F_2X,GZ)$ are natural in $X\in\bfC'$.
\end{Cor}

\noindent{\bf Proof.}
If $\gamma_X\colon F_2X\to GF_1X$ are natural in $X$, then for any morphism $g\colon X\to Y$ in $\bfC'$ we have $\gamma_Y\cdot F_2g=GF_1g\cdot\gamma_X$, hence $(F_2g)^*\cdot\gamma_Y^*\cdot G=(\gamma_Y\cdot F_2g)^*\cdot G=(GF_1g\cdot\gamma_X)^*\cdot G=\gamma_X^*\cdot G\cdot(F_1g)^*$. This implies that the diagram
\begin{align} \label{F1F2XYZ}
\xymatrix{
 \Hom(F_1Y,Z)\ar[r]^G\ar[d]^{(F_1g)^*}&\Hom(GF_1Y,GZ)\ar[r]^{\gamma_Y^*}&\Hom(F_2Y,GZ)\ar[d]^{(F_2g)^*} \\
 \Hom(F_1X,Z)\ar[r]^G&\Hom(GF_1X,GZ)\ar[r]^{\gamma_X^*}&\Hom(F_2X,GZ)
}
\end{align}
is commutative, which means exactly that $\alpha_{X,Z}$ is natural in $X$. Conversely, let the diagram~\eqref{F1F2XYZ} be commutative. Denote $\delta_{Z}=\gamma_X^*\cdot G\cdot(F_1g)^*=\gamma_X^*\cdot(GF_1g)^*\cdot G$, this is the diagonal arrow $\Hom(F_1X,Z)\to\Hom(F_2Y,GZ)$ in the diagram~\eqref{F1F2XYZ}. Note that $\delta_{Z}$ are natural in $Z$, so by virtue of Corollary~\ref{CorYLG} there exist unique morphisms $f_{X,Y}\colon F_2Y\to F_1X$ such that $\delta_{Z}=f_{X,Y}^*\cdot G$. Hence $\gamma_Y\cdot F_2g=f_{X,Y}=GF_1g\cdot\gamma_X$. \qed 

\bb{Universal morphisms.}
Recall that {\it universal morphism from an object $X\in\bfC$ to a functor $G\colon\bfD\to\bfC$} is a pair $(R,\eta)$, where $R\in\bfD$ and $\eta\colon X\to GR$ are an object and a morphism such that for any $Z\in\bfD$ and $f\colon X\to GZ$ there is a unique $h\colon R\to Z$ such that $f=Gh\cdot\eta$.

Analogously, {\it universal morphism from a functor $F\colon\bfC\to\bfD$ to an object $Z\in\bfD$} is a pair $(P,\varepsilon)$, where $P\in\bfC$ and $\varepsilon\colon FP\to Z$ are such that for any $X\in\bfC$ and $h\colon FX\to Z$ there is a unique $f\colon X\to P$ such that $h=\varepsilon\cdot Ff$. This notion can be obtained from the previous one by the categorical duality.

In terms of p.~\ref{bbLYG} the condition that $(R,\eta)$ is a universal morphism from an object $X\in\bfC$ to a functor $G\colon\bfD\to\bfC$ is equivalent to the bijectivity of the natural transformation $\eta^*\cdot G\colon\Hom(R,-)\to\Hom\big(X,G(-)\big)$. Due to Corollary~\ref{CorYLG} we can reformulate it as follows (see~\cite[\S~3.2, Prop.~1]{Mcl}).

\begin{Cor} \label{CorPropYL}
 Any bijections $\Hom(R,Z)\isoright\Hom\big(X,GZ\big)$ natural in $Z\in\bfD$ have the form $\eta^*\cdot G$ for a unique universal morphism $(R,\eta)$ from $X\in\bfC$ to $G\colon\bfD\to\bfC$.
\end{Cor}

\begin{Cor} \label{CorYLUM}
Consider functors $G\colon\bfD\to\bfC$, $F_1\colon\bfC'\to\bfD$ and $F_2\colon\bfC'\to\bfC$. Any bijections $\beta_{X,Z}\colon\Hom(F_1X,Z)\isoright\Hom\big(F_2X,GZ\big)$ natural in $X\in\bfC'$ and $Z\in\bfD$ have the form $\eta_X^*\cdot G$ for a unique natural transformation $\eta\colon F_2\to GF_1$ such that $(F_1X,\eta_X)$ is a universal morphism from $F_2X\in\bfC'$ to $G\colon\bfD\to\bfC$.
\end{Cor}

\noindent{\bf Proof.} This statement follows directly from Corollaries~\ref{CorPropYL} and~\ref{CorYLGF1F2}. \qed

\bb{The case $K=\Hom(X,-)$.}
If the functor $K$ has the form $K=\Hom(X,-)$ for some $X\in\bfD$, then Yoneda Lemma~\ref{LemYL} is a particular case of Corollary~\ref{CorYLG} for $\bfC=\bfD$ and $G=\id_{\bfD}$. It can be written in the following form (see~\cite[\S~3.2]{Mcl}).

\begin{Cor} \label{CorYL}
 Any natural transformation $\Hom(R,-)\to\Hom(X,-)$ has the form $\gamma^*$ for a unique morphism $\gamma\colon X\to R$.
\end{Cor}

Note that $(R,\gamma)$ is a universal morphism form $R$ to $\id_\bfD$ iff $\gamma$ is an isomorphism. Hence Corollary~\ref{CorPropYL} takes the following from.

\begin{Cor} \label{CorYL1}
 The natural transformation $\gamma^*\colon\Hom(R,-)\to\Hom(X,-)$ is an isomorphism of functors iff $\gamma\colon X\to R$ is an isomorphism in $\bfD$.
\end{Cor}

\begin{Rem} \normalfont
 \label{RemUnit}
The universal morphism is unique up to an isomorphism: if $(R,\eta)$ and $(R',\eta')$ are two universal morphisms from $X$ to $G$, then there is a unique isomorphism $\gamma\colon R\isoright R'$ such that $G\gamma\cdot\eta=\eta'$. This fact is deduced in~\cite[\S~3.1]{Mcl} by using comma categories, however it also follows from (the corollaries of) Yoneda Lemma. To show it one needs to apply Corollaries~\ref{CorYL} and \ref{CorYL1} to the composition of $(\eta')^*\cdot G\colon\Hom(R',-)\isoright\Hom\big(X,G(-)\big)$ with the inverse of $\eta^*\cdot G\colon\Hom(R,-)\isoright\Hom\big(X,G(-)\big)$ and then to take into account Corollary~\ref{CorYLG}. Dually, a universal morphism $(P,\varepsilon)$ from $F\colon\bfC\to\bfD$ to an object $Z\in\bfD$ is unique up to a unique isomorphism.
\end{Rem}

Let $F_1$ and $F_2$ be functors $\bfC'\to\bfD$. Corollaries~\ref{CorYLG}, \ref{CorYLGF1F2} and \ref{CorYLUM} imply the following statements.

\begin{Cor} \label{CorYL2}
 If maps $\beta_{X,Z}\colon\Hom(F_1X,Z)\to\Hom(F_2X,Z)$ are natural in $X\in\bfC'$ and $Z\in\bfD$, then $\beta_{X,Z}=\gamma_X^*$ for a unique natural transformation $\gamma\colon F_2\to F_1$.
\end{Cor}

\begin{Cor} \label{CorYL3}
Any bijections $\beta_{X,Z}\colon\Hom(F_1X,Z)\isoright\Hom(F_2X,Z)$ natural in $X\in\bfC'$ and $Z\in\bfD$ have the form $\beta_{X,Z}=\gamma_X^*$ for a unique natural isomorphism $\gamma\colon F_2\isoright F_1$.
\end{Cor}

\subsection{Relative right and left adjoint functors}
\label{secRelAdj}

Now we introduce a notion of relative right/left adjunctions and its version `with a parameter' which are necessary for an appropriate generalisation of the (co)closed monoidal categories.

\bb{Relative adjunctions.}
Let $\bfC$ and $\bfD$ be categories. Recall that a functor $F\colon\bfC\to\bfD$ is called {\it left adjoint} for $G\colon\bfD\to\bfC$ and $G$ is called {\it right adjoint} for $F$ when there exist bijections
\begin{align} \label{alphaXZ}
 \alpha_{X,Z}\colon\Hom_\bfD(FX,Z)\isoright\Hom_\bfC(X,GZ)
\end{align}
natural in $X\in\bfC$ and $Z\in\bfD$. It gives a natural isomorphism $\alpha$ between the bifunctors defined as compositions 
\begin{align}
 &\Hom_\bfD\big(F(-),-\big)\colon \bfC^\op\times\bfD\xrightarrow{F^\op\times\id_\bfD}\bfD^\op\times\bfD\xrightarrow{\Hom_\bfD}\Set, \label{HomF} \\
 &\Hom_\bfC\big(-,G(-)\big)\colon \bfC^\op\times\bfD\xrightarrow{\id_{\bfC^\op}\times G} \bfC^\op\times\bfC\xrightarrow{\Hom_\bfC}\Set, \label{HomG}
\end{align}
where $F^\op\colon\bfC^\op\to\bfD^\op$ is the functor opposite to $F$. The isomorphism $\alpha$ is called {\it adjunction bijection} for the adjoints $F,G$. The triple $(F,G,\alpha)$ is called {\it adjunction}~\cite[\S~4.1]{Mcl}. 

Consider a full subcategory $\bfD'\subset\bfD$. If we restrict the functors~\eqref{HomF}, \eqref{HomG} to the subcategory $\bfC^\op\times\bfD'\subset\bfC^\op\times\bfD$, then we obtain isomorphisms~\eqref{alphaXZ} natural in $X\in\bfC$ and $Z\in\bfD'$. In this case it is sufficient to define the functor $G$ on the subcategory $\bfD'$, however the functors $F\colon\bfC\to\bfD$ and $G\colon\bfD'\to\bfC$ are not adjoint to each other, because the domain of such $G$ does not coincide with the codomain of $F$ (we can not reduce this situation to the adjunction case unless $F(\bfC)\subset\bfD'$, when one can consider $F$ as a functor $\bfC\to\bfD'$).

A generalisation of the notion of adjoint functor appropriate for this situation was introduced by Ulmer in~\cite[Def.~2.3]{Ulm}. Namely, let $\bfB,\bfC,\bfD$ be categories and $J\colon\bfB\to\bfD$, $F\colon\bfC\to\bfD$, $G\colon\bfB\to\bfC$ be functors. A functor $G$ is called {\it right adjoint for $F$ relative to $J$} iff $\Hom_\bfD\big(F(-),J(-)\big)$ and $\Hom_\bfD\big(-,G(-)\big)$ are isomorphic as functors $\bfC^\op\times\bfB\to\Set$. Similarly relative left adjoint is defined~\cite[Def.~2.2]{Ulm}. In this work we consider the case of fully faithful functor $J$ only. For this case we can interpret $\bfB$ as a subcategory of $\bfD$ and $J$ as the corresponding category embedding (up to category equivalence the notion of subcategory essentially coincides with the notion of fully faithful functor). Let us formulate the definition of relative right/left adjoint functor for this particular case in details.

\begin{Def} \normalfont
 Let $F\colon\bfC\to\bfD$ be a functor. A functor $G\colon\bfD'\to\bfC$ is called {\it right adjoint for the functor $F$ relative to} the full subcategory $\bfD'\subset\bfD$ iff there exist isomorphisms of the form~\eqref{alphaXZ} natural in $X\in\bfC$ and $Z\in\bfD'$. These isomorphisms give a natural isomorphism $\alpha$ between the functors
\begin{align}
 & \bfC^\op\times\bfD'\xrightarrow{F^\op\times E_{\bfD'}} \bfD^\op\times\bfD\xrightarrow{\Hom_\bfD}\Set, \label{HomFr} \\
 & \bfC^\op\times\bfD'\xrightarrow{\id_{\bfC^\op}\times G} \bfC^\op\times\bfC\xrightarrow{\Hom_\bfC}\Set; \label{HomGr}
\end{align}
then $(F,G,\alpha)$ is called {\it relative right adjunction} and $\alpha$ is called {\it right adjunction bijection}. In this case we say that $F$ {\it has a right adjoint relative to $\bfD'$}.

Let $G\colon\bfD\to\bfC$ be a functor. A functor $F\colon\bfC'\to\bfD$ is called {\it left adjoint for the functor $G$ relative to} the full subcategory $\bfC'\subset\bfC$ iff there exist isomorphisms of the form~\eqref{alphaXZ} natural in $X\in\bfC'$ and $Z\in\bfD$. These isomorphisms give a natural isomorphism $\alpha$ between the functors
\begin{align}
 & (\bfC')^\op\times\bfD\xrightarrow{F^\op\times\id_\bfD}\bfD^\op\times\bfD\xrightarrow{\Hom_\bfD}\Set, \label{HomFl} \\
 & (\bfC')^\op\times\bfD\xrightarrow{E_{(\bfC')^\op}\times G} \bfC^\op\times\bfC\xrightarrow{\Hom_\bfC}\Set; \label{HomGl}
\end{align}
then $(F,G,\alpha)$ is called {\it relative left adjunction} and $\alpha$ is called {\it left adjunction bijection}. In this case we say that $G$ {\it has a left adjoint relative to $\bfC'$}.
\end{Def}

Note that $G\colon\bfD'\to\bfC$ is a right adjoint for $F\colon\bfC\to\bfD$ relative to $\bfD'\subset\bfD$ iff the opposite functor $G^\op\colon(\bfD')^\op\to\bfC^\op$ is a left adjoint for $F^\op\colon\bfC^\op\to\bfD^\op$ relative to $(\bfD')^\op\subset\bfD^\op$.

\begin{War} \normalfont
In contrast to the case of usual adjunction the fact the $F$ is a relative left adjoint for $G$ does not mean that $G$ is the relative right adjoint for $F$ and vice versa.
\end{War}

Let $\bfD''$ be a full subcategory of $\bfD'$. We see from the definition that if $G\colon\bfD'\to\bfC$ is a right adjoint for $F\colon\bfC\to\bfD$ relative to $\bfD'$, then the restriction of $G$ to $\bfD''$ is a right adjoint for $F$ relative to $\bfD''$. In particular, if $F$ has a usual right adjoint (i.e. relative to the whole $\bfD$), then it has a right adjoint relative to any subcategory $\bfD'\subset\bfD$. The same we can say about relative left adjoint functors.

Sometimes a right/left adjoint for a fixed functor does not exists, but there exists a right/left adjoint relative to some full subcategory. Below we show how to determine the maximal full subcategory such that a right/left adjoint exists relative to this subcategory.

A typical example is the forgetful functor $G\colon\FVect\to\Set$ from the category of finite-dimensional vector spaces to the category of sets: the usual left adjoint does not exists, but we have a left adjoint relative to the subcategory of the finite sets. 

In applications we sometimes have situations, when we are not interested if a usual right/left adjoint exists, we only need to establish the existence of the right/left adjoint relative to a fixed full subcategory.

\bb{Units and counits for relative adjunctions.}
Let us generalise some properties and `equivalent definitions' of the adjoint functors (\cite[\S~4.1, Theorems 1,~2]{Mcl}) for the case of relative adjoints. In some form they were formulated in~\cite{Ulm} without proof. For completeness we briefly add the proofs by following ideas of~\cite{Mcl} and by applying the corollaries of Yoneda Lemma.

\begin{Th} \label{ThUnit}
 Consider a functor $F\colon\bfC\to\bfD$ and a full subcategory $\bfD'\subset\bfD$.
\begin{itemize}
 \item[$(1).$] Let $G\colon\bfD'\to\bfC$ be a right adjoint for $F$ relative to $\bfD'$. The right adjunction bijections $\alpha_{X,Z}\colon\Hom(FX,Z)\isoright\Hom(X,GZ)$ are inverse to
\begin{align} \label{alphaeps}
 &\alpha^{-1}(f)=\varepsilon_Z\cdot Ff, &&&&f\colon X\to GZ, &&X\in\bfC, &&Z\in\bfD',
\end{align}
where $\varepsilon_Z:=\alpha^{-1}(\id_{GZ})\colon FGZ\to Z$ are natural in $Z\in\bfD'$. The pair $(GZ,\varepsilon_Z)$ is a universal morphism from $F$ to $Z\in\bfD'$.
 \item[$(2).$] Suppose that for any $Z\in\bfD'$ there exists a universal morphism $(X_Z,\varepsilon_Z)$ from $F$ to $Z$. Then $F$ has a right adjoint relative to $\bfD'$. This relative right adjoint can be defined as the unique functor $G\colon\bfD'\to\bfC$ such that $GZ=X_Z$ $\;\forall\,Z\in\bfD'$ and $\varepsilon_Z$ form a natural transformation from $FG$ to $E_{\bfD'}$. The corresponding right adjunction bijections can be written via its inverse by the formula~\eqref{alphaeps}.
\end{itemize}
 Dually, consider a functor $G\colon\bfD\to\bfC$ and a full subcategory $\bfC'\subset\bfC$.
\begin{itemize}
 \item[$(1').$] Let $F\colon\bfC'\to\bfD$ be a left adjoint for $G$ relative to $\bfC'\subset\bfC$. The left adjunction bijections $\alpha_{X,Z}\colon\Hom(FX,Z)\isoright\Hom(X,GZ)$ can be written as
\begin{align} \label{alphaeta}
 &\alpha(h)=Gh\cdot\eta_X, &&&&h\colon FX\to Z, &&X\in\bfC', &&Z\in\bfD,
\end{align}
where $\eta_X:=\alpha(\id_{FX})\colon X\to GFX$ are natural in $X\in\bfC'$, and $(FX,\eta_X)$ is a universal morphism from $X\in\bfC'$ to $G$.
 \item[$(2').$] Suppose that for any $X\in\bfC'$ there exists a universal morphism $(Z_X,\eta_X)$ from $X$ to $G$. Then $G$ has a left adjoint relative to $\bfC'$. This relative left adjoint can be defined as the unique functor $F\colon\bfC'\to\bfD$ such that $FX=Z_X$ $\;\forall\,X\in\bfC'$ and $\eta_X$ form a natural transformation from $E_{\bfC'}$ to $GF$. The corresponding left adjunction bijections have the form~\eqref{alphaeta}.
\end{itemize}
\end{Th}

\noindent{\bf Proof.} The statements $(1)$ and $(2)$ are equivalent to $(1')$ and $(2')$ via the categorical duality, so it is sufficient to prove the latter two.

$(1')$: This is a particular case of Corollary~\ref{CorYLUM} for $F_1=F$ and $F_2=E_{\bfC'}$.

$(2')$: The functor $F\colon\bfC'\to\bfD$ is constructed as follows. Define it on objects $X\in\bfC'$ as $FX=Z_X$. For a morphism $g\colon X\to X'$ in $\bfC'$ we define $Fg$ as the unique morphism $FX\to FX'$ such that the diagram
\begin{align} \label{etaNat}
\xymatrix{
 X\ar[r]^{\eta_X}\ar[d]^{g} & GFX\ar[d]^{GFg} \\
 X'\ar[r]^{\eta_{X'}} & GFX'
}
\end{align}
is commutative. By construction the functor $F$ is a unique functor $\bfC'\to\bfD$ such that $FX=Z_X$ $\;\forall\,X\in\bfC'$ and the morphisms $\eta_X\colon X\to GFX$ are natural in $X\in\bfC'$. Define the maps $\alpha_{X,Z}\colon\Hom(FX,Z)\to\Hom(X,GZ)$ by the formula~\eqref{alphaeta}, which implies that $\alpha_{X,Z}$ are natural in $Z\in\bfD$. Their naturality in $X\in\bfC'$ follows from the naturality of $\eta$. The bijectivity of the map $\alpha_{X,Z}$ follows from the universality of $(FX,\eta_X)$, so $F$ is a left adjoint for $G$ relative to $\bfC'$. \qed

The natural transformation $\eta\colon E_{\bfC'}\to GF$ with the components $\eta_X:=\alpha_{X,FX}(\id_{FX})$ is called {\it unit} for the relative left adjunction $(F,G,\alpha)$. Dually, the natural transformation $\varepsilon\colon FG\to E_{\bfD'}$ defined in Theorem~\ref{ThUnit} is called {\it counit} for the relative right adjunction $(F,G,\alpha)$. Note that in general we do not have a unit for a relative right adjunction, nor a counit for a relative left adjunction.

We see from Theorem~\ref{ThUnit} that the maximal full subcategory $\bfD'\subset\bfD$ such that the functor $F\colon\bfC\to\bfD$ has a right adjoint relative to $\bfD'$ can be found as follows. It consists of all objects $Z\in\bfD$ such that there is a universal morphism from $F$ to $Z$. Analogously, the existence of a universal morphism from $X$ to $G\colon\bfD\to\bfC$ defines the maximal full subcategory $\bfC'\subset\bfC$ such that $G$ has a left adjoint relative to $\bfC'$.

\bb{Relative adjunctions with a parameter.}
For a bifunctor $F\colon\bfC\times\bfP\to\bfD$ we can consider a family of functors $F_Y=F(-,Y)\colon\bfC\to\bfD$ parametrised by objects of the category $\bfP$. If $F_Y$ has a right adjoint $G_Y$ for any $Y\in\bfP$, then there is a unique functor $G\colon\bfP^\op\times\bfD\to\bfC$ such that $G(Y,f)=G_Y(f)$ for any morphism $f$ in $\bfD$ and the corresponding adjunction bijections
\begin{align} \label{alphaXYZ}
 \alpha_{X,Y,Z}\colon\Hom_\bfD\big(F(X,Y),Z\big)\isoright\Hom_\bfC\big(X,G(Y,Z)\big)
\end{align}
are natural in $X\in\bfC$, $Y\in\bfP$, $Z\in\bfD$ (see~\cite[\S~4.7, Th.~3]{Mcl}). This fact is immediately generalised to the case of relative adjunctions.

\begin{Th} \label{ThParam}
 Let $F\colon\bfC\times\bfP\to\bfD$ be a bifunctor such that the functor $F(-,Y)\colon\bfC\to\bfD$ has a right adjoint relative to a full subcategory $\bfD'\subset\bfD$ for each $Y\in\bfP$, so we have a relative right adjunction $\big(F(-,Y),G_Y,\alpha_Y\big)$ for each $Y\in\bfP$, where $G_Y\colon\bfD'\to\bfC$ is the relative right adjoint and $\alpha_Y$ is the corresponding right adjunction bijection. Then there exists a unique functor $G\colon\bfP^\op\times\bfD'\to\bfC$ such that $G(Y,-)=G_Y$ and the components $\alpha_{X,Y,Z}=(\alpha_Y)_{X,Z}$ give the isomorphism of the form~\eqref{alphaXYZ} natural in $X\in\bfC$, $Y\in\bfP$, $Z\in\bfD'$.

 Let $G\colon\bfD\times\bfP\to\bfC$ be a bifunctor such that each functor $G(-,Y)\colon\bfD\to\bfC$ has a left adjoint relative to a full subcategory $\bfC'\subset\bfC$, so we have a relative left adjunction $\big(F_Y,G(-,Y),\alpha_Y\big)$ for each $Y\in\bfP$, where $F_Y\colon\bfC'\to\bfD$. Then there exists a unique functor $F\colon\bfP^\op\times\bfC'\to\bfD$ such that $F(Y,-)=F_Y$ and the components $\alpha_{X,Y,Z}=(\alpha_Y)_{X,Z}$ give the isomorphism of the form
\begin{align} \label{alphaXYZl}
 \alpha_{X,Y,Z}\colon\Hom_\bfD\big(F(Y,X),Z\big)\isoright\Hom_\bfC\big(X,G(Z,Y)\big)
\end{align}
natural in $X\in\bfC'$, $Y\in\bfP$, $Z\in\bfD$.
\end{Th}

\noindent{\bf Proof.} Let us prove the second part (the first part is dual). By having the functors $F_Y=F(Y,-)\colon\bfC'\to\bfD$ for $Y\in\bfP$ we need to construct a bifunctor $F\colon\bfP^\op\times\bfC'\to\bfD$ such that~\eqref{alphaXYZl} is natural in $X\in\bfC'$, $Y\in\bfP$ and $Z\in\bfD$. This amounts to construction of functors $\wt F_X=F(-,X)\colon\bfP^\op\to\bfD$ satisfying the following conditions: $F_Y(X)=\wt F_X(Y)$, the bijection~\eqref{alphaXYZl} is natural in $Y\in\bfP$ (for each fixed $X\in\bfC'$, $Z\in\bfD$) and the diagram
\begin{align} \label{FFfg}
\xymatrix{
 F(Y',X)\ar[rr]^{F(Y',f)}\ar[d]^{F(g,X)}& & F(Y',X')\ar[d]^{F(g,X')} \\
 F(Y,X)\ar[rr]^{F(Y,f)}& & F(Y,X') \\
}
\end{align}
commutes for any morphisms $f\colon X\to X'$ in $\bfC'$ and $g\colon Y\to Y'$ in $\bfP$ (see~\cite[\S~2.3, Prop.~1,~2]{Mcl}). For each such $g\colon Y\to Y'$ consider the morphism $\beta_{X,Z}(g)$ defined by the commutative diagram
\begin{align} \label{alphaY}
\xymatrix{
\Hom\big(F(Y,X),Z\big)\ar[rr]^{\alpha_{X,Y,Z}}\ar[d]^{\beta_{X,Z}(g)} &&\Hom\big(X,G(Z,Y)\big)\ar[d]^{G(Z,g)_*} \\
\Hom\big(F(Y',X),Z\big)\ar[rr]^{\alpha_{X,Y',Z}} && \Hom\big(X,G(Z,Y')\big)
}
\end{align}
It is natural in $X\in\bfC'$ and $Z\in\bfD$, since $\alpha_{X,Y,Z}$ is so. Hence by virtue of Corollary~\ref{CorYL2} there exist unique morphisms $F(g,X)\colon F(Y',X)\to F(Y,X)$ natural in $X\in\bfC'$ such that $\beta_{X,Z}(g)=F(g,X)^*$ for any $Z\in\bfD$. Note that $\beta_{X,Z}(g_1\cdot g_2)=\beta_{X,Z}(g_1)\cdot\beta_{X,Z}(g_2)$ for any morphisms $g_1$ and $g_2$ in $\bfP$, so the uniqueness of $F(g,X)$ implies that $\wt F_X=F(-,X)$ is a functor $\bfP^\op\to\bfD$ for each $X\in\bfC'$. We see from the diagram~\eqref{alphaY} that it makes $\alpha_{X,Y,Z}$ to be natural in $Y\in\bfP$. The uniqueness of the morphisms implies also that such functor is unique. The naturality of $F(g,X)$ in $X\in\bfC'$ means exactly the commutativity of the diagram~\eqref{FFfg}. \qed

\begin{Def} \normalfont
 A {\it right adjunction with a parameter} for a bifunctor $F\colon\bfC\times\bfP\to\bfD$ {\it relative} to $\bfD'\subset\bfD$ is a triple $(F,G,\alpha)$, where $G\colon\bfP^\op\times\bfD'\to\bfC$ is a functor and $\alpha$ is an isomorphism of the form~\eqref{alphaXYZ} natural in $X\in\bfC$, $Y\in\bfP$, $Z\in\bfD'$.
 A {\it left adjunction with a parameter} for a bifunctor $G\colon\bfD\times\bfP\to\bfC$ {\it relative} to $\bfC'\subset\bfC$ is a triple $(F,G,\alpha)$, where $F\colon\bfP^\op\times\bfC'\to\bfD$ is a functor and $\alpha$ is an isomorphism of the form~\eqref{alphaXYZl} natural in $X\in\bfC'$, $Y\in\bfP$, $Z\in\bfD$.
\end{Def}

Theorem~\ref{ThParam} implies the following fact: a right/left adjunction with a parameter for a bifunctor $F\colon\bfC\times\bfP\to\bfD$ relative to $\bfD'\subset\bfD$ exists iff the functor $F(-,Y)$ has a right/left adjoint relative to $\bfD'$ for each object $Y\in\bfP$.

The counits of the relative right adjunctions $\big(F(-,Y),G(Y,-),\alpha\big)$ have components
\begin{align} \label{varepsilonYZ}
 \varepsilon_{Y,Z}=\alpha^{-1}(\id_{G(Y,Z)})\colon F\big(G(Y,Z),Y\big)\to Z.
\end{align}
These are morphisms in $\bfD$ natural in $Z\in\bfD'$. They are dinatural in $Y\in\bfP$ in the sense of~\cite[\S~9.4]{Mcl}.

The units of the relative left adjunctions $\big(F(Y,-),G(-,Y),\alpha\big)$ have components
\begin{align} \label{etaYZ}
 \eta_{Y,X}=\alpha(\id_{F(Y,X)})\colon X\to G\big(F(Y,X),Y\big).
\end{align}
These are morphisms in $\bfC$ natural in $X\in\bfC'$ and dinatural in $Y\in\bfP$.

Note that relative right/left adjunction is a particular case of the relative right/left adjunction with a parameter. This is the case of the category $\bfP$ which has a unique object and a unique arrow.

\bb{Uniqueness of relative adjunctions.}
If a usual right/left adjoint to a fixed functor exists, then it is unique up to an isomorphism~\cite[\S~4.1, Corollary~1]{Mcl}. Moreover, this isomorphism is unique and it relates the corresponding adjunction bijections. The same still holds for the case of relative adjunctions. We prove this for more general case: for the relative adjunctions with a parameter.

\begin{Th} \label{ThUniq}
 If there exist two right adjunctions with a parameter $(F,G,\alpha)$ and $(F,G',\alpha')$ for $F\colon\bfC\times\bfP\to\bfD$ relative to a full subcategory $\bfD'\subset\bfD$, then there exists a unique natural isomorphism $\gamma\colon G\isoright G'$ such that $(\gamma_{Y,Z})_*\cdot\alpha_{X,Y,Z}=\alpha'_{X,Y,Z}$ for all $X\in\bfC$, $Y\in\bfP$, $Z\in\bfD'$.

 Dually, if there exist two left adjunctions with a parameter $(F,G,\alpha)$ and $(F',G,\alpha')$ for $G\colon\bfD\times\bfP\to\bfC$ relative to a full subcategory $\bfC'\subset\bfC$, then there exists a unique natural isomorphism $\gamma\colon F\isoright F'$ such that $\alpha_{X,Y,Z}\cdot(\gamma_{Y,X})^*=\alpha'_{X,Y,Z}$ for all $X\in\bfC'$, $Y\in\bfP$, $Z\in\bfD$.
\end{Th}

\noindent{\bf Proof.} As above we prove the `left' version.
Consider the compositions
\begin{align}
 \Hom\big(F'(Y,X),Z\big)\xrightarrow{\alpha'_{X,Y,Z}}\Hom\big(X,G(Z,Y)\big)\xrightarrow{\alpha_{X,Y,Z}^{-1}}\Hom\big(F(Y,X),Z\big).
\end{align}
These are bijections $\beta_{Y,X,Z}\colon\Hom\big(F'(Y,X),Z\big)\isoright\Hom\big(F(Y,X),Z\big)$
natural in $X\in\bfC'$, $Y\in\bfP$ and $Z\in\bfD$. By virtue of Corollary~\ref{CorYL3} they have the form $\beta_{Y,X,Z}=\gamma_{Y,X}^*$ for a unique natural isomorphism $\gamma\colon F\isoright F'$.
 \qed

\begin{Rem} \normalfont
 The uniqueness property of relative adjunctions (with a parameter) can be rewritten in terms of counits~\eqref{varepsilonYZ} and units~\eqref{etaYZ} respectively. Let $\eta_{Y,X}=\alpha(\id_{F(Y,X)})$ and $\eta_{Y,X}=\alpha'(\id_{F'(Y,X)})$ be units for the relative left adjunctions with a parameter $(F,G,\alpha)$ and $(F',G,\alpha')$ respectively. The natural isomorphism $\gamma\colon F\isoright F'$ from Theorem~\ref{ThUniq} relates these units by the formula $\eta'_{Y,X}=G(\gamma_{Y,X},\id_Y)\cdot\eta_{Y,X}$. The existence and uniqueness of such isomorphisms $\gamma_{Y,X}\colon F(Y,X)\isoright F'(Y,X)$ follows from the fact that a universal morphism from $X$ to $G(-,Y)$ is unique up to a unique isomorphism (see Remark~\ref{RemUnit}). Analogously, the counits $\varepsilon_{Y,Z}=\alpha^{-1}(\id_{G(Y,Z)})$ and $\varepsilon'_{Y,Z}=(\alpha')^{-1}(\id_{G'(Y,Z)})$ are related as $\varepsilon'_{Y,Z}\cdot F(\gamma_{Z,Y},\id_Y)= \varepsilon_{Y,Z}$ for unique $\gamma_{Z,Y}$.
\end{Rem}

\subsection{Relatively closed monoidal categories}
\label{secRelCl}

Representations can be defined in a closed monoidal category. But sometimes we need to consider more general situation. The notion of adjunction (with a parameter) relative to a full subcategory allows us to generalise closed monoidal category to the case sufficient for our purposes.

\bb{Generalisations of closed and coclosed monoidal categories.}
Let $\bfC=(\bfC,\otimes)$ be a symmetric monoidal category with a monoidal product $-\otimes-\colon\bfC\times\bfC\to\bfC$.
As in~\cite{Sqrt} we suppose that the monoidal categories are strict. We denote the unit object of $\bfC$ by $I_\bfC$. Recall that $\bfC$ is called {\it closed monoidal category} iff the functor $-\otimes V\colon\bfC\to\bfC$ has a right adjoint for each $V\in\bfC$. This right adjoint is denoted as $h\mapsto\bfhom(V,h)$ and gives a bifunctor $\bfhom\colon\bfC^\op\times\bfC\to\bfC$ called {\it internal hom-functor}. Since $\bfC$ is symmetric, the right adjoint functors for $-\otimes V$ and $V\otimes-$ coincide with each other. Dually, $\bfC$ is called {\it coclosed monoidal category} iff $-\otimes V$ has a left adjoint for each $V\in\bfC$. In this case we obtain {\it internal cohom-functor} $\bfcohom\colon\bfC^\op\times\bfC\to\bfC$.

For purposes of Quantum Representation Theory we considered more general situation in~\cite{Sqrt}, when the functors $-\otimes V$ have right/left adjoints for the objects $V$ of a full subcategory $\bfP\subset\bfC$  (we do not always suppose that this subcategory is monoidal). For this situation we introduced representations/corepresentations of a monoid/comonoid in $\bfC$ on an object $V\in\bfP$. However the condition for these functors to have adjoints is too strong. Here we consider a weaker condition by requiring the existence of right/left adjoints relative to $\bfP\subset\bfC$.

\begin{Def} \normalfont
 The symmetric monoidal category $\bfC=(\bfC,\otimes)$ is called {\it closed/coclosed relative to} a full subcategory $\bfP\subset\bfC$ or {\it relatively closed/coclosed with a parametrising subcategory} $\bfP\subset\bfC$ iff the functor $-\otimes V\colon\bfC\to\bfC$ has a right/left adjoint relative to $\bfP\subset\bfC$ for each object $V\in\bfP$.
\end{Def}

The following facts are proved by restricting the corresponding relative adjoint functors.

\begin{Prop} \label{PropRelclosed}
 Let $\bfC=(\bfC,\otimes)$ be a symmetric monoidal category and $\bfP\subset\bfC$ be its full subcategory.
\begin{itemize}
 \item[{\normalfont(1).}] If the monoidal category $\bfC$ is closed/coclosed relative to $\bfP$, then it is closed/coclosed relative to any full subcategory $\bfP'\subset\bfP$.
 \item[{\normalfont(2).}] Let $\bfC'\subset\bfC$ be a full subcategory such that $\bfP\subset\bfC'$. If the functor $-\otimes V\colon\bfC\to\bfC$ has a right/left adjoint relative to $\bfC'$ for each $V\in\bfP$, then $(\bfC,\otimes)$ is closed/coclosed relative to $\bfP$.
 \item[{\normalfont(3).}] In particular, if the functor $-\otimes V\colon\bfC\to\bfC$ has a right/left adjoint for each $V\in\bfP$, then the monoidal category $\bfC$ is closed/coclosed relative to $\bfP$.
\end{itemize}
\end{Prop}

The condition supposed in the point~(3) of Prop.~\ref{PropRelclosed} is often fulfilled, however this stronger condition is not needed to construct the representation theory. The subcategory $\bfP\subset\bfC$ is usually chosen as bigger as possible such that the functor $-\otimes V$ has a right/left (relative) adjoint for all $V\in\bfP$. In some cases it is easier to check that these functors have adjoints relative to $\bfP$ than relative to the whole $\bfC$.

The simplest non-trivial example is the category of vector spaces $\bfC=\Vect$ with the standard tensor product $\otimes$, where the parametrising subcategory is $\bfP=\FVect$. The functor $-\otimes V\colon\Vect\to\Vect$ has a left adjoint for any $V\in\FVect$ (see~\cite[\S~3.1.11]{Sqrt}), so due to the point~(3) of Prop.~\ref{PropRelclosed} the symmetric monoidal category $(\Vect,\otimes)$ is coclosed relative to $\FVect$.

\bb{Internal hom.}
Due to Theorem~\ref{ThParam} (for the case $\bfD=\bfC$) a symmetric monoidal category $\bfC=(\bfC,\otimes)$ is relatively closed with the parametrising subcategory $\bfP\subset\bfC$ iff there exists a bifunctor $\bfhom\colon\bfP^\op\times\bfP\to\bfC$ with the bijections
\begin{align} \label{homDef}
 \theta=\theta_{X,V,W}\colon\Hom\big(X,\bfhom(V,W)\big)\isoright\Hom(X\otimes V,W)
\end{align}
natural in $X\in\bfC$, $V,W\in\bfP$. Or, equivalently, iff there exists a right adjunction with a parameter $(-\otimes-,\bfhom,\theta^{-1})$ relative to the full subcategory $\bfD'=\bfP$. The bifunctor $\bfhom$ is called {\it (relative) internal hom-functor}. The values $\bfhom(V,W)$ on the objects $V,W\in\bfP$ are called {\it internal hom-objects}.

One can consider the bifunctor $\bfhom$ with the natural isomorphism $\theta$ as an additional structure on the symmetric monoidal category $\bfC$. However,
Theorem~\ref{ThUniq} implies that the relative right adjunction with a parameter $\big(-\otimes-,\bfhom,\theta^{-1}\big)$ is unique up to a unique natural isomorphism. This means that the choice of the certain pair $(\bfhom,\theta)$ is not essential, all the choices are equivalent.

Dually, a symmetric monoidal category $\bfC=(\bfC,\otimes)$ is relatively coclosed with the parametrising subcategory $\bfP\subset\bfC$ iff there exists a bifunctor $\bfcohom\colon\bfP^\op\times\bfP\to\bfC$, called {\it (relative) internal cohom-functor}, with the bijections
\begin{align} \label{cohomDef}
 \vartheta=\vartheta_{W,V,Z}\colon\Hom\big(\bfcohom(V,W),Z\big)\isoright\Hom(W,Z\otimes V)
\end{align}
natural in $V,W\in\bfP$, $Z\in\bfC$. This means that $(-\otimes-,\bfcohom,\vartheta)$ is a left adjunction with a parameter relative to the full subcategory $\bfC'=\bfP$; it is also unique up to an isomorphism. The values $\bfcohom(V,W)$ on the objects $V,W\in\bfP$ are called {\it internal cohom-objects}.

\begin{Prop} \label{PropRelclosedhom}
 Let $\bfC=(\bfC,\otimes)$ be a symmetric monoidal category and $\bfC'\subset\bfC$ be its full monoidal subcategory. Consider also a full subcategory $\bfP\subset\bfC'$.
\begin{itemize}
 \item[{\normalfont(1).}] If $\bfC$ is closed relative to $\bfP$ and $\bfhom(V,W)\in\bfC'$ $\;\forall\,V,W\in\bfP$, then $(\bfC',\otimes)$ is also closed relative to $\bfP$. The internal hom-functor $\bfhom\colon\bfP^\op\times\bfP\to\bfC$ equals to the composition of $\bfhom\colon\bfP^\op\times\bfP\to\bfC'$ with the embedding $\bfC'\hookrightarrow\bfC$.
\item[{\normalfont(2).}] If $\bfC$ is coclosed relative to $\bfP$ and $\bfcohom(V,W)\in\bfC'$ $\;\forall\,V,W\in\bfP$, then $(\bfC',\otimes)$ is also coclosed relative to $\bfP$. The bifunctor $\bfcohom\colon\bfP^\op\times\bfP\to\bfC$ equals to the composition of $\bfcohom\colon\bfP^\op\times\bfP\to\bfC'$ with $\bfC'\hookrightarrow\bfC$.
\end{itemize}
\end{Prop}

\noindent{\bf Proof.} The condition of the point~(1) imply that the bifunctor $\bfhom\colon\bfP^\op\times\bfP\to\bfC$ has the form $\bfhom=E_{\bfC'}\cdot G$ for some $G\colon\bfP^\op\times\bfP\to\bfC'$. For any $V\in\bfP$ the functor $G(V,-)\colon\bfP\to\bfC'$ is a right adjoint for $-\otimes V\colon\bfC'\to\bfC'$ relative to $\bfP$. \qed

\bb{Composition morphisms.}
Let $\bfC=(\bfC,\otimes)$ be a relatively closed monoidal category with the parametrising subcategory $\bfP$. Counits~\eqref{varepsilonYZ} for $(-\otimes-,\bfhom,\theta^{-1})$ are morphisms
\begin{align} \label{evalYZ}
 &\eval_{V,W}=\theta(\id_{\bfhom(V,W)})\colon\bfhom(V,W)\otimes V\to W, &&V,W\in\bfP.
\end{align}
called {\it evaluations}.
The formula~\eqref{alphaeps} takes the form 
\begin{align} \label{thetaEval}
 &\theta_{X,V,W}(f)=\eval_{V,W}\cdot(f\otimes\id_V), &&X\in\bfC,\;\; V,W\in\bfP,\quad f\colon X\to\bfhom(V,W).
\end{align}
This implies the formula
\begin{align} \label{thetafg}
 \theta_{\wt X,V,W}(f\cdot g)=
 \theta_{X,V,W}(f)\cdot(g\otimes\id_V),
\end{align}
where $X,\wt X\in\bfC$, $V,W\in\bfP$ and $f\colon X\to\bfhom(V,W)$, $g\colon\wt X\to X$ are morphisms in $\bfC$ (see~\cite[Prop.~3.1]{Sqrt}).

For any closed monoidal category one can define internal compositions (see~\cite[\S~6.1]{Bor2}). Let us generalise them for the relative case. Consider the morphisms
\begin{align} \label{evaleval}
 \bfhom(V,W)\otimes\bfhom(U,V)\otimes U\xrightarrow{\id\otimes\eval_{U,V}}
 \bfhom(V,W)\otimes V\xrightarrow{\eval_{V,W}}W,
\end{align}
where $U,V,W\in\bfP$.
By applying the map $\theta^{-1}_{X,U,W}$, where $X=\bfhom(V,W)\otimes\bfhom(U,V)$, to the morphism~\eqref{evaleval} we obtain
\begin{align} \label{cUVW}
 &c_{U,V,W}\colon\bfhom(V,W)\otimes\bfhom(U,V)\to\bfhom(U,W).
\end{align}
It is called {\it (internal) composition morphism}. It is the unique morphism of the form~\eqref{cUVW} such that the diagram
\begin{align} \label{cUVWdiag}
 \xymatrix{
 \bfhom(V,W)\otimes\bfhom(U,V)\otimes U\ar[rr]^{\qquad\id\otimes\eval_{U,V}}\ar[d]^{c_{U,V,W}\otimes\id_U} &&
 \bfhom(V,W)\otimes V\ar[d]^{\eval_{V,W}} \\
 \bfhom(U,W)\otimes U\ar[rr]^{\qquad\eval_{U,W}} && W 
}
\end{align}
is commutative.

In the case of relatively coclosed monoidal category $\bfC=(\bfC,\otimes)$ with parametrising subcategory $\bfP$ we have units~\eqref{etaYZ} of the form
\begin{align} \label{coevXY}
 \coev_{V,W}=\vartheta(\id_{\bfhom(V,W)})\colon W\to\bfcohom(V,W)\otimes V, &&V,W\in\bfP.
\end{align}
They are called {\it coevaluations}.
We get the formula
\begin{align} \label{varthetacoev}
 &\vartheta(f)=(f\otimes\id_V)\cdot\coev_{V,W}, &&f\colon\bfcohom(V,W)\to Z.
\end{align}
The {\it (internal) cocomposition morphism} is 
\begin{align} \label{dWVU}
 d_{W,V,U}=\vartheta^{-1}\big((\id\otimes\coev_{U,V})\cdot\coev_{V,W}\big)\colon\bfcohom(U,W)\to\bfcohom(V,W)\otimes\bfcohom(U,V),
\end{align}
where $U,V,W\in\bfP$. This is a unique morphism making the following diagram commutative:
\begin{align} \label{dWVUdiag}
\xymatrix{
W\ar[rrr]^{\coev_{V,W}\qquad}\ar[d]^{\coev_{U,W}} &&& \bfcohom(V,W)\otimes V\ar[d]^{\id\otimes\coev_{U,V}} \\
 \bfcohom(U,W)\otimes U\ar[rrr]^{d_{W,V,U}\otimes\id\qquad}  &&&
 \bfcohom(V,W)\otimes\bfcohom(U,V)\otimes U 
}
\end{align}

\begin{Prop} \label{PropAss}
The composition~\eqref{cUVW} is associative: the diagram
\begin{align} \label{cAss}
 \xymatrix{
 \bfhom(V,W)\otimes\bfhom(U,V)\otimes\bfhom(T,U)\ar[rr]^{\qquad c_{U,V,W}\otimes\id}\ar[d]^{\id\otimes c_{T,U,V}} &&
 \bfhom(U,W)\otimes\bfhom(T,U)\ar[d]^{c_{T,U,W}} \\
 \bfhom(V,W)\otimes\bfhom(T,V)\ar[rr]^{\qquad c_{T,V,W}} &&
 \bfhom(T,W)
}
\end{align}
commutes for any $T,U,V,W\in\bfP$. Dually, cocomposition~\eqref{dWVU} is coassociative in the sense of commutativity of the diagram
\begin{align*}
\xymatrix{
 \bfcohom(T,W)\ar[rr]^{d_{W,U,T}\qquad\quad}\ar[d]^{d_{W,V,T}} &&
 \bfcohom(U,W)\otimes\bfcohom(T,U)\ar[d]^{d_{W,V,U}\otimes\id} \\
 \bfcohom(V,W)\otimes\bfcohom(T,V)\ar[rr]^{\id\otimes d_{V,U,T}\quad\qquad} &&
 \bfcohom(V,W)\otimes\bfcohom(U,V)\otimes\bfcohom(T,U)
}
\end{align*}
\end{Prop}

\noindent{\bf Proof.} By drawing the commutative diagram~\eqref{cUVWdiag} for each composition morphism in~\eqref{cAss} one yields
\begin{multline*}
 \eval_{T,W}\cdot(c_{T,U,W}\otimes\id_T)\cdot(c_{U,V,W}\otimes\id\otimes\id_T)=\eval_{V,W}\cdot(\id\otimes\eval_{U,V})\cdot(\id\otimes\id\otimes\eval_{T,U})=\\
 \eval_{T,W}\cdot(c_{T,V,W}\otimes\id_T)\cdot(\id\otimes c_{T,U,V}\otimes\id_T).
\end{multline*}
This implies commutativity of the diagram~\eqref{cAss} due to the formula~\eqref{thetaEval} and bijectivity of $\theta$ (or, equivalently, due to the universality of $\big(\bfhom(T,W),\eval_{T,W}\big)$). The coassociativity of the cocomposition then follows from the duality principle. \qed

\bb{Internal end.} \label{bbIntend}
 Consider the bijection~\eqref{homDef} with the unit object $X=I_\bfC$ and $W=V$. Denote the image of the identification isomorphism $I\otimes V=V$ under $\theta^{-1}$ by
\begin{align} \label{uV}
 u_V:=\theta^{-1}(\id_V)\colon I_\bfC\to\bfhom(V,V), && V\in\bfP.
\end{align}
This is a unique morphism making the following diagram commute:
\begin{align} \label{uVDiag}
\xymatrix{
 I_\bfC\otimes V\ar[rr]^{u_V\otimes\id_V\quad}\ar@{=}[drr]
 &&\bfhom(V,V)\otimes V\ar[d]^{\eval_{V,V}} \\
 && V
}
\end{align}

\begin{Lem} \label{Lemuc}
 For any $V,W\in\bfP$ we have the commutative diagrams
\begin{align} \label{LemucDiag}
\xymatrix{
 I\otimes\bfhom(V,W)\ar[d]_{u_W\otimes\id}\ar@{=}[r] & \bfhom(V,W)\ar@{=}[r] & \bfhom(V,W)\otimes I\ar[d]^{\id\otimes u_V} \\
 \bfhom(W,W)\otimes\bfhom(V,W)\ar[ru]_{c_{V,W,W}} & & \bfhom(V,W)\otimes\bfhom(V,V)\ar[lu]^{c_{V,V,W}}
}
\end{align}
\end{Lem}

\noindent{\bf Proof.} The commutativity of these diagrams can be shown by means of the formula~\eqref{thetafg}. Let us consider the right diagram~\eqref{LemucDiag}. By applying the bijection $\theta$ one yields
\begin{multline*}
 \theta\big(c_{V,V,W}\cdot(\id\otimes u_V)\big)=\theta(c_{V,V,W})\cdot(\id\otimes u_V\otimes\id)=\eval_{V,W}\cdot(\id\otimes\eval_{V,V})\cdot(\id\otimes u_V\otimes\id)=\\
 \eval_{V,W}\cdot\big(\id\otimes\theta(u_V)\big)=\eval_{V,W}\cdot\big(\id\otimes\id\big)=\theta(\id).
\end{multline*}
The commutativity of the left diagram is proved similarly. \qed

Denote $\bfend(V)=\bfhom(V,V)$. Proposition~\ref{PropAss} for the case $U=V=W$ and Lemma~\ref{Lemuc} for $V=W$ imply that the morphisms $c_V:=c_{V,V,V}\colon\bfend(V)\otimes\bfend(V)\to\bfend(V)$ and $u_V\colon I_\bfC\to\bfend(V)$ give a structure of monoid in $\bfC=(\bfC,\otimes)$ on the object $\bfend(V)\in\bfC$.

Dually, by applying the inverse of the bijection~\eqref{cohomDef}, where $W=V$ and $Z=I_\bfC$, to the identification isomorphism $V=I\otimes V$ we obtain the morphism
\begin{align} \label{vV}
 &v_V:=\vartheta^{-1}(\id_V)\colon\bfcohom(V,V)\to I_\bfC, && V\in\bfP.
\end{align}
This is a unique morphism making the following diagram commute:
\begin{align} \label{vVDiag}
\xymatrix{
 V\ar[rr]^{\coev_{V,V}\qquad}\ar@{=}[drr]
 &&\bfcohom(V,V)\otimes V\ar[d]^{v_V\otimes\id_V} \\
 && I_\bfC\otimes V
}
\end{align}

A symmetric monoidal category $\bfC=(\bfC,\otimes)$ is closed relative to $\bfP$ iff its opposite category $\bfC^\op=(\bfC^\op,\otimes)$ is coclosed relative to $\bfP^\op$. The hom-object $\bfhom(V,W)$ for $\bfC$ is the cohom-object $\bfcohom(V,W)$ for $\bfC^\op$. The morphisms $\theta_{X,V,W}$ and $\eval_{V,W}$ for $\bfC$ are $\vartheta_{W,V,X}$ and $\coev_{V,W}$ for $\bfC^\op$. The composition $c_{U,V,W}$ and morphism $u_V$ for $\bfC$ are the cocomposition $d_{W,V,U}$ and morphism $v_V$ for $\bfC^\op$. By taking into account this duality and argument mentioned above one can state the following.

\begin{Prop}
 If $\bfC=(\bfC,\otimes)$ is relatively closed monoidal category with parametrising subcategory $\bfP\subset\bfC$, then for each $V\in\bfP$ we have the monoid $\bfend(V)=\big(\bfend(V),c_V,u_V\big)$ in $(\bfC,\otimes)$, where $c_V:=c_{V,V,V}$ and $u_V$ is~\eqref{uV}.

 Dually, if $\bfC=(\bfC,\otimes)$ is relatively coclosed monoidal category with parametrising subcategory $\bfP\subset\bfC$, then for each $V\in\bfP$ the morphisms $d_V:=d_{V,V,V}$ and~\eqref{vV} equip the object $\bfcoend(V):=\bfcohom(V,V)$ with a structure of comonoid $\bfcoend(V)=\big(\bfcoend(V),d_V,v_V\big)$ in $(\bfC,\otimes)$.
\end{Prop}

The diagrams~\eqref{cUVWdiag} for $U=V=W$ and~\eqref{uVDiag} imply that $\eval_{V,V}\colon\bfend(V)\otimes V\to V$ is a (left) action of the monoid $\bfend(V)$ on $V$. It also follows from the diagrams~\eqref{dWVUdiag} for $U=V=W$ and \eqref{vVDiag} that $\coev_{V,V}\colon V\to\bfcoend(V)\otimes V$ is a (left) coaction of the comonoid $\bfcoend(V)$ on $V$.

\bb{The hom-functor on morphisms.} The hom-functor as well as cohom-functor is a bifunctor contravariant in the first and covariant in the second argument, i.e. a functor $\bfC^\op\times\bfC\to\bfC$. Let us calculate their values on morphisms in terms of evaluation and coevaluation.

\begin{Prop} \label{Propthetahom}
 Let $f\colon V'\to V$ and $g\colon W\to W'$ be morphisms in $\bfP$. In the relatively closed case we have
\begin{align} \label{thetahom}
 \theta\big(\bfhom(f,g)\big)=g\cdot\eval_{V,W}\cdot(\id_{\bfhom(V,W)}\otimes f).
\end{align}
For the relatively coclosed case the following formula holds:
\begin{align} \label{thetacohom}
 \vartheta\big(\bfcohom(f,g)\big)=(\id_{\bfcohom(V',W')}\otimes f)\cdot\coev_{V',W'}\cdot g.
\end{align}
\end{Prop}

\noindent{\bf Proof.} Repeats the proof of~\cite[Prop.~3.4]{Sqrt}. \qed

\bb{The case $V=I_\bfC$.} Suppose that $I=I_\bfC\in\bfP$. Then we can define the hom-objects $\bfhom(I,W)\in\bfC$ and the monoid $\bfend(I)\in\Mon(\bfC)$ (in the relatively closed case) or the cohom-objects $\bfcohom(I,W)\in\bfC$ and the comonoid $\bfcoend(I)\in\Comon(\bfC)$ (in the relatively coclosed case). We can describe them as follows (see~\cite[Prop.~6.1.8]{Bor2} for the usual closed case).

\begin{Prop} \label{ProphomIW}
 The object $\bfhom(I,W)$ is naturally isomorphic to $W$ via the isomorphism
\begin{align} \label{evalIW}
 \bfhom(I,W)=\bfhom(I,W)\otimes I\xrightarrow{\eval_{I,W}} W.
\end{align}
In particular, $\eval_{I,I}$ gives the identification of objects $\bfend(I)=I$. As a monoid it coincides with $\bfend(I)=(I,\id_{I},\id_{I})$.
The object $\bfcohom(I,W)$ is naturally isomorphic to $W$ via
\begin{align} \label{coevalIW}
 \bfcohom(I,W)=W\xrightarrow{\coev_{I,W}}\bfcohom(I,W)\otimes I.
\end{align}
In particular, $\coev_{I,I}$ gives the identification $\bfcoend(I)=I$. As a comonoid it has the form $\bfcoend(I)=(I,\id_{I},\id_{I})$.
\end{Prop}

\noindent{\bf Proof.} We prove it for the relatively coclosed case. The morphisms~\eqref{coevalIW} are natural in $W\in\bfP$, so we need to prove that they are isomorphisms. Due to the formula~\eqref{varthetacoev} the natural isomorphism
\begin{align} \label{thetaXIW}
 \beta_{W,Z}=\vartheta_{W,I,Z}\colon\Hom\big(\bfcohom(I,W),Z\big)\isoright\Hom(W,Z\otimes I)=\Hom(W,Z)
\end{align}
has the form $\beta_{W,Z}=(\coev_{I,W})^*$. Then, Corollary~\ref{CorYL1} implies that the coevaluations $\coev_{I,W}$ are all isomorphisms.
To show that the last sentence of Prop.~\ref{ProphomIW} is also valid we note that all the arrows in the diagrams~\eqref{dWVUdiag} and \eqref{vVDiag} for $U=V=W=I$ are identifications $\id_I$. This implies that $d_I=\id_I$ and $v_I=\id_I$. \qed

\bb{The case of coreflective relatively coclosed monoidal subcategory.}
 Recall that a subcategory $\bfC'\subset\bfC$ is called {\it coreflective} iff the inclusion $E_{\bfC'}\colon\bfC'\hookrightarrow\bfC$ has a right adjoint $G\colon\bfC\to\bfC'$ (see~\cite[\S~4.3]{Mcl}). Due to Theorem~\ref{ThUnit} this condition means exactly that for any $X\in\bfC$ there exists a universal arrow $(GX,\varepsilon_X)$ from $E_{\bfC'}$ to $X$, where the universality means that for any object $X'\in\bfC'$ and morphism $f\colon X'\to X$ in $\bfC$ there is a unique $h\colon X'\to GX$ such that the diagram
\begin{align} \label{fgeps}
\xymatrix{
X'\ar[r]^h\ar[dr]_f & GX\ar[d]^{\varepsilon_X}\\
 & X
}
\end{align}
is commutative.
We choose the adjunction such that $GX'=X'$, $\eta_{X'}=\varepsilon_{X'}=\id_{X'}$ $\;\forall\,X'\in\bfP$. Since $\varepsilon_X$ is natural in $X$, this choice implies that the morphisms in~\eqref{fgeps} are related as $h=Gf$.

Suppose that the subcategory $\bfC'\subset\bfC$ is monoidal with respect to $\otimes$. In particular, we have $GX\otimes GY\in\bfC'$ for any objects $X,Y\in\bfC$. The commutative diagram~\eqref{fgeps} for $f=\varepsilon_X\otimes\varepsilon_Y\colon GX\otimes GY\to X\otimes Y$ takes the form
\begin{align} \label{GXGY}
\xymatrix{
GX\otimes GY\ar[r]^{\phi_{X,Y}}\ar[dr]_{\varepsilon_X\otimes\varepsilon_Y} & G(X\otimes Y)\ar[d]^{\varepsilon_{X\otimes Y}} \\
& X\otimes Y
}
\end{align}
where $\phi_{X,Y}=G(\varepsilon_X\otimes\varepsilon_Y)$ (this is lax monoidal structure morphisms of the right adjoint $G$). The following fact will be needed in Subsection~\ref{secQR}.

\begin{Th} \label{ThCoref}
 Let $\bfC'$ be a coreflective monoidal subcategory of $(\bfC,\otimes)$ and let $(\bfC',\otimes)$ be coclosed relative to a full subcategory $\bfP\subset\bfC'$. Suppose that $\phi_{X,W}\colon GX\otimes W\to G(X\otimes W)$ is an isomorphism for any $X\in\bfC$ and $W\in\bfP$. Then $(\bfC,\otimes)$ is also coclosed relative to $\bfP$. The internal cohom-functor $\bfcohom\colon\bfP^\op\times\bfP\to\bfC$ is the composition of $\bfcohom\colon\bfP^\op\times\bfP\to\bfC'$ with the embedding $E_{\bfC'}\colon\bfC'\hookrightarrow\bfC$.
\end{Th}

\noindent{\bf Proof.} Let $V,W\in\bfP$ and $X\in\bfC$. The isomorphism $GX\otimes W\cong G(X\otimes W)$, universality of $\big(G(X\otimes W),\varepsilon_{X\otimes W}\big)$ and commutativity of the diagram~\eqref{GXGY} for $Y=W$ imply that for any $f\colon V\to X\otimes W$ there exists a unique $f'\colon V\to GX\otimes W$ such that $(\varepsilon_X\otimes\id_W)\cdot f'=f$. Since $GX\in\bfC'$ and $\bfC'$ is coclosed relative to $\bfP$, there is a unique $h'\colon\bfcohom(W,V)\to GX$ such that $(h'\otimes\id_W)\cdot\coev_{W,V}=f'$, so we have the commutative diagrams
\begin{align*}
\xymatrix{
V\ar[d]_{\coev_{W,V}}\ar[rr]^{f}\ar[drr]^{f'} && X\otimes W \\
\bfcohom(W,V)\otimes W\ar[rr]_{\qquad h'\otimes\id_W} && GX\otimes W\ar[u]_{\varepsilon_X\otimes\id_W}
} \qquad\qquad
\xymatrix{
V\ar[d]_{\coev_{W,V}}\ar[r]^{f} & X\otimes W \\
\bfcohom(W,V)\otimes W\ar[ru]_{\quad h\otimes\id_W}
}
\end{align*}
where $h=\varepsilon_X\cdot h'$. Thus we proved that for any $f\colon V\to X\otimes W$ there exists a morphism $h\colon\bfcohom(W,V)\to X$ making the second diagram commute. Now we only need to prove that it is unique. Let $g\colon\bfcohom(W,V)\to X$ be such that $(g\otimes\id_W)\cdot\coev_{W,V}=f$. Since $\bfcohom(W,V)\in\bfC'$, it factors as $g=\varepsilon_X\cdot g'$ for some $g'\colon\bfcohom(W,V)\to GX$, so we obtain $(\varepsilon_X\otimes\id_W)\cdot(g'\otimes\id_W)\cdot\coev_{W,V}=f=(\varepsilon_X\otimes\id_W)\cdot f'$. Due to the uniqueness of $f'$ we derive $(g'\otimes\id_W)\cdot\coev_{W,V}=f'=(h'\otimes\id_W)\cdot\coev_{W,V}$. By taking into account the uniqueness of $h'$ we get $g'=h'$ and, hence, $g=h$. \qed

\begin{Rem} \normalfont
 By the duality principle one can state an analogous theorem for a reflective monoidal subcategory $(\bfC',\otimes)\subset(\bfC,\otimes)$ closed relative to $\bfP$. In this case the inclusion functor $E_{\bfC'}\colon\bfP\hookrightarrow\bfC$ has a left adjoint functor $F\colon\bfC\to\bfP$, which is colax monoidal. If $\phi_{X,W}\colon F(X\otimes W)\to FX\otimes FW$ are isomorphisms for all $X\in\bfC$ and $W\in\bfP$, then $(\bfC,\otimes)$ is also closed relative to $\bfP$.
\end{Rem}

\begin{Rem} \normalfont
\label{RemCountEx}
 Note that not every symmetric monoidal category is (co)closed relative to a (co)closed monoidal subcategory. Consider a counterexample. Let $\bfC$ be the partial order of subsets $X\subset\RR$ with respect to inclusion. Define the monoidal product as union: $X\otimes Y=X\cup Y$. The subcategory $\bfC'\subset\bfC$ consisting of open subsets is monoidal. It has a monoidal subcategory $\bfP$ of the intervals $(0,x)$, where $x\in\RR_{\ge0}$. By identifying $(0,x)$ with $x$ we can interpret $\bfP$ as the standard normal order on $\RR_{\ge0}$ with the monoidal product $x\otimes y=\max(x,y)$. The monoidal category $(\bfP,\otimes)$ is coclosed with the cohom-functor 
\begin{align} \label{cohomxyP}
 \bfcohom(x,y)=
\begin{cases}
 0, & x\ge y; \\
 y, & x<y.
\end{cases}
\end{align}
Nevertheless, the monoidal category $(\bfC',\otimes)$ is not coclosed relative to $\bfP$. Theorem~\ref{ThCoref} is not applicable here, since the subcategory $\bfP\subset\bfC'$ is not coreflective. The whole category $\bfC$ is coclosed relative to $\bfP$ but with another cohom-functor
\begin{align} \label{cohomxyC}
 \bfcohom(x,y)=\{a>0\mid x\le a<y\}=
\begin{cases}
 \emptyset, & x\ge y; \\
 (0,y), & x=0; \\
 [x,y), & 0<x<y.
\end{cases}
\end{align}
This gives two counterexamples: for $\bfC'\subset\bfC$ and $\bfP\subset\bfC$. The former one shows that Prop.~\ref{PropRelclosedhom} does not work without the conditions on the cohom-functor, while the latter one tells us that it is not necessary condition and the internal (co)hom-functor for a monoidal subcategory may differ from the internal (co)hom-functor for the whole category.
\end{Rem}

\section{General representation theory}
\label{secGRT}

Here we generalise the approach described in~\cite[\S~3]{Sqrt} for the case of relatively (co)closed monoidal categories.

\subsection{Representations in relatively closed monoidal categories}
\label{secRep}

By following the ideas of~\cite[\S~3.2]{Sqrt} we first consider the categories of (co)representations.

\bb{Category of representations.}
Let $\MM=(X,\mu_X,\eta_X)\in\Mon(\bfC)$ be monoid in a relatively closed monoidal category $\bfC=(\bfC,\otimes,I)$ with a parametrising subcategory $\bfP$, where $X\in\bfC$. 

\begin{Def} \label{DefRep} \normalfont 
 A {\it representation} of the monoid $\MM$ on an object $V\in\bfP$ is a morphism $\rho\colon\MM\to\bfend(V)$ in the category $\Mon(\bfC)$. This is a morphism $\rho\colon X\to\bfend(V)$ in $\bfC$ such that the diagrams
\begin{align} \label{rhoDiag}
\xymatrix{
 X\otimes X\ar[rr]^{\mu_X}\ar[d]_{\rho\otimes\rho}& & X\ar[d]^\rho \\
 \bfend(V)\otimes\bfend(V)\ar[rr]^{\qquad c_{V}}& &\bfend(V)
}\qquad\qquad
\xymatrix{
I\ar[r]^{\eta_X}\ar[dr]_{u_V} & X\ar[d]^\rho \\
 & \bfend(V)
}
\end{align}
are commutative. Let $\rho\colon\MM\to\bfend(V)$ and $\tau\colon\MM\to\bfend(W)$ be representations of the same monoid $\MM$ on the objects $V\in\bfP$ and $W\in\bfP$ respectively. A morphism between these representations is a morphism $V\to W$ in $\bfC$ making the following diagram commute:
\begin{align} \label{rhof}
\xymatrix{
X\ar[rr]^{\rho}\ar[d]_{\tau} && \bfhom(V,V)\ar[d]^{\bfhom(\id_V,f)} \\
\bfhom(W,W)\ar[rr]^{\bfhom(f,\id_W)} && \bfhom(V,W)
}
\end{align}
\end{Def}

Denote by $\Rep_\bfP(\MM)$ the category of representations of the monoid $\MM$ with their morphisms. The objects of this category are the pairs $(V,\rho)$, where $V\in\bfP$ and $\rho\colon\MM\to\bfend(V)$ is a representation of $\MM$ on $V$. A morphism $(V,\rho)\to(W,\tau)$ in $\Rep_\bfP(\MM)$ is a morphism $f\colon V\to W$ in $\bfP$ respecting the representations in the sense of diagram~\eqref{rhof}.

\begin{Prop} \label{PropRepIso}
 The forgetful functor
\begin{align} \label{Gff}
 &G\colon\Rep_\bfP(\MM)\to\bfP, &&(V,\rho)\mapsto V, &&f\mapsto f,
\end{align}
reflects isomorphisms: if $(V,\rho)$, $(W,\rho)$ are objects of the category $\Rep_\bfP(\MM)$ and a morphism $f\colon(V,\rho)\to(W,\tau)$ in $\Rep_\bfP(\MM)$ is given by an isomorphism $f\colon V\isoright W$ in $\bfP$, then $f\colon(V,\rho)\to(W,\tau)$ is an isomorphism in $\Rep_\bfP(\MM)$.
\end{Prop}

\noindent{\bf Proof.} Let $f^{-1}\colon W\to V$ be the inverse morphism for $f$. Then, by composing the diagram~\eqref{rhof} with $\bfhom(f^{-1},f^{-1})\colon\bfhom(V,W)\to\bfhom(W,V)$ we see that $f^{-1}$ gives a morphism $(W,\tau)\to(V,\rho)$ in $\Rep_\bfP(\MM)$. \qed

\begin{Rem} \normalfont
 Consider the composition functor
\begin{align} \label{FGGhom}
 \Rep_\bfP(\MM)^\op\times\Rep_\bfP(\MM)\xrightarrow{G^\op\times G}\bfP^\op\times\bfP\xrightarrow{\bfhom}\bfC,
\end{align}
where $G$ is the forgetful functor~\eqref{Gff}. The commutative diagram~\eqref{rhof} means exactly that the formula $\alpha_{(V,\rho)}=\rho\colon X\to\bfhom(V,V)$ defines a dinatural transformation $\alpha$ from the constant functor $X$ to~\eqref{FGGhom}.
\end{Rem}

Let $\bfC=(\bfC,\otimes,I)$ be a relatively coclosed monoidal category with the parametrising subcategory $\bfP$. Consider a comonoid $\OO=(X,\Delta_X,\varepsilon_X)\in\Comon(\bfC)$, where $X\in\bfC$.

\begin{Def} \normalfont
 A {\it corepresentation} of the comonoid $\OO$ on an object $V\in\bfP$ is a morphism $\omega\colon\bfcoend(V)\to\OO$ in the category $\Comon(\bfC)$. This is a morphism $\omega\colon\bfcoend(V)\to X$ in $\bfC$ such that the diagrams
\begin{align*}
\xymatrix{
 \bfcoend(V)\ar[d]_{d_{V}}\ar[rr]^{\qquad\omega}& & X\ar[d]^{\Delta_X} \\
 \bfcoend(V)\otimes\bfcoend(V)\ar[rr]^{\qquad\qquad\omega\otimes\omega}& &X\otimes X
}\qquad\qquad
\xymatrix{
\bfcoend(V)\ar[r]^{\qquad\omega}\ar[dr]_{v_V} & X\ar[d]^{\varepsilon_X} \\
 &I 
}
\end{align*}
are commutative. Denote by $\Corep_\bfP(\OO)$ the category whose objects are the pairs $(V,\omega)$, where $V\in\bfP$ and $\omega$ is a corepresentation of $\OO$ on $V$, and morphisms $(V,\omega)\to(W,\nu)$ are morphisms $f\colon V\to W$ in $\bfP$ such that
\begin{align} \label{omegaf}
\xymatrix{
 \bfcohom(W,V)\ar[rrr]^{\bfcohom(f,\id_V)}\ar[d]^{\bfcohom(\id_W,f)} &&& \bfcohom(V,V)\ar[d]^{\omega} \\
\bfcohom(W,W)\ar[rrr]^{\qquad\nu} &&& X
}
\end{align}
\end{Def}

The categorical duality allows us to pass from relatively closed categories to relatively coclosed categories (see p.~\ref{bbIntend}). For any monoidal category $\bfC=(\bfC,\otimes)$ we have the correspondence
\begin{align} \label{MonComon}
 &\Mon(\bfC)=\Comon(\bfC^\op)^\op,
\end{align}
which identifies the monoid $\MM=(X,\mu_X,\eta_X)$ in $\bfC$ with the comonoid $\OO=(X,\Delta_X,\varepsilon_X)$ in the opposite monoidal category $\bfC^\op=(\bfC^\op,\otimes)$, where $\Delta_X\colon X\to X\otimes X$ and $\varepsilon_X\colon X\to I_\bfC$ are the morphisms $\mu_X\colon X\otimes X\to X$ and $\eta_X\colon  I_\bfC\to X$ regarded as morphisms in $\bfC^\op$.
If $\bfC=(\bfC,\otimes)$ is symmetric and closed relative to $\bfP$, then the categories of representations of the monoid $\MM$ and corepresentations of the corresponding comonoid $\OO$ are opposite to each other. A representation $\rho\colon\MM\to\bfend(V)$ on $V\in\bfP$ is identified with the corepresentation $\omega\colon\bfcoend(V)\to\OO$ by the functor~\eqref{MonComon} as morphisms of monoids and comonoids respectively, where $\bfend(V)$ and $\bfcoend(V)$ are end-monoid and coend-monoid for $\bfC$ and $\bfC^\op$ respectively identified by~\eqref{MonComon}. Moreover, if $\tau\colon\MM\to\bfend(W)$ is a representation of $\MM$ on $W\in\bfP$ and $\nu\colon\bfcoend(W)\to\OO$ is the corresponding corepresentation, then a morphism $f\colon V\to W$ in $\bfC$ is a morphism $(V,\rho)\to(W,\tau)$ in $\Rep_\bfP(\MM)$ iff it is a morphism $(W,\nu)\to(V,\omega)$ in $\Corep_\bfP(\OO)$. Thus we have the isomorphism of categories
\begin{align} \label{RepCorep}
 &\Rep_\bfP(\MM)=\Corep_{\bfP^\op}(\OO)^\op,
\end{align}
which identifies $(V,\rho)$ with $(V,\omega)$ and $f\colon(V,\rho)\to(W,\tau)$ with $f\colon(W,\nu)\to(V,\omega)$.

\bb{Representations as actions.} 
Recall that for any (symmetric) monoidal category $\bfC=(\bfC,\otimes)$ the pairs $(V,a)$ of an object $V\in\bfC$ with an action $a\colon X\otimes V\to V$ of a monoid $\MM=(X,\mu_X,\eta_X)\in\Mon(\bfC)$ form a category denoted by $\Lact(\MM)$. For a comonoid $\OO=(X,\Delta_X,\varepsilon_X)\in\Comon(\bfC)$ the pairs $(V,\delta)$ of $V\in\bfC$ with a coaction $\delta\colon V\to X\otimes V$ form the category $\Lcoact(\OO)$.

\begin{Th} \label{ThRepLact}
 Fix a parametrising subcategory $\bfP\subset\bfC$. If $\bfC$ is relatively closed, then for any $V\in\bfP$ the bijection
\begin{align} \label{thetaXVV}
 &\theta_{X,V,V}\colon\Hom\big(X,\bfend(V)\big)\isoright\Hom\big(X\otimes V,V\big)
\end{align}
establishes one-to-one correspondence between representations $\rho\colon\MM\to\bfend(V)$ and actions $a\colon X\otimes V\to V$ of $\MM$ on $V$, i.e. $\rho\in\Hom\big(X,\bfend(V)\big)$ is a representation iff $a=\theta(\rho)$ is an action. A morphism $f\colon V\to W$ is a morphism of representations $(V,\rho)\to(W,\tau)$ iff it is a morphism of actions $(V,a)\to(W,b)$, where $b=\theta(\tau)$. Thus we can consider $\Rep_\bfP(\MM)$ as a full subcategory of $\Lact(\MM)$ via the fully faithful functor
\begin{align} \label{RepLact}
 &\Rep_\bfP(\MM)\hookrightarrow\Lact(\MM), &&(V,\rho)\mapsto\big(V,\theta(\rho)\big).
\end{align}

Analogously, if $\bfC$ is relatively coclosed, then the bijections
\begin{align} \label{thetaP}
 &\vartheta_{V,V,X}\colon\Hom\big(\bfcoend(V),X\big)\cong\Hom\big(V,X\otimes V\big), &&V\in\bfP,
\end{align}
give the fully faithful functor
\begin{align} \label{CorepLcoact}
 &\Corep_\bfP(\OO)\hookrightarrow\Lcoact(\OO), &&(V,\omega)\mapsto\big(V,\vartheta(\omega)\big).
\end{align}
If $\bfC$ is closed or coclosed (i.e. $\bfP=\bfC$), then the functor~\eqref{RepLact} or, respectively, \eqref{CorepLcoact} is an equivalence (more exactly, isomorphism) of categories.
\end{Th}

\noindent{\bf Proof.} Repeats the proof of~\cite[Lemmas~3.3, 3.5]{Sqrt}. \qed

In the relatively closed case the identity morphism
\begin{align}\label{repId}
 X=\bfend(V)\xrightarrow{\id}\bfend(V)
\end{align}
is a representation of the monoid $\bfend(V)=(X,c_V,u_V)$ on $V\in\bfP$.
The corresponding action is $\theta(\id)=\eval_{V,V}\colon\bfend(V)\otimes V\to V$, which was mentioned in p.~\ref{bbIntend}.
In the relatively coclosed case we have an analogical corepresentation
\begin{align}\label{corepId}
 \bfcoend(V)\xrightarrow{\id}\bfcoend(V)=X
\end{align}
of the comonoid $\bfcoend(V)=(X,d_V,v_V)$ on the object $V\in\bfP$. It corresponds to the coaction $\vartheta(\id)=\coev_{V,V}\colon V\to\bfcoend(V)\otimes V$.

\subsection{Tensor product of representations}
\label{secTP}

The monoidal product of representations (actions) of a monoid $\MM$ is introduced by means of an additional structure on $\MM$, namely by a structure of comonoid which turn $\MM$ into a bimonoid. In~\cite{Sqrt} we considered monoidal product in $\Rep_\bfP(\MM)$ via a monoidal product in $\Lact(\MM)$. In the present work we define and study the former monoidal product without the notion of (co)action, but then we show that these two ways give the same result. Note also that here we consider more general case then in~\cite{Sqrt}, the case of relatively (co)closed monoidal categories.

In this subsection we always suppose that the monoidal category $\bfC=(\bfC,\otimes,I)$ is symmetric. The corresponding symmetric structure is a natural isomorphism, which we denote by $\sigma$. It has components $\sigma=\sigma_{X,Y}\colon X\otimes Y\isoright Y\otimes X$. Let $X_1,\ldots,X_n$ be objects of $\bfC$. We use the following notation: for $a=1,\ldots,n-1$ we denote by $\sigma^{(a,a+1)}$ the isomorphism
 $X_1\otimes\cdots\otimes X_a\otimes X_{a+1}\otimes\cdots\otimes X_n\isoright X_1\otimes\cdots\otimes X_{a+1}\otimes X_a\otimes\cdots\otimes X_n$, where $\sigma$
acts on the $a$-th and $(a+1)$-st factors, that is
\begin{multline*}
 \sigma^{(a,a+1)}=\id^{\otimes(a-1)}\otimes\sigma\otimes\id^{\otimes(n-a-1)}=\id_{X_1}\otimes\cdots\otimes\id_{X_{a-1}}\otimes\sigma_{X_a,X_{a+1}}\otimes\id_{X_{a+2}}\otimes\cdots\otimes\id_{X_n}.
\end{multline*}

\bbo{Monoidal structure on the category of monoids.} \label{bbMonC}
Recall that the symmetric structure of $\bfC$ allows to define a monoidal structure on the category of monoids $\Mon(\bfC)$. Let us describe this structure in details.

The monoidal product of two monoids $(X,\mu_X,\eta_X),(Y,\mu_Y,\eta_Y)\in\Mon(\bfC)$ is the triple $(X\otimes Y,\mu_{X\otimes Y},\eta_{X\otimes Y})$, where $\mu_{X\otimes Y}$ and $\eta_{X\otimes Y}$ are defined by the commutative diagrams
\begin{align} \label{muXYDef}
\xymatrix{X\otimes Y\otimes X\otimes Y\ar[d]_{\sigma^{(23)}}\ar[rrd]^{\mu_{X\otimes Y}}   \\
 X\otimes X\otimes Y\otimes Y\ar[rr]^{\mu_X\otimes\mu_Y} && X\otimes Y
}  \qquad\qquad
\xymatrix{I\ar@{=}[d]\ar[rrd]^{\eta_{X\otimes Y}} \\
 I\otimes I\ar[rr]^{\eta_X\otimes\eta_Y}&& X\otimes Y
}
\end{align}
One can check that it is also a monoid in $\bfC$. The monoidal product of morphisms in $\Mon(\bfC)$ is obtained as follows.

\begin{Prop} \label{PropfgTP}
 Let $(X,\mu_X,\eta_X),(Y,\mu_Y,\eta_Y),(X',\mu_{X'},\eta_{X'}),(Y',\mu_{Y'},\eta_{Y'})\in\Mon(\bfC)$. If $f\colon X\to Y$ and $g\colon Y\to Y'$ are morphisms of monoids $(X,\mu_X,\eta_X)\to(X',\mu_{X'},\eta_{X'})$ and $(Y,\mu_Y,\eta_Y)\to(Y',\mu_{Y'},\eta_{Y'})$, then $f\otimes g\colon X\otimes Y\to X'\otimes Y'$ is a morphism
 $$(X\otimes Y,\mu_{X\otimes Y},\eta_{X\otimes Y})\to(X'\otimes Y',\mu_{X'\otimes Y'},\eta_{X'\otimes Y'}).$$
\end{Prop}

\noindent{\bf Proof.} It follows from the commutative diagrams
\begin{align*}
 \xymatrix{
 X\otimes Y\otimes X\otimes Y\ar[r]^{\sigma^{(23)}}\ar[d]^{f\otimes g\otimes f\otimes g} &
 X\otimes X\otimes Y\otimes Y\ar[rr]^{\qquad\mu_X\otimes\mu_Y}\ar[d]^{f\otimes f\otimes g\otimes g} && X\otimes Y\ar[d]^{f\otimes g} && I\otimes I\ar[ll]_{\eta_X\otimes\eta_Y}\ar[lld]^{\eta_{X'}\otimes\eta_{Y'}} \\
 X'\otimes Y'\otimes X'\otimes Y'\ar[r]^{\sigma^{(23)}} &
 X'\otimes X'\otimes Y'\otimes Y'\ar[rr]^{\qquad\mu_{X'}\otimes\mu_{Y'}} && X'\otimes Y'
}
\end{align*}
and the definitions~\eqref{muXYDef}. \qed

Thus we have a bifunctor $-\otimes-\colon\Mon(\bfC)\times\Mon(\bfC)\to\Mon(\bfC)$. It is easy to check that $\Mon(\bfC)=\big(\Mon(\bfC),\otimes\big)$ is a symmetric monoidal category with the unit object $I_{\Mon(\bfC)}=(I,\id_I,\id_I)$. The category of comonoids $\Comon(\bfC)$ is equipped with a structure of symmetric monoidal category in the dual way.

\bb{Natural transformation $\pi$.}
 Until the end of this subsection we suppose that the subcategory $\bfP\subset\bfC$ is monoidal (this condition is reduced to $I\in\bfP$, $V\otimes V'\in\bfP\;\forall\,V,V'\in\bfP$, since $\bfP$ is a full subcategory of $\bfC$). Then for any $V,V',W,W'\in\bfP$ we have the hom-object $\bfhom(V\otimes V',W\otimes W')$. In particular, we have the monoids $\bfend(V\otimes V')$ and $\bfend(I)$. By Prop.~\ref{ProphomIW} the latter monoid is the unit object: $\bfend(I)=(I,\id_I,\id_I)$. We are going to construct some useful morphisms of monoids $\bfend(V)\otimes\bfend(V')\to\bfend(V\otimes V')$.

Let us use the notation $[V,W]:=\bfhom(V,W)$ for briefness.
Consider the composition
\begin{align} \label{evevsigma}
 [V,W]\otimes[V',W']\otimes V\otimes V'\xrightarrow{\sigma^{(23)}}
 [V,W]\otimes V\otimes[V',W']\otimes V'\xrightarrow{\eval\otimes\eval}
 W\otimes W',
\end{align}
where $V,V',W,W'\in\bfP$. By applying the bijection
\begin{multline*}
 \theta^{-1}\colon\Hom\big([V,W]\otimes[V',W']\otimes V\otimes V',W\otimes W'\big)\isoright \\
 \Hom\big([V,W]\otimes[V',W'],[V\otimes V',W\otimes W']\big)
\end{multline*}
to~\eqref{evevsigma} we get the morphism $\pi=\pi^{V,V'}_{W,W'}=\theta^{-1}\big((\eval_{V,W}\otimes\eval_{V',W'})\cdot\sigma^{(23)}\big)$ of the form
\begin{align} \label{piDef}
 \pi\colon[V,W]\otimes[V',W']\to[V\otimes V',W\otimes W'].
\end{align}
This is a unique morphism making the diagram
\begin{align} \label{piDiag}
\xymatrix{
 [V,W]\otimes[V',W']\otimes V\otimes V'\ar[rr]^{\sigma^{(23)}}\ar[d]^{\pi\otimes\id} &&
 [V,W]\otimes V\otimes[V',W']\otimes V'\ar[d]^{\eval\otimes\eval} \\
 [V\otimes V',W\otimes W']\otimes V\otimes V'\ar[rr]^{\eval} &&
 W\otimes W'
}
\end{align}
commute.

\begin{Prop} \label{ProppiNat}
 The morphisms $\pi^{V,V'}_{W,W'}$ are natural in all four arguments $V,V',W,W'\in\bfP$.
\end{Prop}

\noindent{\bf Proof.} We need to prove the commutativity of the diagram
\begin{align} \label{piNatDiag}
\xymatrix{
 [V,W]\otimes[V',W']\ar[r]^\pi\ar[d]_{\bfhom(f,g)\otimes\bfhom(f',g')} &
 [V\otimes V',W\otimes W']\ar[d]^{\bfhom(f\otimes f',g\otimes g')} \\
 [U,Z]\otimes[U',Z']\ar[r]^\pi &
 [U\otimes U', Z\otimes Z']
}
\end{align}
for arbitrary morphisms $f\colon U\to V$, $f\colon U'\to V'$, $g\colon W\to Z$ and $g'\colon W'\to Z'$ in $\bfP$. Let us apply $\theta$ to the compositions in this diagram. This amounts to monoidal multiplication of each object and arrow in the diagram by $U\otimes U'$ from the right and composition with $\eval\colon[U\otimes U',Z\otimes Z']\otimes U\otimes U'\to Z\otimes Z'$. By using the commutativity of~\eqref{piDiag} and naturality of $\sigma$ we get the equivalent diagram
\begin{align*}
\xymatrix{
 [V,W]\otimes[V',W']\otimes U\otimes U'\ar[rr]^{\pi\otimes\id_U\otimes\id_{U'}}
\ar[d]_{\sigma^{(23)}} &&
 [V\otimes V',W\otimes W']\otimes U\otimes U'
\ar[d]^{\bfhom(f\otimes f',g\otimes g')\otimes\id_{U\otimes U'}} \\
 [V,W]\otimes U\otimes[V',W']\otimes U'
\ar[d]_{\bfhom(f,g)\otimes\id_U\otimes\bfhom(f',g')\otimes\id_{U'}} &&
 [U\otimes U', Z\otimes Z']\otimes U\otimes U'\ar[d]^{\eval} \\
 [U,Z]\otimes U\otimes[U',Z']\otimes U'\ar[rr]^{\eval\otimes\eval}&&
Z\otimes Z'
}
\end{align*}
By virtue of Prop.~\ref{Propthetahom} we obtain
\begin{align*}
 &\eval\cdot\big(\bfhom(f\otimes f',g\otimes g')\otimes\id_{U\otimes U'}\big)=\theta\big(\bfhom(f\otimes f',g\otimes g')\big)=(g\otimes g')\cdot\eval\cdot(\id\otimes f\otimes f'), \\[10pt]
&(\eval\otimes\eval)\cdot\big(\bfhom(f,g)\otimes\id_U\otimes\bfhom(f',g')\otimes\id_{U'}\big)=\theta\big(\bfhom(f,g)\big)\otimes\theta\big(\bfhom(f',g')\big)=\\
&\big(g\cdot\eval\cdot(\id\otimes f)\big)\otimes\big(g'\cdot\eval\cdot(\id\otimes f')\big),
\end{align*}
so the diagram~\eqref{piNatDiag} is equivalent to
\begin{align} \label{piNatDiagEval}
\xymatrix{
 [V,W]\otimes[V',W']\otimes U\otimes U'\ar[rr]^{\id\otimes\id\otimes f\otimes f'}
\ar[d]_{\sigma^{(23)}} &&
 [V,W]\otimes[V',W']\otimes V\otimes V'
\ar[d]^{\pi\otimes\id\otimes\id} \\
 [V,W]\otimes U\otimes[V',W']\otimes U'
\ar[d]_{\id\otimes f\otimes\id\otimes f'} &&
 [V\otimes V', W\otimes W']\otimes V\otimes V'\ar[d]^{\eval} \\
 [V,W]\otimes V\otimes[V',W']\otimes V'\ar[rr]^{\eval\otimes\eval}&&
W\otimes W'\ar[r]^{g\otimes g'} & Z\otimes Z'
}
\end{align}
By drawing the arrow $\sigma^{(23)}$ from the top right corner to the down left one we see that the commutativity of the diagram~\eqref{piNatDiagEval} follows from the naturality of $\sigma$ and commutativity of~\eqref{piDiag}. \qed

\bb{Properties of $\pi$.}
The natural transformation $\pi$ is agreed with the internal composition~\eqref{cUVW} in the following sense.

\begin{Prop} \label{PropUVW}
 For any objects $U,V,W,U',V',W'\in\bfP$ we get a commutative diagram
\begin{align*}
\xymatrix{
 [V,W]\otimes[V',W']\otimes[U,V]\otimes[U',V']\ar[r]^{\sigma^{(23)}}\ar[d]^{\pi\otimes\pi} &
 [V,W]\otimes[U,V]\otimes[V',W']\otimes[U',V']\ar[d]^{c_{U,V,W}\otimes c_{U',V',W'}} \\
 [V\otimes V',W\otimes W']\otimes[U\otimes U',V\otimes V']\ar[dr]_{c_{U\otimes U',V\otimes V',W\otimes W'}\qquad} &
 [U,W]\otimes[U',W']\ar[d]^{\pi}  \\
& [U\otimes U',W\otimes W'] 
}
\end{align*}
where $c_{U,V,W}$ are the internal compositions~\eqref{cUVW}.
\end{Prop}

\noindent{\bf Proof.} We need to show the equality of two morphisms $X\to[U\otimes U',W\otimes W']$, where $X=[V,W]\otimes[V',W']\otimes[U,V]\otimes[U',V']$. One can show that their images under $\theta_{X,U\otimes U',W\otimes W'}$ coincide with each other by using the naturality of $\sigma$ and the diagrams~\eqref{piDiag},~\eqref{cUVWdiag} in an appropriate way. However, the proof in terms of diagrams would take too much space. We present it more compactly. By using the formula~\eqref{thetafg} and the definition of $c_{U\otimes U',V\otimes V',W\otimes W'}$ we obtain
\begin{multline}
 \theta\big(c_{U\otimes U',V\otimes V',W\otimes W'}\cdot(\pi\otimes\pi)\big)=\theta(c_{U\otimes U',V\otimes V',W\otimes W'})\cdot(\pi\otimes\pi\otimes\id_{U\otimes U'})=\\
\eval_{V\otimes V',W\otimes W'}\cdot(\id\otimes\eval_{U\otimes U',V\otimes V'})\cdot(\pi\otimes\pi\otimes\id_{U\otimes U'}). \label{PropUVWp1}
\end{multline}
Note that $(\id\otimes\eval_{U\otimes U',V\otimes V'})\cdot(\pi\otimes\id)=\pi\otimes\eval_{U\otimes U',V\otimes V'}=(\pi\otimes\id_{V\otimes V'})\cdot(\id\otimes\eval_{U\otimes U',V\otimes V'})$, so~\eqref{PropUVWp1} equals
\begin{multline}
 \eval_{V\otimes V',W\otimes W'}\cdot(\pi\otimes\id_{V\otimes V'})\cdot(\id\otimes\eval_{U\otimes U',V\otimes V'})\cdot(\id\otimes\pi\otimes\id_{U\otimes U'})=\\
 (\eval_{V,W}\otimes\eval_{V',W'})\cdot\sigma^{(23)}\cdot(\id\otimes\id\otimes\eval_{U,V}\otimes\eval_{U',V'})\cdot\sigma^{(45)}, \label{PropUVWp2}
\end{multline}
where we used the commutative diagram~\eqref{piDiag}. Then, by taking into account the equality $\sigma^{(23)}\cdot(\id\otimes\id\otimes\eval_{U,V}\otimes\eval_{U',V'})=(\id\otimes\eval_{U,V}\otimes\id\otimes\eval_{U',V'})\cdot\sigma^{(34)}\cdot\sigma^{(23)}$ and again the definition of internal composition we can rewrite~\eqref{PropUVWp2} in the form
\begin{multline*}
 (\eval_{U,W}\otimes\eval_{U',W'})\cdot(c_{U,V,W}\otimes\id_U\otimes c_{U',V',W'}\otimes\id_{U'})\cdot\sigma^{(34)}\cdot\sigma^{(23)}\cdot\sigma^{(45)}= \\
 (\eval_{U,W}\otimes\eval_{U',W'})\cdot\sigma^{(23)}\cdot(c_{U,V,W}\otimes c_{U',V',W'}\otimes\id_U\otimes\id_{U'})\cdot\sigma^{(23)}= \\
\theta(\pi)\cdot(c_{U,V,W}\otimes c_{U',V',W'}\otimes\id_{U\otimes U'})\cdot(\sigma^{(23)}\otimes\id_{U\otimes U'})=
\theta\big(\pi\cdot(c_{U,V,W}\otimes c_{U',V',W'})\cdot\sigma^{(23)}\big).
\end{multline*}
Since $\theta$ is a bijection, the proposition is proved. \qed

Let us establish `associativity' of the natural transformation $\pi$.

\begin{Prop} \label{PropAsspi}
 For any $V,W,V',W',V'',W''\in\bfP$ we have the commutative diagram
\begin{align}
\xymatrix{
 [V,W]\otimes[V',W']\otimes[V'',W'']\ar[rr]^{\pi\otimes\id}\ar[d]^{\id\otimes\pi} &&
 [V\otimes V',W\otimes W']\otimes[V'',W'']\ar[d]^{\pi} \\
 [V,W]\otimes[V'\otimes V'',W'\otimes W'']\ar[rr]^{\pi} &&
 [V\otimes V'\otimes V'',W\otimes W'\otimes W'']
}
\end{align}
\end{Prop}

\noindent{\bf Proof.} By using the diagram~\eqref{piDiag} one can show that the morphism $\theta\big(\pi\cdot(\pi\otimes\id)\big)$ coincides with the composition $[V,W]\otimes[V',W']\otimes[V'',W'']\otimes V\otimes V'\otimes V''\xrightarrow{\sigma^{(23)}\sigma^{(45)}\sigma^{(34)}}[V,W]\otimes V\otimes[V',W']\otimes V'\otimes[V'',W'']\otimes V''\xrightarrow{\eval_{V,W}\otimes\eval_{V',W'}\otimes\eval_{V'',W''}}W\otimes W'\otimes W''$. Analogously, one derives that
$\theta\big(\pi\cdot(\id\otimes\pi)\big)$ coincides with the same composition. \qed

\begin{Prop} \label{Proppisigma}
The natural transformation $\pi$ is in agreement with the symmetric structure $\sigma$ in the sense that the diagram
\begin{align} \label{pisigmaDiag}
\xymatrix{
 [V,W]\otimes[V',W']\ar[rr]^{\pi}\ar[d]^{\sigma} &&
 [V\otimes V',W\otimes W']\ar[d]^{\bfhom(\sigma_{V',V},\sigma_{W,W'})} \\
 [V',W']\otimes[V,W]\ar[rr]^{\pi} &&
[V'\otimes V,W'\otimes W] 
}
\end{align}
commutes, where
\begin{gather}
 \sigma=\sigma_{[V,W],[V',W']}\colon[V,W]\otimes[V',W']\to[V',W']\otimes[V,W],  \label{sigma1} \\
 \sigma_{V',V}\colon V'\otimes V\to V\otimes V', \quad\quad \sigma_{W,W'}\colon W\otimes W'\to W'\otimes W  \label{sigma2}
\end{gather}
are the corresponding components of the natural transformation $\sigma$.
\end{Prop}

\noindent{\bf Proof.} Consider the morphism $\sigma^{-1}=\sigma_{[V',W'],[V,W]}$ inverse to~\eqref{sigma1}. Together with the morphisms~\eqref{sigma2} and naturality of $\theta^{-1}$ it gives the commutative diagram
\begin{align} \label{HomsigmathetaDiag}
\xymatrix{
 \Hom\big([V,W]\otimes[V',W']\otimes V\otimes V',W\otimes W'\big)\ar[r]^{\theta^{-1}} \ar[d]_{(\sigma^{-1}\otimes\sigma_{V',V})^*\cdot(\sigma_{W,W'})_*} &
 \Hom\big([V,W]\otimes[V',W'],[V\otimes V',W\otimes W']\big)
\ar[d]_{(\sigma^{-1})^*\cdot\bfhom(\sigma_{V',V},\sigma_{W,W'})_*} \\
 \Hom\big([V',W']\otimes[V,W]\otimes V'\otimes V,W'\otimes W\big)\ar[r]^{\theta^{-1}} &
 \Hom\big([V',W']\otimes[V,W],[V'\otimes V,W'\otimes W]\big)
}
\end{align}
Recall that $\pi=\theta\big((\eval\otimes\eval)\cdot\sigma^{(23)}\big)$. The commutativity of the diagram
\begin{align}
 \xymatrix{
 [V,W]\otimes[V',W']\otimes V\otimes V'\ar[r]^{\sigma^{(23)}} \ar[d]^{\sigma\otimes\sigma_{V,V'}} &
[V,W]\otimes V\otimes[V',W']\otimes V'\ar[rr]^{\qquad\qquad\eval\otimes\eval} \ar[d]^{\sigma_{[V,W]\otimes V,[V',W']\otimes V'}} &&
 W\otimes W'\ar[d]^{\sigma_{W,W'}} \\
 [V',W']\otimes[V,W]\otimes V'\otimes V\ar[r]^{\sigma^{(23)}} &
[V',W']\otimes V'\otimes[V,W]\otimes V\ar[rr]^{\qquad\qquad\eval\otimes\eval} &&
 W'\otimes W
}
\end{align}
implies that the left vertical arrow of the diagram~\eqref{HomsigmathetaDiag} maps the morphism $(\eval\otimes\eval)\cdot\sigma^{(23)}$ to $(\eval\otimes\eval)\cdot\sigma^{(23)}$. Hence, by applying~\eqref{HomsigmathetaDiag} to this morphism we obtain the needed equality $\pi=\bfhom(\sigma_{V',V},\sigma_{W,W'})\cdot\pi\cdot\sigma^{-1}$. \qed

\bb{The morphisms $\pi$ for internal ends.}
If $W=V$ and $W'=V'$, then the morphism~\eqref{piDef} takes the form
\begin{align} \label{piVV}
 \pi\colon\bfend(V)\otimes\bfend(V')\to\bfend(V\otimes V').
\end{align}
The domain object of this morphism is equipped with a structure of monoid as a monoidal product of two monoids; namely, the multiplication morphism for this monoid is $(c_V\otimes c_{V'})\cdot\sigma^{(23)}$, while the unit is $u_V\otimes u_{V'}$. The codomain of~\eqref{piVV} is a monoid with the structure $(c_{V\otimes V'},u_{V\otimes V'})$. Prop.~\ref{PropUVW} for $U=W=V$ and $U'=W'=V'$ means exactly that~\eqref{piVV} preserves the multiplications. Let us prove that it also preserves the units.

\begin{Lem} \label{Lempiu}
The following diagram is commutative:
\begin{align}
\xymatrix{
 I\otimes I\ar[rr]^{u_V\otimes u_{V'}\qquad}\ar@{=}[d] && \bfend(V)\otimes\bfend(V')\ar[d]^{\pi} \\
 I\ar[rr]^{u_{V\otimes V'}\qquad} && \bfend(V\otimes V')
}
\end{align}
\end{Lem}

\noindent{\bf Proof.} Due to the formula~\eqref{thetaEval} and commutativity of the diagram~\eqref{piDiag} we obtain
\begin{multline} \label{LempiuPr}
 \theta\big(\pi\cdot(u_V\otimes u_{V'})\big)=\eval_{V\otimes V',V\otimes V'}\cdot(\pi\otimes\id_{V\otimes V'})\cdot(u_V\otimes u_{V'}\otimes\id_{V\otimes V'})= \\
 (\eval\otimes\eval)\cdot\sigma^{(23)}\cdot(u_V\otimes u_{V'}\otimes\id_V\otimes\id_{V'})=(\eval\otimes\eval)\cdot(u_V\otimes\id_V\otimes u_{V'}\otimes\id_{V'})\cdot\sigma^{(23)}. 
\end{multline}
From the diagram~\eqref{uVDiag} we see that the right hand side of the formula~\eqref{LempiuPr} is the identification $I\otimes I\otimes V\otimes V'=V\otimes V'$, so it equals $\theta(u_{V\otimes V'})$. \qed

\begin{Prop} \label{Proppi}
 The morphisms~\eqref{piVV} are morphisms of monoids.
\end{Prop}

\noindent{\bf Proof.} This is a direct consequence of Prop.~\ref{PropUVW} and Lemma~\ref{Lempiu}. \qed
 
\begin{Prop} \label{ProppiII}
 Both morphisms
\begin{align*}
 &\pi^{I,V}_{I,W}\colon[I,I]\otimes[V,W]\to[I\otimes V,I\otimes W], &
 &\pi^{V,I}_{W,I}\colon[V,W]\otimes[I,I]\to[V\otimes I,W\otimes I]
\end{align*}
are equal to $\id_{[V,W]}$. In particular, the morphism~\eqref{piVV} for $V=I$ or $V'=I$ takes the form of the identification
 $\bfend(I)\otimes\bfend(V')=I\otimes\bfend(V')=\bfend(V')$ or, respectively, it coincides with
 $\bfend(V)\otimes\bfend(I)=\bfend(V)\otimes I=\bfend(V)$.
\end{Prop}

\noindent{\bf Proof.} By substituting $V'=W'=I$ (or $V=W=I$) to the diagram~\eqref{piDiag} and using the identifications $I\otimes V=V=I\otimes V$, $[I,I]=I$ we see that $\pi=\id_{[V,W]}$ makes it commute. Then the statement of the proposition follows from the uniqueness of $\pi$. \qed

\begin{Lem} \label{Lempisigma}
 The morphisms~\eqref{piVV} make the following diagram commute:
\begin{align}
\xymatrix{
 \bfend(V)\otimes\bfend(V')\ar[d]_{\sigma}\ar[r]^{\pi} &
\bfend(V\otimes V')\ar[d]_{\bfhom(\sigma_{V',V},\sigma_{V,V'})}\ar[rrrd]^{\bfhom(\id,\sigma_{V,V'})} \\
\bfend(V')\otimes\bfend(V)\ar[r]_{\pi} & 
\bfend(V'\otimes V)\ar[rrr]_{\bfhom(\sigma_{V,V'},\id)\qquad} &&& [V\otimes V',V'\otimes V]
}
\end{align}
\end{Lem}

\noindent{\bf Proof.} This follows from Prop.~\ref{Proppisigma} and the equalities $\bfhom(\sigma_{V,V'},\id)\cdot\bfhom(\sigma_{V',V},\sigma_{V,V'})=\bfhom(\sigma_{V,V'}\cdot\sigma_{V',V},\sigma_{V,V'})=\bfhom(\id,\sigma_{V,V'})$. \qed

\bb{Natural transformation $\varkappa$.}
In the relatively coclosed case we can consider cohom-objects $\bfcohom(V\otimes V',W\otimes W')$ for any objects $V,V',W,W'\in\bfP$. In particular, the have comonoids $\bfcoend(V\otimes V')$ and $\bfcoend(I)=I_{\Comon(\bfC)}=(I,\id_I,\id_I)$.

To define morphisms dual to~\eqref{piDef} we apply the bijection
\begin{multline*}
 \vartheta^{-1}\colon\Hom\big(W\otimes W',\bfcohom(V,W)\otimes\bfcohom(V',W')\otimes V\otimes V'\big)\isoright \\
 \Hom\big(\bfcohom(V\otimes V',W\otimes W'),\bfcohom(V,W)\otimes\bfcohom(V',W')\big)
\end{multline*}
to the composition
\begin{multline*}
 W\otimes W'\xrightarrow{\coev\otimes\coev}
\bfcohom(V,W)\otimes V\otimes\bfcohom(V',W')\otimes V'\xrightarrow{\sigma^{(23)}} \\
 \bfcohom(V,W)\otimes\bfcohom(V',W')\otimes V\otimes V'.
\end{multline*}
This gives the morphism
$\varkappa=\varkappa^{V,V'}_{W,W'}=\theta^{-1} \big(\sigma^{(23)}\cdot(\coev\otimes\coev)\big) \colon\bfcohom(V\otimes V',W\otimes W')\to\bfcohom(V,W)\otimes\bfcohom(V',W').$
By applying duality principle to Prop.~\ref{ProppiNat} we conclude that these morphisms are also natural in $V,V',W,W'\in\bfP$.

The natural transformation $\varkappa$ is coassociative in the sense similar to Prop.~\ref{PropAsspi}, but with the opposite arrows: $(\varkappa\otimes\id)\cdot\varkappa=(\id\otimes\varkappa)\cdot\varkappa$.

In the case $W=V$ and $W'=V'$ it reads
\begin{align} \label{varkVV}
 \varkappa\colon\bfcoend(V\otimes V')\to\bfcoend(V)\otimes\bfcoend(V').
\end{align}

\begin{Prop} \label{Propvark}
 The morphisms~\eqref{varkVV} are morphisms of comonoids.
\end{Prop}

\noindent{\bf Proof.} It follows from Prop.~\ref{Proppi} by the duality principle. \qed

\bb{Tensor product of representations of bimonoids.} \label{bbTPBiMon}
A bimonoid in the symmetric monoidal category $\bfC$ is $\BB=(X,\mu_X,\eta_X,\Delta_X,\varepsilon_X)\in\Bimon(\bfC)$, where $\mu_X$, $\eta_X$, $\Delta_X$ and $\varepsilon_X$ are multiplication, unit, comultiplication and counit respectively on the object $X\in\bfC$. By considering $\BB$ as a monoid $\MM=(X,\mu_X,\eta_X)$ we obtain the category $\Lact(\BB)=\Lact(\MM)$. If $\bfC$ is closed relative to $\bfP$, then we have the subcategory $\Rep_\bfP(\BB)=\Rep_\bfP(\MM)\subset\Lact(\BB)$. Here we describe how the comonoid structure $(\Delta_X,\varepsilon_X)$ gives a monoidal structure on the category $\Rep_\bfP(\BB)$.

Let $\rho\colon X\to\bfend(V)$ and $\rho'\colon X\to\bfend(V')$ be two representations of the bimonoid $\BB$ (i.e. of the monoid $\MM$) on the objects $V,V'\in\bfP$. Define a morphism $\rho''\colon X\to\bfend(V\otimes V')$ as the composition
\begin{align} \label{rhoTP}
 X\xrightarrow{\Delta_X}X\otimes X\xrightarrow{\rho\otimes\rho'}\bfend(V)\otimes\bfend(V')
\xrightarrow{\pi}\bfend(V\otimes V').
\end{align}

\begin{Lem} \label{LemrhoTP}
 The morphism~\eqref{rhoTP} is a representation of $\BB$ on $V\otimes V'$.
\end{Lem}

\noindent{\bf Proof.} Since $\BB$ is a bimonoid, $\Delta_X$ is a morphism of monoids $\MM\to\MM\otimes\MM$. From Prop.~\ref{PropfgTP} we obtain the monoid morphism $\MM\otimes\MM\xrightarrow{\rho\otimes\rho'}\bfend(V)\otimes\bfend(V')$. Due to Prop.~\ref{Proppi} the formula~\eqref{rhoTP} gives a morphism $\MM\to\bfend(V\otimes V')$ in $\Mon(\bfC)$. \qed

We call the representation~\eqref{rhoTP} the {\it tensor product} of the representations $\rho$, $\rho'$.

\begin{Lem} \label{Lemrhotau}
 Let $\rho$, $\rho'$, $\tau$, $\tau'$ be representations of $\BB$ on $V,V',W,W'\in\bfP$. Denote the tensor products of $\rho$, $\rho'$ and of $\tau$, $\tau'$ by $\rho''$ and by $\tau''$ respectively. If morphisms $f\colon V\to W$ and $f'\colon V'\to W'$ are morphisms of representations, i.e. morphisms $f\colon (V,\rho)\to(W,\tau)$ and $f'\colon(V',\rho')\to(W',\tau')$ in $\Rep_\bfP(\BB)$, then $f\otimes f'\colon V\otimes V'\to W\otimes W'$ is a morphism $f\otimes f'\colon (V\otimes V',\rho'')\to(W\otimes W',\tau'')$ in $\Rep_\bfP(\BB)$.
\end{Lem}

\noindent{\bf Proof.} By definition~\ref{DefRep} the fact that $f$ and $f'$ are morphisms of representations means the equalities $\bfhom(\id,f)\cdot\rho=\bfhom(f,\id)\cdot\tau$ and $\bfhom(\id,f')\cdot\rho'=\bfhom(f',\id)\cdot\tau'$. By taking into account these equalities and the naturality of $\pi$ (Prop.~\ref{ProppiNat}) we obtain the commutative diagram
\begin{align*}
\xymatrix{
 X\otimes X\ar[rrrr]^{\rho\otimes\rho'}\ar[d]^{\tau\otimes\tau'} &&&&
[V,V]\otimes[V',V']\ar[d]^{\bfhom(\id,f)\otimes\bfhom(\id,f')}\ar[r]^{\pi}  &
 \bfend(V\otimes V')\ar[dd]^{\quad\bfhom(\id,f\otimes f')} \\
[W,W]\otimes[W',W']\ar[d]^{\pi}\ar[rrrr]^{\bfhom(f,\id)\otimes\bfhom(f',\id)} &&&&
[V,W]\otimes[V',W']\ar[dr]^{\pi} \\
\bfend(W\otimes W')\ar[rrrrr]^{\bfhom(f\otimes f',\id)} &&&&&
[V\otimes V',W\otimes W']
}
\end{align*}
Composition with $\Delta_X\colon X\to X\otimes X$ gives $\bfhom(\id,f\otimes f')\cdot\rho''=\bfhom(f\otimes f',\id)\cdot\tau''$. \qed

\begin{Th} \label{ThrhoTP}
 For any bimonoid $\BB$ in a symmetric monoidal category $\bfC=(\bfC,\otimes)$ which is closed relative to a monoidal subcategory $\bfP\subset\bfC$ the category $\Rep_\bfP(\BB)$ is equipped with the monoidal structure in the following way. The monoidal product of objects $(V,\rho)$ and $(V',\rho')$ is the object $(V'',\rho'')$, where $V''=V\otimes V'$ and $\rho''\colon\BB\to\bfend(V'')$ is defined by~\eqref{rhoTP}. The monoidal product of morphisms $f$ and $g$ in $\Rep_\bfP(\BB)$ coincides with their monoidal product in $\bfP$. The unit object of $\Rep_\bfP(\BB)$ is $(I_\bfC,\varepsilon_X)$. If $\BB$ is cocommutative (i.e. $\sigma\cdot\Delta_X=\Delta_X$), then the monoidal category $\Rep_\bfP(\BB)$ is symmetric and $\sigma_{(V,\rho),(V',\rho')}=\sigma_{V,V'}$.
\end{Th}

\noindent{\bf Proof.} The monoidal operation on objects and on morphisms of $\Rep_\bfP(\BB)$ is defined correctly due to Lemmas~\ref{LemrhoTP} and \ref{Lemrhotau} respectively. Let us check its associativity. Consider representations $\rho,\rho',\wt\rho$ of $\BB$ on $V,V',\wt V\in\bfP$ respectively and let $\rho''$ be the tensor product of $\rho$ and $\rho'$. The sequential monoidal products of objects $(V,\rho)$, $(V',\rho')$ and $(\wt V,\wt\rho)$ is given by the tensor product of the representations $\rho''$ and $\wt\rho$, this is a representation of $\BB$ on $V\otimes V'\otimes\wt V$ of the form
\begin{multline}
 \pi\cdot(\rho''\otimes\wt\rho)\cdot\Delta_X=
\pi\cdot(\id\otimes\wt\rho)\cdot(\pi\otimes\id)\cdot(\rho\otimes\rho'\otimes\id)\cdot(\Delta_X\otimes\id)\cdot\Delta_X= \\
\pi\cdot(\pi\otimes\id)\cdot(\rho\otimes\rho'\otimes\wt\rho)\cdot(\Delta_X\otimes\id)\cdot\Delta_X. \label{rhorhowt}
\end{multline}
By taking the tensor products in another order we obtain
\begin{multline}
\pi\cdot(\rho\otimes\wh\rho)\cdot\Delta_X=
\pi\cdot(\rho\otimes\id)\cdot(\id\otimes\pi)\cdot(\id\otimes\rho'\otimes\wt\rho)\cdot(\id\otimes\Delta_X)\cdot\Delta_X= \\
\pi\cdot(\id\otimes\pi)\cdot(\rho\otimes\rho'\otimes\wt\rho)\cdot(\id\otimes\Delta_X)\cdot\Delta_X, \label{rhorhowh}
\end{multline}
where $\wh\rho$ is the tensor product of $\rho'$ and $\wt\rho$. The representation~\eqref{rhorhowt} coincides with~\eqref{rhorhowh} by virtue of coassociativity of $\Delta_X$ and Prop.~\ref{PropAsspi}. The associativity of the monoidal operation on morphisms follows from the associativity of the monoidal product in $\bfP$.

By definition of bimonoid the counit $\varepsilon_X$ is in agreement with the monoid structure $(\mu_X,\eta_X)$, i.e. the diagrams
\begin{align*}
\xymatrix{
X\otimes X\ar[rrr]^{\mu_X}\ar[d]_{\varepsilon_X\otimes\varepsilon_X} &&& X\ar[d]^{\varepsilon_X} \\
I\otimes I\ar@{=}[r] &I\ar[rr]^{c_I=\id_I} && I
} \qquad\qquad
\xymatrix{
 I\ar[r]^{\eta_X}\ar[dr]_{u_I=\id_I} & X\ar[d]^{\varepsilon_X} \\
 & I
}
\end{align*}
commute. This means exactly that $\varepsilon_X\colon X\to I=\bfend(I)$ is a representation of $\BB$ on the object $I$. Let us check that $(I,\varepsilon_X)\in\Rep_\bfP(\BB)$ is a unit object. The tensor product of a representation $\rho\colon X\to\bfend(V)$ with the representation $\varepsilon_X\colon X\to I=\bfend(I)$ equals $\pi^{I,V}_{I,V}\cdot(\varepsilon_X\otimes\rho)\cdot\Delta_X=\pi^{I,V}_{I,V}\cdot(\id_I\otimes\rho)\cdot(\varepsilon_X\otimes\id_X)\cdot\Delta_X$. Due to $(\varepsilon_X\otimes\id_X)\cdot\Delta_X=\id_X$ and Prop.~\ref{ProppiII} it coincides with $\rho$, so that $(I,\varepsilon_X)\otimes(V,\rho)=(V,\rho)$. The identification $(V,\rho)\otimes(I,\varepsilon_X)=(V,\rho)$ is checked in the same way.

If $\BB$ is cocommutative, then due to the naturality of $\sigma$ and Lemma~\ref{Lempisigma} the diagram
\begin{align*}
\xymatrix{
X\ar[r]^{\Delta_X\quad}\ar[dr]_{\Delta_X} &
 X\otimes X\ar[r]^{\rho\otimes\rho'\qquad\quad}\ar[d]_{\sigma} &
 \bfend(V)\otimes\bfend(V')\ar[d]_{\sigma}\ar[r]^{\quad\pi} &
\bfend(V\otimes V')\ar[rrd]^{\bfhom(\id,\sigma_{V,V'})} \\
 & X\otimes X\ar[r]^{\rho'\otimes\rho\qquad\quad} &
\bfend(V')\otimes\bfend(V)\ar[r]_{\quad\pi} & 
\bfend(V'\otimes V)\ar[rr]_{\bfhom(\sigma_{V,V'},\id)\quad} && [V\otimes V',V'\otimes V]
}
\end{align*}
commutes, so $\sigma_{V,V'}\colon V\otimes V'\isoright V'\otimes V$ is a morphism $(V,\rho)\otimes(V',\rho')\to(V',\rho')\otimes(V,\rho)$ in $\Rep_\bfP(\BB)$. This is an isomorphism in $\Rep_\bfP(\BB)$, since its inverse $\sigma_{V',V}$ is a morphism $(V',\rho')\otimes(V,\rho)\to(V,\rho)\otimes(V',\rho')$. \qed

Let us draw the dual picture. The concept of bimonoid is self-dual. By considering $\BB$ as a comonoid $\OO=(X,\Delta_X,\varepsilon_X)$ we obtain the category $\Lcoact(\BB)=\Lcoact(\OO)$. If $\bfC$ is coclosed relative to $\bfP$, then we define $\Corep_\bfP(\BB)=\Corep_\bfP(\OO)\subset\Lcoact(\BB)$. This subcategory is a monoidal category whose monoidal structure is given by the morphisms $\mu_X,\eta_X$. To describe this structure we dualise Lemma~\ref{LemrhoTP} and Theorem~\ref{ThrhoTP}.

\begin{Lem}
Let $\omega\colon\bfcoend(V)\to X$ and $\omega'\colon\bfcoend(V')\to X$ be corepresentations of the bimonoid $\BB$ (as a comonoid $\OO$) on the objects $V\in\bfP$ and $V'\in\bfP$ respectively. Then the morphism $\omega''\colon\bfcoend(V\otimes V')\to X$ defined as the composition
\begin{align} \label{omegaTP}
 \bfcoend(V\otimes V')\xrightarrow{\varkappa}
\bfcoend(V)\otimes\bfcoend(V')\xrightarrow{\omega\otimes\omega'}
X\otimes X\xrightarrow{\mu_X}X
\end{align}
is a corepresentation of $\BB$ on $V\otimes V'$.
\end{Lem}

The corepresentation~\eqref{omegaTP} is called the {\it tensor product} of the corepresentations $\omega,\omega'$.

\begin{Th} \label{ThomegaTP}
 For any bimonoid $\BB$ in a symmetric monoidal category $\bfC=(\bfC,\otimes)$ which is coclosed relative to a monoidal subcategory $\bfP\subset\bfC$ the category $\Corep_\bfP(\BB)$ is equipped with the monoidal structure. The monoidal product on objects has the form $(V,\omega)\otimes(V',\omega')=(V'',\omega'')$, where $V''=V\otimes V'$ and $\omega''\colon\bfcoend(V'')\to\BB$ is the corepresentation~\eqref{omegaTP}. On morphisms it is given by the monoidal product $\otimes$ of $\bfC$. The unit object of $\Corep_\bfP(\BB)$ is $(I_\bfC,\eta_X)$. If $\BB$ is commutative (i.e. $\mu_X\cdot\sigma=\mu_X$), then the monoidal category $\Corep_\bfP(\BB)$ is symmetric and $\sigma_{(V,\omega),(V',\omega')}=\sigma_{V,V'}$.
\end{Th}

\bbo{Tensor product of actions.}
Since Theorem~\ref{ThRepLact} identifies the representations with actions, Lemma~\ref{LemrhoTP} allows to define tensor product of actions of the bimonoid $\BB$ on objects of $\bfP$ in the relatively closed case. This operation can be written explicitly via the coproduct $\Delta_X$. In the relatively coclosed case we have tensor product of coactions of the bimonoid $\BB$ on objects of $\bfP$, which can be written via $\mu_X$.

\begin{Prop} \label{PropAcTP}
 Let $\rho\colon X\to\bfend(V)$ and $\rho'\colon X\to\bfend(V')$ be representations of $\BB$ on $V,V'\in\bfP$ and $\rho''\colon X\to\bfend(V'')$ be their tensor product~\eqref{rhoTP}, where $V''=V\otimes V'$. Consider the corresponding actions $a=\theta(\rho)\colon X\otimes V\to V$, $a'=\theta(\rho')\colon X\otimes V'\to V'$ and $a''=\theta(\rho'')\colon X\otimes V''\to V''$. Then $a''$ coincides with the composition
\begin{align} \label{AcTP}
 X\otimes V\otimes V'\xrightarrow{\Delta_X\otimes\id\otimes\id}
 X\otimes X\otimes V\otimes V'\xrightarrow{\sigma^{(23)}}
 X\otimes V\otimes X\otimes V'\xrightarrow{a\otimes a'} V\otimes V'.
\end{align}

Dually, let $\omega\colon\bfcoend(V)\to X$ and $\omega'\colon\bfcoend(V')\to X$ be corepresentations of the bimonoid $\BB$ on the objects $V,V'$. Let $\delta=\vartheta(\omega)\colon V\to X\otimes V$ and $\delta'=\vartheta(\omega')\colon V'\to X\otimes V'$. Then the coaction $\delta''=\vartheta(\omega'')\colon V''\to X\otimes V''$ corresponding to the tensor product~\eqref{omegaTP} coincides with
\begin{align} \label{CoacTP}
V\otimes V' \xrightarrow{\delta\otimes\delta'} X\otimes V\otimes X\otimes V \xrightarrow{\sigma^{(23)}} X\otimes X\otimes V\otimes V'\xrightarrow{\mu_X\otimes\id\otimes\id}
 X\otimes V\otimes V'.
\end{align}
\end{Prop}

\noindent{\bf Proof.} By sequentially applying~\eqref{thetafg}, definition of $\pi$, naturality of $\sigma$ and the formula~\eqref{thetaEval} we derive
\begin{multline} \label{PropAcTPeq1}
 \theta\big(\pi\cdot(\rho\otimes\rho')\big)=
\theta(\pi)\cdot(\rho\otimes\rho'\otimes\id_{V\otimes V'})=
 (\eval\otimes\eval)\cdot\sigma^{(23)}\cdot (\rho\otimes\rho'\otimes\id_V\otimes\id_{V'})= \\
 (\eval\otimes\eval)\cdot (\rho\otimes\id_V\otimes\rho'\otimes\id_{V'})\cdot\sigma^{(23)}=
 (\theta(\rho)\otimes\theta(\rho'))\cdot\sigma^{(23)}=(a\otimes a')\cdot\sigma^{(23)}.
\end{multline}
By using~\eqref{thetafg} again we obtain
\begin{align} \label{PropAcTPeq2}
 a''=\theta(\rho'')=\theta\big(\pi\cdot(\rho\otimes\rho')\cdot\Delta_X\big)=
\theta\big(\pi\cdot(\rho\otimes\rho')\big)\cdot(\Delta_X\otimes\id_{V\otimes V'}),
\end{align}
Substitution of~\eqref{PropAcTPeq1} to \eqref{PropAcTPeq2} yields~\eqref{AcTP}. \qed

\begin{Rem} \label{RemRepLact} \normalfont
 If $\bfC=\bfP$ is a closed monoidal category, then the equivalence~\eqref{RepLact} induces a monoidal structure from the category $\Rep_\bfC(\BB)$ to $\Lact(\BB)$. It is claimed in~\cite{P} (without proof) that for arbitrary symmetric monoidal category $\bfC$ and a bimonoid $\BB\in\Bimon(\bfC)$ the category $\Lact(\BB)$ is monoidal with respect to the monoidal product defined by the formula~\eqref{AcTP}. Dually, $\Lcoact(\BB)$ is a monoidal category with the structure defined by~\eqref{CoacTP}. By this logic, Prop.~\ref{PropAcTP} means that the embeddings~\eqref{RepLact} and \eqref{CorepLcoact} are strict monoidal functors.
\end{Rem}

\subsection{Translation of representations by a functor}
\label{secTrans}

If a functor between two relatively (co)closed monoidal categories respects their structures in an appropriate way, then this functor induces functors between the categories of (co)monoids, bimonoids, (co)actions and (co)representations. We describe the functor which translates (co)representations in details by following ideas of~\cite[\S~3.3]{Sqrt} (we adapt them to the case of relatively (co)closed monoidal categories). As above we do not use the actions and their connection with representations.
For the case of bimonoid we study the translation of the tensor product of (co)representations and get a monoidal structure of the corresponding induced functor.

\bb{Translation of internal hom.} \label{bbTrhom}
Let $\bfC=(\bfC,\otimes)$ and $\bfD=(\bfD,\odot)$ be symmetric monoidal categories. Let $\bfP\subset\bfC$ and $\bfQ\subset\bfD$ be their full subcategories. Suppose that $\bfC$ and $\bfD$ are relatively closed with the parametrising categories $\bfP$ and $\bfQ$ respectively.

Let $F=(F,\varphi,\phi)\colon(\bfC,\otimes)\to(\bfD,\odot)$ be a lax monoidal functor: it consists of a functor $F\colon\bfC\to\bfD$, a morphism $\varphi\colon I_\bfD\to FI_\bfC$ and a natural transformation $\phi$ with the components $\phi_{X,Y}\colon FX\odot FY\to F(X\otimes Y)$ satisfying necessary conditions. 

Suppose $F(\bfP)\subset\bfQ$. Then we can define the hom-object $\bfhom(FV,FW)$ in $\bfD$ for any $V,W\in\bfP$ and relate it with the image $F\big(\bfhom(V,W)\big)$ by a natural morphism as follows (see~\cite[Prop.~6.4.5]{Bor2} for the closed case).

Let us apply the right adjunction bijection
\begin{align*}
 \theta^{-1}\colon\Hom\big(F\big(\bfhom(V,W)\big)\odot FV,FW\big)\isoright
\Hom\big(F\big(\bfhom(V,W)\big),\bfhom(FV,FW)\big)
\end{align*}
to the composition
\begin{align} \label{phiFeval}
 F\big(\bfhom(V,W)\big)\odot FV\xrightarrow{\phi}
 F\big(\bfhom(V,W)\otimes V\big)\xrightarrow{F(\eval_{V,W})}FW
\end{align}
and denote the result by $\Phi_{V,W}=\theta^{-1}\big(F(\eval_{V,W})\cdot\phi_{\bfhom(V,W),V}\big)$. This is a unique morphism
\begin{align} \label{PhiDef}
 \Phi_{V,W}\colon F\big(\bfhom(V,W)\big)\to\bfhom(FV,FW)
\end{align}
such that the diagram
\begin{align} \label{PhiDiag}
\xymatrix{
 F\big(\bfhom(V,W)\big)\odot FV\ar[rr]^{\Phi_{V,W}\odot\id_{FV}}\ar[d]_\phi &&
 \bfhom(FV,FW)\odot FV\ar[d]^{\eval_{FV,FW}} \\
 F\big(\bfhom(V,W)\otimes FV\big)\ar[rr]^{F(\eval_{VW})} && FW
}
\end{align}
commutes.

\begin{Prop} \label{PropPhiphi}
 For any objects $X\in\bfC$, $V,W\in\bfP$ and morphism $f\colon X\to\bfhom(V,W)$ in $\bfC$ we have the formula
\begin{align} \label{Phif}
 \theta\big(\Phi_{V,W}\cdot F(f)\big)=F\big(\theta(f)\big)\cdot\phi_{X,V}.
\end{align}
\end{Prop}

\noindent{\bf Proof.} The formula~\eqref{Phif} follows from the commutativity of the diagram~\eqref{PhiDiag} and the naturality of $\phi_{X,Y}$ in the first argument (see the proof of~\cite[Prop.~3.9]{Sqrt}). \qed

In the closed case the naturality of $\Phi$ is due to the point~(1) of~\cite[Prop.~6.4.5]{Bor2}, which was claimed without a (detailed) proof.
For the relative case we deduce the naturality of $\Phi$ by using Prop.~\ref{PropPhiphi} and the naturality of $\phi_{X,Y}$ in $Y\in\bfP$.

\begin{Prop} \label{PropPhiNat}
 The morphisms~\eqref{PhiDef} are natural in $V,W\in\bfP$, i.e. for any morphisms $f\colon V'\to V$ and $g\colon W\to W'$ in $\bfP$ we have the commutative diagram
\begin{align} \label{PhiNat}
\xymatrix{
 F\big(\bfhom(V,W)\big)\ar[r]^{\Phi_{V,W}}\ar[d]_{F\big(\bfhom(f,g)\big)} &
 \bfhom(FV,FW)\ar[d]^{\bfhom(Ff,Fg)} \\
 F\big(\bfhom(V',W')\big)\ar[r]^{\Phi_{V',W'}} & \bfhom(FV',FW')
}
\end{align}
\end{Prop}

\noindent{\bf Proof.} By using the formulae~\eqref{Phif}, \eqref{thetahom} and \eqref{thetafg} we obtain
\begin{align*}
 &\theta\Big(\Phi_{V',W'}\cdot F\big(\bfhom(f,g)\big)\Big)=
 F\Big(\theta\big(\bfhom(f,g)\big)\Big)\cdot\phi_{X,V'}=Fg\cdot F(\eval)\cdot F(\id\otimes f)\cdot\phi_{X,V}, \\
 &\theta\big(\bfhom(Ff,Fg)\cdot\Phi_{V,W}\big)=
 \theta\big(\bfhom(Ff,Fg)\big)\cdot(\Phi_{V,W}\odot\id)=
 Fg\cdot\eval\cdot(\Phi_{V,W}\odot Ff),
\end{align*}
where $X=\bfhom(V,W)$. The right hand sides of these equalities coincide with each other due to the naturality of $\phi$ and the commutativity of~\eqref{PhiDiag}:
\begin{align*}
 &F(\id\otimes f)\cdot\phi_{X,V}=\phi_{X,V'}\cdot(\id\odot Ff),
 &&\eval\cdot(\Phi_{V,W}\odot Ff)=F(\eval)\cdot\phi_{X,V'}\cdot(\id\odot Ff).
\end{align*}
Then the bijectivity of $\theta$ implies that the diagram~\eqref{PhiNat} commutes. \qed

\begin{Rem} \normalfont
 Even if the functor $F$ is strong monoidal, it is not guaranteed that the natural transformation $\Phi$ is an isomorphism. A counterexample is the monoidal embedding $\bfP\hookrightarrow\bfC$ described in Remark~\ref{RemCountEx}: the formulae~\eqref{cohomxyP} and \eqref{cohomxyC} imply that the components $\Phi_{x,y}$ are not isomorphisms if $0<x<y$.
\end{Rem}

\bbo{Translation of representations.} \label{bbFrep}
For a monoid $\MM=(X,\mu_X,\eta_X)\in\Mon(\bfC)$ define morphisms $\mu_{FX}\colon FX\odot FX\to FX$ and $\eta_{FX}\colon I_\bfD\to FX$ as the compositions
\begin{align} \label{muFX}
 &FX\odot FX\xrightarrow{\phi_{X,X}}F(X\otimes X)\xrightarrow{F\mu_X} FX &&\text{and} && I_\bfD\xrightarrow{\varphi}FI_\bfC\xrightarrow{F\eta_X} FX
\end{align}
respectively.
Any lax monoidal functor $F\colon\bfC\to\bfD$ induces the functor 
\begin{align} \label{MonF}
 \Mon(F)\colon\Mon(\bfC)\to\Mon(\bfD)
\end{align}
which translates $\MM$ to the monoid  $\Mon(F)\MM=(FX,\mu_{FX},\eta_{FX})\in\Mon(\bfD)$.
If $F$ is faithful, then so is $\Mon(F)$. If $F$ is strong monoidal and fully faithful, then $\Mon(F)$ is fully faithful as well.

For any object $V\in\bfP$ the monoid $\bfend(V)$ is translated to the monoid $\Mon(F)\big(\bfend(V)\big)$, which is the object $F\big(\bfend(V)\big)\in\bfD$ with the multiplication $F(c_{V})\cdot\phi_{\bfend(V),\bfend(V)}$ and unit $F(u_V)\cdot\varphi$.

Substitution $W=V\in\bfP$ to~\eqref{PhiDef} gives the morphism
\begin{align} \label{PhiVV}
 \Phi_{V,V}\colon F\big(\bfend(V)\big)\to\bfend(FV)
\end{align}
in $\bfD$, where the end-object $\bfend(FV)\in\bfD$ is defined due to $FV\in\bfQ$.
The first sentence of the following statement is a particular case of the points~(2) and (3) of~\cite[Prop.~6.4.5]{Bor2} generalised to the relative case.

\begin{Prop} \label{PropFend}
 The morphism~\eqref{PhiVV} is a morphism in $\Mon(\bfD)$ from the monoid $\Mon(F)\big(\bfend(V)\big)$ to $\bfend(FV)=\big(\bfend(FV),c_{FV},u_{FV}\big)$. In particular, if~\eqref{PhiVV} is an isomorphism in $\bfD$, then it gives an isomorphism of monoids $\Mon(F)\big(\bfend(V)\big)\cong\bfend(FV)$.
\end{Prop}

\noindent{\bf Proof.} The first sentence can be reformulated as the commutativity of the diagrams
\begin{gather*}
\xymatrix{
 F\big(\bfend(V)\big)\odot F\big(\bfend(V)\big)\ar[r]^{\quad\phi}\ar[d]^{\Phi_{V,V}\odot\Phi_{V,V}} & F\big(\bfend(V)\otimes\bfend(V)\big)\ar[r]^{\quad\qquad F(c_V)} & F\big(\bfend(V)\big)\ar[d]^{\Phi_{V,V}} \\
 \bfend(FV)\odot\bfend(FV)\ar[rr]^{c_{FV}}&& \bfend(FV)
} \\
\xymatrix{
I_\bfD\ar[r]^{\varphi}\ar[drr]_{u_{FV}} & FI_\bfC\ar[r]^{F(u_V)\qquad} & F\big(\bfend(V)\big)\ar[d]^{\Phi_{V,V}} \\
 && \bfend(FV)
}
\end{gather*}
By using the formulae~\eqref{thetafg}, \eqref{Phif},~$\phi\cdot(\varphi\odot\id_{FV})=\id_{FV}$ and $\theta(u_V)=\id_V$ we obtain
\begin{multline*}
 \theta\big(\Phi_{V,V}\cdot F(u_V)\cdot\varphi\big)=
\theta\big(\Phi_{V,V}\cdot F(u_V)\big)\cdot(\varphi\odot\id_{FV})= F\big(\theta(u_V)\big)\big)\cdot\phi\cdot(\varphi\odot\id_{FV})= \\
\id_{FV}=\theta(u_{FV}).
\end{multline*}
This implies the commutativity of the second diagram. The commutativity of the first one is checked in the same manner:
\begin{multline*}
 \theta\big(\Phi_{V,V}\cdot F(c_V)\cdot\phi\big)=
 \theta\big(\Phi_{V,V}\cdot F(c_V)\big)\cdot(\phi\odot\id)=
 F\big(\theta(c_V)\big)\cdot\phi\cdot(\phi\odot\id)= \\
 F(\eval)\cdot F(\id\otimes\eval)\cdot\phi\cdot(\id\odot\phi)=
 F(\eval)\cdot\phi\cdot\big(\id\odot F(\eval)\big)\cdot(\id\odot\phi)= \\
 \theta(\Phi_{V,V})\cdot\big(\id\odot \theta(\Phi_{V,V})\big)=
 \eval\cdot(\Phi_{V,V}\odot\id)\cdot(\id\odot\eval)\cdot(\id\odot\Phi_{V,V}\odot\id)= \\
 \eval\cdot(\id\odot\eval)\cdot(\Phi_{V,V}\odot\id)\cdot(\id\odot\Phi_{V,V}\odot\id)=
 \theta(c_{FV})\cdot(\Phi_{V,V}\odot\Phi_{V,V}\odot\id)=
 \theta\big(c_{FV}\cdot(\Phi_{V,V}\odot\Phi_{V,V})\big),
\end{multline*}
where we also used $\theta(c_V)=\eval\cdot(\id\otimes\eval)$, $\phi\cdot(\phi\odot\id)=\phi\cdot(\id\odot\phi)$, the formula~\eqref{thetaEval} and the naturality of $\phi$. The second sentence of Prop.~\ref{PropFend} follows from the fact that the forgetful functor $\Mon(\bfC)\to\bfC$ reflects isomorphisms (see e.g.~\cite[Prop.~2.1]{Sqrt}). \qed

Let us describe how the lax monoidal functor $F$ translates representations of the monoid $\MM$ on $V$ to representations of the monoid $\wt\MM=\Mon(F)\MM$ on $FV$.

\begin{Th} \label{ThFRep}
Let $\bfC=(\bfC,\otimes)$ and $\bfD=(\bfD,\odot)$ be closed relative to $\bfP$ and $\bfQ$. Any lax monoidal functor $F\colon(\bfC,\otimes)\to(\bfD,\odot)$ satisfying $F(\bfP)\subset\bfQ$ induces the functor
\begin{align} \label{FRep}
 F_{\bfP,\MM}\colon\Rep_\bfP(\MM)\to\Rep_\bfQ(\wt\MM),
\end{align}
which maps $(V,\rho)\in\Rep_\bfP(\MM)$ to $(FV,\wt\rho)$, where $\wt\rho\colon\wt\MM\to\bfend(FV)$ is the composition
\begin{align} \label{FrhoPhi}
 FX\xrightarrow{F\rho}F\big(\bfend(V)\big)\xrightarrow{\Phi_{V,V}}\bfend(FV).
\end{align}
The image of a morphism $f\colon(V,\rho)\to(W,\tau)$ under the functor $F_{\bfP,\MM}$ is the morphism $Ff\colon(FV,\wt\rho)\to(FW,\wt\tau)$, where $(W,\tau)\in\Rep_\bfP(\MM)$ and $\wt\tau=\Phi_{W,W}\cdot F\tau$.
\end{Th}

\noindent{\bf Proof.} The functor $\Mon(F)$ maps the monoid morphism $\rho\colon\MM\to\bfend(V)$ to the monoid morphism $F\rho\colon\wt\MM\to\Mon(F)\big(\bfend(V)\big)$. Hence Prop.~\ref{PropFend} implies that~\eqref{FrhoPhi} is a monoid morphism $\wt\MM\to\bfend(FV)$. The commutativity of the diagram
\begin{align}
\xymatrix{
FX\ar[rr]^{F\rho}\ar[d]_{F\tau} &&
 F\big([V,V]\big)\ar[d]^{F\big(\bfhom(\id,f)\big)}\ar[rr]^{\Phi_{V,V}} &&
 [FV,FV]\ar[dd]^{\bfhom(\id,Ff)} \\
F\big([W,W]\big)\ar[rr]^{F\big(\bfhom(f,\id)\big)}\ar[d]^{\Phi_{W,W}} &&
 F\big([V,W]\big)\ar[drr]^{\Phi_{V,W}} \\
[FW,FW]\ar[rrrr]^{\bfhom(Ff,\id)} &&&& [FV,FW]
}
\end{align}
follows from the commutativity of~\eqref{rhof} and the naturality of $\Phi$ (Prop.~\ref{PropPhiNat}). \qed

\bb{Translation of corepresentations.} \label{bbTrCorep}
Let $\bfC=(\bfC,\otimes)$ and $\bfD=(\bfD,\odot)$ be symmetric relatively coclosed monoidal categories with the parametrising subcategories $\bfP$ and $\bfQ$ respectively. Let $F=(F,\varphi,\phi)\colon(\bfC,\otimes)\to(\bfD,\odot)$ be a colax monoidal functor satisfying $F(\bfP)\subset\bfQ$.

By considering the relatively closed monoidal categories $(\bfC^\op,\otimes)$ and $(\bfD^\op,\odot)$ with the parametrising subcategories $\bfP^\op\subset\bfC^\op$ and $\bfQ^\op\subset\bfD^\op$ and by using the identification functors~\eqref{MonComon} and \eqref{RepCorep} we can apply the results of pp.~\ref{bbTrhom} and \ref{bbFrep} to the lax monoidal functor $F^\op\colon(\bfC^\op,\otimes)\to(\bfD^\op,\odot)$.

The colax monoidal functor $F$ induces the functor
\begin{align} \label{ComonF}
 \Comon(F)\colon\Comon(\bfC)\to\Comon(\bfD)
\end{align}
which translates a comonoid $(X,\Delta_X,\varepsilon_X)\in\Comon(\bfC)$ to $(FX,\Delta_{FX},\varepsilon_{FX})\in\Comon(\bfD)$, where
\begin{align} \label{DeltaFX}
 &\Delta_{FX}=\phi_{X,X}\cdot F\Delta_X\colon FX\to FX\odot FX, &&\varepsilon_{FX}=\varphi\cdot F\varepsilon_X\colon FX\to I_\bfD.
\end{align}
The monoidal structure of $F$ defines the natural transformation $\Phi$ with components
\begin{align}
 \Phi_{V,W}=\vartheta^{-1}\big(\phi_{\bfcohom(V,W),V}\cdot F(\coev_{V,W})\big)\colon
\bfcohom(FV,FW)\to F\big(\bfcohom(V,W)\big).
\end{align}
Some of these components give monoid morphisms:
\begin{align} \label{PhiVVco}
 \Phi_{V,V}\colon\bfcoend(FV)\to\Comon(F)\big(\bfcoend(V)\big).
\end{align}
If $\Phi_{V,V}\colon\bfcoend(V)\to F\big(\bfcoend(V)\big)$ is an isomorphism in $\bfD$, then~\eqref{PhiVVco} is an isomorphism of comonoids $\bfcoend(FV)\cong\Comon(F)\big(\bfcoend(V)\big)$.

\begin{Th} \label{ThCorep}
Let $\bfC$ and $\bfD$ be coclosed relative to $\bfP$ and $\bfQ$. Let $F\colon(\bfC,\otimes)\to(\bfD,\odot)$ be a colax monoidal functor such that $F(\bfP)\subset\bfQ$. Let $\OO\in\Comon(\bfC)$ and $\wt\OO=\Comon(F)\OO$. Then $F$ induces the functor
\begin{align} \label{FCorep}
 F^\OO_{\bfP}\colon\Corep_\bfP(\OO)\to\Corep_\bfQ(\wt\OO),
\end{align}
which maps $(V,\omega)\in\Corep_\bfP(\OO)$ to $(FV,\wt\omega)$, where $\wt\omega$ is the composition
\begin{align} \label{FomegaPhi}
 \bfcoend(FV)\xrightarrow{\Phi_{V,V}}F\big(\bfcoend(V)\big)\xrightarrow{F\omega} FX.
\end{align}
A morphism $f\colon(V,\omega)\to(W,\nu)$ is mapped to the morphism $Ff\colon(FV,\wt\omega)\to(FW,\wt\nu)$, where $W\in\bfP$, $\nu$ is a corepresentation of $\OO$ on $W$ and $\wt\nu=F\nu\cdot\Phi_{W,W}$.
\end{Th}

\bbo{Translation of actions.}
For any monoid $\MM\in\Mon(\bfC)$ a lax monoidal functor $F=(F,\varphi,\phi)\colon\bfC\to\bfD$ induces a functor
\begin{align} \label{FLact}
 F_\MM\colon\Lact(\MM)\to\Lact(\wt\MM),
\end{align}
where $\wt\MM=\Mon(F)\MM$. Each action $a\colon X\otimes V\to V$ of $\MM$ is translated to the action $\wt a\colon FX\odot FV\to FV$ of $\wt\MM$ defined as the composition
\begin{align}
 FX\odot FV\xrightarrow{\phi_{X,V}}F(X\otimes V)\xrightarrow{Fa}FV,
\end{align}
i.e. $(V,a)\mapsto(FV,\wt a)$ and $f\mapsto Ff$ (see e.g.~\cite[\S~2.4.11]{Sqrt} for details).

Dually, a colax monoidal functor $F=(F,\varphi,\phi)\colon\bfC\to\bfD$ induces the functor
\begin{align} \label{FLcoact}
 &F^\OO\colon\Lcoact(\OO)\to\Lcoact(\wt\OO)
\end{align}
for comonoids $\OO\in\Comon(\bfC)$ and $\wt\OO=\Comon(F)\OO\in\Comon(\bfD)$. It translates a coaction $\delta\colon V\to X\otimes V$ to the coaction
\begin{align}
FV\xrightarrow{F\delta}F(X\otimes V)\xrightarrow{\phi_{X,V}} FX\odot FV.
\end{align}

In the situation defined in p.~\ref{bbTrhom} we have subcategories $\Rep_\bfP(\MM)\subset\Lact(\MM)$ and $\Rep_\bfQ(\wt\MM)\subset\Lact(\wt\MM)$. Since $F(\bfP)\subset\bfQ$, the restriction of the functor~\eqref{FLact} gives a functor $\Rep_\bfP(\MM)\to\Rep_\bfQ(\wt\MM)$.

\begin{Prop} \label{PropRepLact}
 The restrictions of the functors~\eqref{FLact} and \eqref{FLcoact} coincide with~\eqref{FRep} and \eqref{FCorep} respectively. In other words the functor diagrams
\begin{align}
\xymatrix{
\Rep_\bfP(\MM)\ar[r]^{F_{\bfP,\MM}}\ar@{^{(}->}[d] &
 \Rep_\bfQ(\wt\MM)\ar@{^{(}->}[d] \\
\Lact(\MM)\ar[r]^{F_\MM} & \Lact(\wt\MM) 
} \qquad\qquad
\xymatrix{
\Corep_\bfP(\OO)\ar[r]^{F^\OO_\bfP}\ar@{^{(}->}[d] &
 \Corep_\bfQ(\wt\OO)\ar@{^{(}->}[d] \\
\Lcoact(\OO)\ar[r]^{F^\OO} & \Lcoact(\wt\OO) 
}
\end{align}
commute, where the vertical arrows `$\hookrightarrow$' are the embedding functors~\eqref{RepLact}, \eqref{CorepLcoact}.
\end{Prop}

\noindent{\bf Proof.} Repeats the proof of~\cite[Prop.~3.11]{Sqrt}. \qed

\begin{Rem} \normalfont \label{RemFend}
Theorems~\ref{ThFRep} and \ref{ThCorep} follows from Prop.~\ref{PropRepLact} and the construction of the functors~\eqref{FLact} and \eqref{FLcoact} (see~\cite[\S~3.3.3]{Sqrt}).
By taking the representation~\eqref{repId} as $\rho$ in Theorem~\ref{ThFRep} we conclude that $\Phi_{V,V}$ is a monoid morphism, so Prop.~\ref{PropFend} is also a consequence the construction of~\eqref{FLact}.
\end{Rem}

\bbo{Translation under a contravariant functor.}
We call a contravariant functor $F\colon\bfC\to\bfD$ {\it lax/colax/strong monoidal} from $(\bfC,\otimes)$ to $(\bfD,\odot)$ iff the corresponding covariant functor $\bar F\colon\bfC^\op\to\bfD$ is equipped with the structure of a lax/colax/strong monoidal functor $(\bar F,\varphi,\phi)\colon(\bfC^\op,\otimes)\to(\bfD,\odot)$. Equivalently, it means that the functor $\bar F^\op\colon\bfC\to\bfD^\op$ is a colax/lax/strong monoidal. Let both monoidal categories $\bfC=(\bfC,\otimes)$ and $\bfD=(\bfD,\odot)$ be symmetric. The monoidal functor $F$ is called {\it symmetric} iff $\bar F\colon\bfC^\op\to\bfD$ (or, equivalently, $\bar F^\op\colon\bfC\to\bfD^\op$) is symmetric.

The results written for covariant monoidal functors can be reformulated for contravariant monoidal functors. We consider only colax monoidal functors applying the results of pp.~\ref{bbTrhom} and \ref{bbFrep} to the categories $\bfC$ and $\bfD^\op$ or, equivalently, the results of p.~\ref{bbTrCorep} to $\bfC^\op$ and $\bfD$.
The case of contravariant lax monoidal functor is completely analogous, while a contravariant strong monoidal functor is a particular case of a contravariant colax monoidal functor.

Suppose that the symmetric monoidal categories $\bfC$ and $\bfD$ are relatively closed and coclosed respectively. Let $F=(F,\varphi,\phi)\colon\bfC\to\bfD$ be a contravariant colax monoidal functor, where $\varphi\colon FI_\bfC\to I_\bfD$ and $\phi_{X,Y}\colon F(X\otimes Y)\to FX\odot FY$ are structure morphisms for the colax monoidal functor $\bar F\colon\bfC^\op\to\bfD$. The contravariant functor $F$ induces a contravariant functor 
\begin{align} \label{ComonFCon}
 \Comon(F)\colon\Mon(\bfC)\to\Comon(\bfD)
\end{align}
corresponding to the composition of the identification~\eqref{MonComon} with the functor~\eqref{ComonF} for covariant functor $\bar F$.
A monoid $\MM=(X,\mu_X,\eta_X)\in\Mon(\bfC)$ is translated under~\eqref{ComonFCon} to the comonoid $\Comon(F)\MM=(FX,\Delta_{FX},\varepsilon_{FX})\in\Comon(\bfD)$, where
 $\Delta_{FX}=\phi_{X,X}\cdot F\mu_X$, $\varepsilon_{FX}=\varphi\cdot F\eta_X$.

The natural transformation $\Phi$ for the contravariant case has the components
\begin{align}
 \Phi_{V,W}=\vartheta^{-1}\big(\phi_{\bfhom(V,W),V}\cdot F(\eval_{V,W})\big)\colon
\bfcohom(FV,FW)\to F\big(\bfhom(V,W)\big),
\end{align}
$V,W,\in\bfP$.
It gives the morphisms of comonoids  $\Phi_{V,V}\colon\bfcoend(FV)\to\Comon(F)\big(\bfend(V)\big)$. If $\Phi_{V,V}$ is an isomorphism in $\bfD$, then it is an isomorphism in $\Comon(\bfD)$.

Let $\MM\in\Mon(\bfC)$ and $\OO\in\Comon(\bfC^\op)$ be a monoid and a comonoid related via~\eqref{MonComon}.
Theorems~\ref{ThFRep} and~\ref{ThCorep} for $\bar F^\op\colon\bfC\to\bfD^\op$ and $\bar F\colon\bfC^\op\to\bfD$ give the covariant functors $(\bar F^\op)_{\bfP,\MM}\colon\Rep_\bfP(\MM)\to\Rep_{\bfQ^\op}(\wt\MM)$ and $(\bar F)^\OO_{\bfP^\op}\colon\Corep_{\bfP^\op}(\OO)\to\Corep_{\bfQ}(\wt\OO)$ respectively, where $\wt\MM=\Mon(\bar F^\op)\MM\in\Mon(\bfD^\op)$, $\wt\OO=\Comon(\bar F)\OO\in\Comon(\bfD)$. By using~\eqref{RepCorep} we obtain the following statement.

\begin{Th} \label{ThRepCon}
Let the monoidal categories $\bfC=(\bfC,\otimes)$ and $\bfD=(\bfD,\odot)$ be closed and coclosed relative to $\bfP$ and $\bfQ$ respectively. Let $\MM=(X,\mu_X,\eta_X)\in\Mon(\bfC)$. Any contravariant colax monoidal functor $F\colon(\bfC,\otimes)\to(\bfD,\odot)$ such that $F(\bfP)\subset\bfQ$ induces the contravariant functor
\begin{align} \label{FRepCon}
 F^\MM_{\bfP}\colon\Rep_\bfP(\MM)\to\Corep_\bfQ(\wt\OO),
\end{align}
where $\wt\OO=\Comon(F)\MM$. It maps an object $(V,\rho)\in\Rep_\bfP(\MM)$ to $(FV,\wt\omega)\in\Corep_\bfQ(\wt\OO)$, where $\wt\omega\colon\bfcoend(FV)\to\wt\OO$ is the corepresentation defined as
\begin{align} \label{rhoomegaCon}
 \bfcoend(FV)\xrightarrow{\Phi_{V,V}}F\big(\bfend(V)\big)\xrightarrow{F\rho} FX.
\end{align}
A morphism $f\colon(V,\rho)\to(W,\tau)$ is mapped to the morphism $Ff\colon(FV,\wt\omega)\to(FW,\wt\nu)$, where $(W,\tau)\in\Rep_\bfP(\MM)$ and $\wt\nu=F\tau\cdot\Phi_{W,W}$.
\end{Th}

\bbo{Translation of tensor product of monoids.}
We denote the symmetric structure of the both symmetric monoidal categories  $\bfC=(\bfC,\otimes)$ and $\bfD=(\bfD,\odot)$ by the same letter $\sigma$.

\begin{Lem} \label{Lemsigmaphi}
 Let $F=(F,\varphi,\phi)\colon(\bfC,\otimes)\to(\bfD,\odot)$ be a lax monoidal functor. Then for any objects $X,Y,Z,W\in\bfC$ the diagram
\begin{align} \label{phiphiDiag}
\xymatrix{
FX\odot FY\odot FZ\odot FW\ar[rr]^{\id\odot\phi_{Y,Z}\odot\id}\ar[d]^{\phi_{X,Y}\odot\phi_{Z,W}} &&
FX\odot F(Y\otimes Z)\odot FW\ar[d]^{\phi\cdot(\phi\odot\id)=\phi\cdot(\id\odot\phi)} \\
F(X\otimes Y)\odot F(Z\otimes W)\ar[rr]^{\phi_{X\otimes Y,Z\otimes W}} &&
F(X\otimes Y\otimes Z\otimes W)
}
\end{align}
is commutative. If $F$ is symmetric, then we also have the following commutative diagram:
\begin{align*}
\xymatrix{
FX\odot FY\odot FZ\odot FW\ar[r]^{\sigma^{(23)}}\ar[d]^{\phi\odot\phi} &
FX\odot FZ\odot FY\odot FW\ar[r]^{\phi\odot\phi} &
F(X\otimes Z)\odot F(Y\otimes W)\ar[d]^{\phi} \\
F(X\otimes Y)\odot F(Z\otimes W)\ar[r]^{\phi} &
F(X\otimes Y\otimes Z\otimes W)\ar[r]^{F(\sigma^{(23)})} &
F(X\otimes Z\otimes Y\otimes W)
}
\end{align*}
\end{Lem}

\noindent{\bf Proof.} Due to the identity $\phi\cdot(\phi\odot\id)=\phi\cdot(\id\odot\phi)$ we obtain $\phi\cdot(\id\odot\phi)\cdot(\id\odot\phi\odot\id)= \phi\cdot(\id\odot\phi)\cdot(\id\odot\id\odot\phi)= \phi\cdot(\phi\odot\id)\cdot(\id\odot\id\odot\phi)=\phi\cdot(\phi\odot\phi)$. Then, by using~\eqref{phiphiDiag} one yields $\phi\cdot(\phi\odot\phi)\cdot\sigma^{(23)}= \phi\cdot(\phi\odot\id)\cdot\big(\id\odot(\phi\cdot\sigma)\odot\id\big)= \phi\cdot(\phi\odot\id)\cdot\big(\id\odot(F\sigma\cdot\phi)\odot\id\big)$. Due to the naturality of $\phi$ this is equal to $\phi\cdot\big(F(\id\otimes\sigma)\odot\id\big)\cdot(\phi\odot\id)\cdot(\id\odot\phi\odot\id)= F(\sigma^{(23)})\cdot\phi\cdot(\phi\odot\id)\cdot(\id\odot\phi\odot\id)= F(\sigma^{(23)})\cdot\phi\cdot(\phi\odot\phi)$. \qed

Note that $\Mon(\bfC)$ and $\Mon(\bfD)$ are both symmetric monoidal categories (see p.~\ref{bbMonC}). Let us describe a monoidal structure of $\Mon(F)\colon\Mon(\bfC)\to\Mon(\bfD)$ induced by the monoidal structure of $F$.

\begin{Prop} \label{PropMonFmon}
 If $F$ is a symmetric lax/strong monoidal functor, then $\Mon(F)$ is also lax/strong monoidal and also symmetric. The monoidal structure is given by the same morphisms $\varphi\colon I_\bfD\to F I_\bfC$ and $\phi_{X,Y}\colon FX\odot FY\to F(X\otimes Y)$ regarded as morphisms $\varphi\colon I_{\Mon(\bfD)}\to\Mon(F)(I_{\Mon(\bfC)})$ and $\phi_{\MM,\MM'}=\phi_{X,Y}\colon\wt\MM\odot\wt\MM'\to\Mon(F)(\MM\otimes\MM')$, where $\MM=(X,\mu_X,\eta_X)\in\Mon(\bfC)$, $\MM'=(Y,\mu_Y,\eta_Y)\in\Mon(\bfC)$, $\wt\MM=\Mon(F)\MM$ and $\wt\MM'=\Mon(F)\MM'$.

 Dually, if $F$ is a symmetric colax/strong monoidal functor, then $\Comon(F)$ is also colax/strong monoidal and also symmetric. Its monoidal structure is given by the morphisms $\varphi\colon I_\bfD\to F I_\bfC$ and $\phi_{X,Y}\colon F(X\otimes Y)\to FX\odot FY$ regarded as morphisms between comonoids.
\end{Prop}

\noindent{\bf Proof.} We need to check that $\phi_{X,Y}\colon FX\odot FY\to F(X\otimes Y)$ is a morphism from the monoid $\wt\MM\odot\wt\MM'=\big(FX\odot FY,(\mu_{FX}\odot\mu_{FY})\cdot\sigma^{(23)},\eta_{FX}\odot\eta_{FY}\big)$ to the monoid $\Mon(F)(\MM\otimes\MM')=\big(F(X\otimes Y),F(\mu_X\otimes\mu_Y)\cdot F(\sigma)\cdot\phi,F(\eta_X\otimes\eta_Y)\cdot\varphi\big)$, where $\mu_{FX}=F\mu_X\cdot\phi$, $\eta_{FX}=F\eta_X\cdot\varphi$. This amounts to the commutative diagrams
\begin{align} \label{phiXYXY}
\xymatrix{
FX\odot FY\odot FX\odot FY \ar[rr]^{\quad(\phi\odot\phi)\cdot\sigma^{(23)}} \ar[d]^{\phi\odot\phi} &&
 F(X\otimes X)\odot(Y\otimes Y)\ar[d]^{\phi} \ar[rr]^{\qquad F\mu_X\odot F\mu_Y} &&
FX\odot FY\ar[d]^{\phi} \\
F(X\otimes Y)\odot F(X\otimes Y) \ar[rr]^{\quad F(\sigma^{(23)})\cdot\phi} &&
F(X\otimes X\otimes Y\otimes Y) \ar[rr]^{\qquad F(\mu_X\otimes\mu_Y)} &&
F(X\otimes Y)
}
\end{align}
\begin{align} \label{phiIXY}
\xymatrix{
 I_\bfD\odot I_\bfD\ar@{=}[d]\ar[rr]^{\varphi\odot\varphi} &&
 FI_\bfC\odot FI_\bfC\ar[rr]^{F\eta_X\odot F\eta_Y}\ar[d]^{\phi_{I_\bfC,I_\bfC}} &&
 FX\odot FY\ar[d]^{\phi_{X,Y}} \\
 I_\bfD\ar[r]^{\varphi} &
 FI_\bfC\ar@{=}[r]&F(I_\bfC\otimes I_\bfC)\ar[rr]^{F(\eta_X\otimes\eta_Y)} &&
 F(X\otimes Y)
}
\end{align}
The commutativity of~\eqref{phiXYXY} follows from Lemma~\ref{Lemsigmaphi} (for $Z=X$, $W=Y$) and the naturality of $\phi$. The diagram~\eqref{phiIXY} is obtained by using the formula $\phi_{I_\bfC,I_\bfC}\cdot(\varphi\odot\id)=\id$. The same formula implies that $\varphi$ is a morphism from the monoid $I_{\Mon(\bfD)}=(I_{\bfD},\id,\id)$ to $\Mon(F)I_{\Mon(\bfC)}=(FI_\bfC,\phi_{I_\bfC,I_\bfC},\varphi)$. The conditions on the monoidal structure of $\Mon(F)$ defined in this way as well as the agreement with the symmetric structure are valid, since they are valid for the structure of $F$. If $F$ is strong monoidal, then $\varphi$ and $\phi_{X,Y}$ are isomorphisms in $\bfD$. They are isomorphisms in $\Mon(\bfD)$, since the forgetful functor $\Mon(\bfD)\to\bfD$ reflects isomorphisms. \qed

\bb{Translation of bimonoids.}
Let $F\colon(\bfC,\otimes)\to(\bfD,\odot)$ be a (covariant) symmetric strong monoidal functor. Let $\BB=(X,\mu_X,\eta_X,\Delta_X,\varepsilon_X)$ be a bimonoid in $(\bfC,\otimes)$. By definition it means that $\MM=(X,\mu_X,\eta_X)$ and $\OO=(X,\Delta_X,\varepsilon_X)$ are a monoid and a comonoid in $\bfC$ such that the comonoid structure gives the morphisms $\Delta_X\colon\MM\to\MM\otimes\MM$ and $\varepsilon_X\colon\MM\to I_{\Mon(\bfC)}$ in $\Mon(\bfC)$. The formulae~\eqref{muFX} give the monoid $\wt\MM=(FX,\mu_{FX},\eta_{FX})=\Mon(F)\MM$ in $\bfD$.
 By considering $F=(F,\varphi,\phi)$ as a colax monoidal functor with the structure morphisms $\varphi^{-1}\colon FI_\bfC\to I_\bfD$ and $\phi^{-1}_{X,Y}\colon F(X\otimes Y)\to FX\odot FY$ we obtain the comonoid $\wt\OO=(FX,\Delta_{FX},\varepsilon_{FX})=\Comon(F)\OO$ in $\bfD$, where $\Delta_{FX}$ and $\varepsilon_{FX}$ are defined by the formulae~\eqref{DeltaFX} but with $\varphi^{-1}$ and $\phi^{-1}$ instead of $\varphi$ and $\phi$, i.e.
\begin{align} \label{DeltaFXB}
 &\Delta_{FX}\colon FX\xrightarrow{F\Delta_X}F(X\otimes X)\xrightarrow{\phi^{-1}_{X,X}} FX\odot FX, &&\varepsilon_{FX}\colon FX\xrightarrow{F\varepsilon_X}FI_\bfC\xrightarrow{\varphi^{-1}} I_\bfD.
\end{align}
By applying Prop.~\ref{PropMonFmon} to the strong monoidal functor $(F,\varphi,\phi)$ we see that $\phi_{X,X}$ and $\varphi$ are monoid isomorphisms with the inverses $\phi^{-1}_{X,X}\colon\Mon(F)(\MM\otimes\MM)\isoright\wt\MM\odot\wt\MM$ and $\varphi^{-1}\colon\Mon(F)I_{\Mon(\bfC)}\isoright I_{\Mon(\bfD)}$, where $\wt\MM=\Mon(F)\MM$. Hence the compositions~\eqref{DeltaFXB} give the monoid morphisms $\Delta_{FX}\colon\wt\MM\to\wt\MM\odot\wt\MM$ and $\varepsilon_{FX}\colon\wt\MM\to I_{\Mon(\bfD)}$, so we get a bimonoid $\wt\BB=(FX,\mu_{FX},\eta_{FX},\Delta_{FX},\varepsilon_{FX})\in\Bimon(\bfD)$. We have proved the following fact.

\begin{Prop}
A (covariant) symmetric strong monoidal functor $F\colon(\bfC,\otimes)\to(\bfD,\odot)$ translates bimonoids to bimonoids. Namely, it induces the symmetric strong monoidal functor
\begin{align}
 \Bimon(F)\colon\Bimon(\bfC,\otimes)\to\Bimon(\bfD,\odot)
\end{align}
agreed with the functors~\eqref{MonF} and~\eqref{ComonF}. It maps $\BB$ to the bimonoid $\wt\BB=\Bimon(F)\BB$ defined by the formulae~\eqref{muFX}, \eqref{DeltaFXB}.

A contravariant symmetric strong monoidal functor $F\colon(\bfC,\otimes)\to(\bfD,\odot)$ gives the contravariant functor
 $\Bimon(F)\colon\Bimon(\bfC,\otimes)\to\Bimon(\bfD,\odot)$
corresponding to covariant functor $\Bimon(\bar F)\colon\Bimon(\bfC^\op,\otimes) \to\Bimon(\bfD,\odot)$. It translates $\BB\in\Bimon(\bfC,\otimes)$ to $\wh\BB=\Bimon(F)\BB=(FX,\mu_{FX},\eta_{FX},\Delta_{FX},\varepsilon_{FX})\in\Bimon(\bfD,\odot)$, where
\begin{align}
 &\mu_{FX}\colon FX\odot FX\xrightarrow{\phi_{X,X}} F(X\otimes X)\xrightarrow{F\Delta_X}FX,
 &&\eta_{FX}\colon I_\bfD\xrightarrow{\varphi}FI_\bfC\xrightarrow{F\varepsilon_X}FX, \label{whBB} \\
 &\Delta_{FX}\colon FX\xrightarrow{F\mu_X}F(X\otimes X)\xrightarrow{\phi_{X,X}^{-1}}FX\odot FX,
 &&\varepsilon_{FX}\colon FX\xrightarrow{F\eta_X}FI_\bfC\xrightarrow{\varphi^{-1}}I_\bfD. \label{whBBD}
\end{align}
\end{Prop}

\bbo{Translation of tensor product of representations.}
Suppose that the parametrising subcategories $\bfP\subset\bfC$ and $\bfQ\subset\bfD$ are monoidal subcategories of the symmetric monoidal categories $(\bfC,\otimes)$ and $(\bfD,\odot)$ respectively and that $F(\bfP)\subset\bfQ$ for a covariant strong monoidal functor $F$. Let $\BB=(X,\mu_X,\eta_X,\Delta_X,\varepsilon_X)\in\Bimon(\bfC,\otimes)$ and $\wt\BB=\Bimon(F)\BB$.

If $\bfC=(\bfC,\otimes)$ and $\bfD=(\bfD,\odot)$ are both relatively closed, then we have the functors~\eqref{FRep}. Denote such functor for a monoid $\MM=(X,\mu_X,\eta_X)$ by
\begin{align} \label{FRepBB}
 F_{\bfP,\BB}\colon\Rep_\bfP(\BB)\to\Rep_\bfQ(\wt\BB).
\end{align}

In the case, when $\bfC$ and $\bfD$ are both relatively coclosed, the functor~\eqref{FCorep} for the comonoid $\OO=(X,\Delta_X,\varepsilon_X)$ and the colax monoidal functor $F=(F,\varphi^{-1},\phi^{-1})$ is defined and has the form
\begin{align} \label{FCorepBB}
 F^\BB_{\bfP}\colon\Corep_\bfP(\BB)\to\Corep_\bfQ(\wt\BB).
\end{align}

If $\bfC=(\bfC,\otimes)$ and $\bfD=(\bfD,\odot)$ are relatively closed and coclosed respectively and $F$ is a contravariant symmetric strong monoidal functor, then we get the contravariant functor
\begin{align} \label{FRepConBB}
 F^\BB_{\bfP}\colon\Rep_\bfP(\BB)\to\Corep_\bfQ(\wh\BB),
\end{align}
where $\wh\BB$ is defined by~\eqref{whBB}, $\eqref{whBBD}$.

Recall that all the categories appeared in~\eqref{FRepBB}, \eqref{FCorepBB} and \eqref{FRepConBB} are monoidal categories (Theorems~\ref{ThrhoTP}, \ref{ThomegaTP}).

\begin{Th} \label{ThFRepBB}
 The functors~\eqref{FRepBB}, \eqref{FCorepBB} and \eqref{FRepConBB} are strong monoidal functors (the corresponding conditions are supposed: the symmetric monoidal categories $\bfC$ and $\bfD$ should be relatively closed/coclosed, the symmetric strong monoidal $F$ should be covariant/contravariant). The monoidal structure of the functors~\eqref{FRepBB}, \eqref{FCorepBB} and \eqref{FRepConBB} is inherited from $F$.
\end{Th}

\noindent{\bf Proof.} Let $\rho$ and $\rho'$ be representations of $\BB$ on the objects $V$ and $V'$. Their tensor product $\rho''=\pi\cdot(\rho\otimes\rho')\cdot\Delta_X$ is a representation of $\BB$ on $V''=V\otimes V'$. The functor $F$ translates the objects $(V,\rho)$, $(V',\rho')$ and $(V'',\rho'')$ of $\Rep_\bfP(\BB)$ to the objects $(FV,\tau)$, $(FV',\tau')$ and $(FV'',\tau'')$ of $\Rep_\bfQ(\wt\BB)$, where $\tau=\Phi_{V,V}\cdot F\rho$, $\tau'=\Phi_{V',V'}\cdot F\rho'$ and $\tau''=\Phi_{V'',V''}\cdot F\rho''$. The tensor product of $\tau$ and $\tau'$ is the representation $\lambda=\pi\cdot(\tau\odot\tau')\cdot\Delta_{FX}$ on $FV\odot FV'$. We need to show that $\phi_{V,V'}\colon FV\odot FV'\to F(V\otimes V')$ is a morphism between representations $\lambda$ and $\tau''$. According to Definition~\ref{DefRep} this amounts to the commutative diagram
\begin{align} \label{lambdatauphi}
\xymatrix{
X\ar[rr]^{\lambda\qquad}\ar[d]_{\tau''} && \bfhom(FV\odot FV',FV\odot FV')\ar[d]^{\bfhom(\id,\phi)} \\
\bfhom(FV'',FV'')\ar[rr]^{\bfhom(\phi,\id)\qquad} && \bfhom(FV\odot FV',FV'')
}
\end{align}
The representations in this diagram can be written in the form
\begin{align*}
 &\lambda=\pi\cdot(\tau\odot\tau')\cdot\Delta_{FX}=
 \pi\cdot(\Phi\odot\Phi)\cdot(F\rho\odot F\rho')\cdot\phi^{-1}\cdot F\Delta_{X}= \\
 &\qquad\qquad \pi\cdot(\Phi\odot\Phi)\cdot\phi^{-1}\cdot F(\rho\otimes\rho')\cdot F\Delta_{X}, \\[7pt]
 &\tau''=\Phi\cdot F\rho''=\Phi\cdot F\pi\cdot F(\rho\otimes\rho')\cdot F\Delta_X,
\end{align*}
so we need to check the equality
$\bfhom(\id,\phi)\cdot\pi\cdot(\Phi\odot\Phi)=\bfhom(\phi,\id)\cdot\Phi\cdot F\pi\cdot\phi$ of two morphisms $F\big(\bfhom(V,V)\big)\odot F\big(\bfhom(V',V')\big)\to\bfhom(FV\odot FV',FV'')$. By using the formulae~\eqref{thetafg},~\eqref{thetahom} and the commutative diagrams~\eqref{piDiag}, \eqref{PhiDiag} we obtain their images under $\theta$:
\begin{multline} \label{lambdaM}
\theta\big(\bfhom(\id,\phi)\big)\cdot(\pi\odot\id)\cdot(\Phi\odot\Phi\odot\id)= \phi\cdot\eval\cdot(\pi\odot\id)\cdot(\Phi\odot\Phi\odot\id)= \\
\phi\cdot(\eval\odot\eval)\cdot\sigma^{(23)}\cdot(\Phi\odot\Phi\odot\id)=
\phi\cdot(\eval\odot\eval)\cdot(\Phi\odot\id\odot\Phi\odot\id)\cdot\sigma^{(23)}= \\
\phi\cdot\big(F(\eval)\odot F(\eval)\big)\cdot(\phi\odot\phi)\cdot\sigma^{(23)}=
F(\eval\otimes\eval)\cdot\phi\cdot(\phi\odot\phi)\cdot\sigma^{(23)},
\end{multline}
\begin{multline} \label{tauM}
\theta\big(\bfhom(\phi,\id)\big) \cdot(\Phi\odot\id)\cdot(F\pi\odot\id)\cdot(\phi\odot\id)= \\
\eval\cdot(\id\odot\phi)\cdot(\Phi\odot\id)\cdot(F\pi\odot\id)\cdot(\phi\odot\id)=
F(\eval)\cdot\phi\cdot(F\pi\odot\id)\cdot(\phi\odot\phi)= \\
F(\eval)\cdot F(\pi\otimes\id)\cdot\phi\cdot(\phi\odot\phi)=
F(\eval\otimes\eval)\cdot F(\sigma^{(23)})\cdot\phi\cdot(\phi\odot\phi).
\end{multline}
The right hand sides of~\eqref{lambdaM} and \eqref{tauM} coincide by virtue of Lemma~\ref{Lemsigmaphi}, so the diagram~\eqref{lambdatauphi} is commutative.

The unit object $I_{\Rep_\bfP(\BB)}=(I_\bfC,\varepsilon_X)$ is translated to the object $(FI_\bfC,\wt\varepsilon_X)$, where
\begin{align}
 \wt\varepsilon_X\colon FX\xrightarrow{ F\varepsilon_X}FI_\bfC=F\big(\bfend(I_\bfC)\big)\xrightarrow{\Phi_{I_\bfC,I_\bfC}}\bfend(FI_\bfC)
\end{align}
is the corresponding representation of $\wt\BB$ on $FI_\bfC$. Let us prove that $\varphi\colon I_\bfD\to FI_\bfC$ gives the morphisms between the representations $\varepsilon_{FX}=\varphi^{-1}\cdot F\varepsilon_X\colon FX\to I_D=\bfend(I_\bfD)$ and $\wt\varepsilon_X$. Due to Prop.~\ref{ProphomIW} the morphism $I_\bfD=\bfhom(I_\bfD,I_\bfD) \xrightarrow{\bfhom(\id,\varphi)}\bfhom(I_\bfD,FI_\bfC)=FI_\bfC$ coincides with $\varphi$, so we obtain
\begin{align*}
 &\bfhom(\id,\varphi)\cdot\varepsilon_{FX}=\varphi\cdot\varphi^{-1}\cdot F\varepsilon_X=F\varepsilon_X, &
 &\bfhom(\varphi,\id)\cdot\wt\varepsilon_X=\bfhom(\varphi,\id)\cdot\Phi\cdot F\varepsilon_X.
\end{align*}
Hence we need only to check that the morphism
\begin{align} \label{Phihomvarphi}
 F\big(\bfhom(I_\bfC,I_\bfC)\big)\xrightarrow{\Phi}\bfhom(FI_\bfC,FI_\bfC)\xrightarrow{\bfhom(\varphi,\id)}\bfhom(I_\bfD,FI_\bfC)
\end{align}
coincides with the identification $F\big(\bfhom(I_\bfC,I_\bfC)\big)=FI_\bfC=\bfhom(I_\bfD,FI_\bfC)$. Since~\eqref{evalIW} is the identification, the morphism~\eqref{Phihomvarphi} coincides with the image of~\eqref{Phihomvarphi} itself under $\theta$. Due to~\eqref{thetafg} and~\eqref{thetahom} this image equals
\begin{multline}
 \theta(\bfhom(\varphi,\id))\cdot(\Phi\odot\id)= \eval\cdot(\id\odot\varphi)\cdot(\Phi\odot\id)=\eval\cdot(\Phi\odot\id)\cdot(\id\odot\varphi)=\\
 F(\eval)\cdot\phi\cdot(\id\odot\varphi)\colon F\big(\bfhom(I_\bfC,I_\bfC)\big)\odot I_\bfD\to FI_\bfC. \label{idFIComp}
\end{multline}
Since $F(\eval)\colon F\big(\bfhom(I_\bfC,I_\bfC)\otimes I_\bfC\big)\to FI_\bfC$ and
 $$F\big(\bfhom(I_\bfC,I_\bfC)\big)\odot I_\bfD\xrightarrow{\id\odot\varphi}F\big(\bfhom(I_\bfC,I_\bfC)\big)\odot FI_\bfC\xrightarrow{\phi} F\big(\bfhom(I_\bfC,I_\bfC)\otimes I_\bfC\big)$$
 are both the identification morphisms, their composition~\eqref{idFIComp} is the identification as well.

By virtue of Prop.~\ref{PropRepIso} the morphisms of representations $\phi_{V,V'}\colon(FV\odot FV',\lambda)\to(FV'',\tau'')$ and $\varphi\colon(I_\bfD,\varepsilon_{FX})\to(FX,\wt\varepsilon_X)$ are isomorphisms in $\Rep_\bfP(\BB)$. \qed

\begin{Rem} \normalfont
Recall that any symmetric strong monoidal functor $F=(F,\varphi,\phi)$ induces strong monoidal functors
\begin{align} \label{FLactB}
 &F_\BB\colon\Lact(\BB)\to\Lact(\wt\BB), &
 &F^\BB\colon\Lcoact(\BB)\to\Lcoact(\wt\BB)
\end{align}
with the (lax) monoidal structure $(\varphi,\phi)$. The functors~\eqref{FRepBB} and \eqref{FCorepBB} are restrictions of~\eqref{FLactB} as monoidal functors (in the sense of Prop.~\ref{PropRepLact} and Remark~\ref{RemRepLact}). These facts can be used for an alternative proof of Theorem~\ref{ThFRepBB}.
\end{Rem}

\section{Examples}
\label{secEx}

The general theory of representations described above was motivated by two types of representations. First, this theory allows to unify the classical representations on vector spaces such as representations of groups and algebras. Second example is quantum representations on quadratic algebras defined in~\cite{Sqrt}.

Let us fix a basic field $\KK$ of characteristic $\chara\KK\ne2$. We suppose that any vector space or algebra is over $\KK$ unless we explicitly specify another field.

\subsection{Representations on vector spaces}
\label{secCR}

The model example of the general representation theory is the representations of algebras. This is the case of the closed monoidal category of vector spaces. To describe the representations of groups we introduce a bigger relatively closed monoidal category which contains sets as well as vector spaces as monoidal subcategories. This gives the unification of representations of groups and algebras. This case is generalised to groups with a geometric structure such as Lie groups, algebraic groups etc (the algebraic case was described in~\cite[\S~4.3]{Sqrt}).

\bb{Representations of algebras.}
The category of vector spaces $(\Vect,\otimes)$ with the standard tensor product $\otimes$ is a closed symmetric monoidal category with the unit object $\KK$. The monoids in $(\Vect,\otimes)$ are algebras: $\Mon(\Vect,\otimes)=(\Alg,\otimes)$. For $V,W\in\Vect$ the object $\bfhom(V,W)$ is a vector space of all the linear operators $V\to W$, while $\bfend(V)$ is the algebra of operators on $V$. A morphism $\rho\colon\gA\to\bfend(V)$ in $\Mon(\Vect,\otimes)=\Alg$ is a representation of an algebra $\gA\in\Alg$ on $V\in\Vect$ in the usual sense.

\bb{Semi-linear sets.}
The category of sets $\Set$ is the category with finite products. The usual monoids and groups are monoids and groups in $\Set$, i.e. monoids and Hopf monoids in $(\Set,\times)$. Notice that the monoidal category $(\Set,\times)$ is closed with the hom-objects $\bfhom(X,Y)=\Hom(X,Y)\in\Set$ and $\bfend(X)=\End(X)\in\Mon(\Set,\times)$, where $X,Y\in\Set$.

To consider representations of monoids and groups on vector spaces in frame of the general representation theory constructed in Section~\ref{secGRT} we define the category $\SLSet$. Its objects are pairs $(X,V)$, where $X\in\Set$ and $V\in\Vect$; denote such object by $X\times V$. Define morphisms by the formula
\begin{align} \label{HomSLSet}
 \Hom_\SLSet(X\times V,Y\times W)=\Hom_\Set(X,Y)\times\Hom_\Set\big(X,\Hom_\Vect(V,W)\big).
\end{align}
A morphism $X\times V\to Y\times W$ is presented by a collection $\big(\varphi,(f_x)_{x\in X}\big)$, where $\varphi\colon X\to Y$ is an arbitrary map and $f_x\colon V\to W$ are linear maps indexed by $x\in X$.

We can consider the object $X\times V\in\SLSet$ as the Cartesian product of the sets $X$ and $V$. A morphism $F=(\varphi,f_x)$ in $\SLSet$ can be considered as the map $F\colon X\times V\to Y\times W$ given by the formula
\begin{align} \label{Fxv}
 F(x,v)=\big(\varphi(x),f_x(v)\big).
\end{align}
Such maps are linear only in the second variable, so we call them {\it semi-linear maps}.
If $(\psi,g_y)\colon Y\times W\to Z\times U$ is another morphism in $\SLSet$ for some $Y\in\Set$ and $U\in\Vect$, then the composition $X\times V\xrightarrow{(\varphi,f_x)}Y\times W\xrightarrow{(\psi,g_y)} Z\times U$ has the form
\begin{align}
 (\psi,g_y)\cdot(\varphi,f_x)=(\psi\cdot\varphi,g_{\varphi(x)}\cdot f_x).
\end{align}
This formula defines the composition in $\SLSet$. Due to the associativity of the composition of maps we indeed defined a category $\SLSet$. It is a concrete category with the forgetful functor $\SLSet\to\Set$ defined by~\eqref{Fxv}. We call its objects {\it semi-linear spaces} or {\it semi-linear sets}.

Let $X,X'\in\Set$ and $V,V'\in\Vect$. We introduce the tensor product of the semi-linear sets $X\times V$ and $X'\times V'$ by the formula
\begin{align} \label{TPObSL}
 (X\times V)\otimes(X'\times V')=(X\times X')\times(V\otimes V').
\end{align}
The tensor product on morphisms is defined as
\begin{align} \label{TPMorSL}
 \big(\varphi,(f_x)_{x\in X}\big)\otimes\big(\varphi',(f'_{x'})_{x'\in X'}\big)=\big(\varphi\times\varphi',(f_x\otimes f'_{x'})_{(x,x')\in X\times X'}\big), 
\end{align}
where $\varphi\colon X\to Y$, $f_x\colon V\to W$, $\varphi'\colon X'\to Y'$, $f'_{x'}\colon V'\to W'$.
This gives a symmetric monoidal category $(\SLSet,\otimes)$ with the unit object $\{0\}\times\KK$, where $\{0\}\in\Set$ and $\KK\in\Vect$ are unit objects of $(\Set,\times)$ and $(\Vect,\otimes)$ respectively. There is one more useful `forgetful' functor $\SLSet\to\Set$, $X\times V\mapsto X$. It is a strict monoidal functor $(\SLSet,\otimes)\to(\Set,\times)$ and hence it induces $\Mon(\SLSet,\otimes)\to\Mon(\Set,\times)$.

An object $\SSS\in\Mon(\SLSet,\otimes)$ is a semi-linear set $X\times V$ with the semi-linear maps
\begin{align}
 &\mu_{X\times V}=\big(\mu_X,(f_{x,y})_{(x,y)\in X\times X}\big)\colon (X\times X)\times(V\otimes V)\to X\times V, \\
 &\eta_{X\times V}=(\eta_X,\eta_V)\colon\{0\}\times\KK\to X\times V
\end{align}
satisfying the following conditions. The triple $(X,\mu_X,\eta_X)$ should be a monoid in $(\Set,\times)$, denote $e=\eta_X(0)\in X$ and $xy=\mu_X(x,y)$. The linear maps $f_{x,y}\colon V\otimes V\to V$ and $\eta_V\colon\KK\to V$ are subjected to $f_{x,yz}\cdot(\id\otimes f_{y,z})=f_{xy,z}\cdot(f_{x,y}\otimes\id)$ and $f_{e,x}\cdot(\eta_V\otimes\id)=\id=f_{x,e}\cdot(\id\otimes\eta_V)$.
We call $\SSS=(X\times V,\mu_{X\times V},\eta_{X\times V})\in\Mon(\SLSet,\otimes)$ a {\it semi-linear} monoid.

The fully faithful strict monoidal functor
\begin{align} \label{VectEmb}
 &(\Vect,\otimes)\hookrightarrow(\SLSet,\otimes), &&V\mapsto \{0\}\times V, 
\end{align}
embeds the category $\Vect$ into $\SLSet$ as a full monoidal subcategory of the monoidal category $(\SLSet,\otimes)$.

We embed $\Set$ into $\SLSet$ by the faithful strict monoidal functor
\begin{align} \label{SetEmb}
 &(\Set,\times)\hookrightarrow(\SLSet,\otimes), &&X\mapsto X\times\KK. 
\end{align}
This embedding functor is not full. It allows us to regard the category $\Set$ as a (not full) monoidal subcategory of $(\SLSet,\otimes)$.

\begin{Rem} \normalfont
 The category $\SLSet$ is also a category with finite products and $(\Set,\times)$ can be considered as a monoidal subcategory of $(\SLSet,\times)$ via another embedding $\Set\to\SLSet$, $X\mapsto X\times0$. The embedding $V\mapsto\{0\}\times V$ does not respect the tensor product $\otimes$ of vector spaces, it is the strict monoidal functor $(\Vect,\oplus)\to(\SLSet,\times)$. However, these embeddings will not give the classical representations on vector spaces.
\end{Rem}

\bbo{Representations of usual monoids and groups.} \label{bbRepUMG}
A classical notion of the representation of a monoid $\MM\in\Mon(\Set,\times)$ is a monoid homomorphism $\rho\colon\MM\to\End(V)$, where $\End(V)$ is the monoid of the linear operators on $V$. When $\MM$ is a group, the image of $\rho$ lies in the group $GL(V)\subset\End(V)$ consisting of the invertible operators on $V$, so the notion of the representation of a monoid on a linear space includes the linear representations of groups.

The $\Mon$-functor~\eqref{MonF} induced by the fully faithful strict monoidal functor~\eqref{VectEmb} is the category embedding $\Alg\hookrightarrow\Mon(\SLSet,\times)$. It identifies an algebra $\gA=(\gA,\mu_V,\eta_V)\in\Alg$ with the semi-linear monoid $\{0\}\times\gA$; in this case $X=\{0\}$, the maps $\mu_X$, $\eta_X$ are trivial and $f_{0,0,0}=\mu_V\colon V\otimes V\to V$, $\eta_V\colon V\to\KK$ are the structure morphisms of the algebra $\gA$.

The functor~\eqref{SetEmb} induces the embedding $\Mon(\Set,\times)\hookrightarrow\Mon(\SLSet,\otimes)$. The monoid $\MM=(X,\mu_X,\eta_X)\in\Mon(\Set,\times)$ is identified with the semi-linear monoid
\begin{align} \label{SSSM}
 &\SSS_\MM=\Big(X\times\KK,\big(\mu_X,(\id_\KK)\big),(\eta_V,\id_\KK)\Big)\in\Mon(\SLSet,\otimes). 
\end{align}
We will see that the representations of a monoid $\MM$ on $V$ are exactly the representations of $\SSS_\MM$ on the object $V$ identified with $\{0\}\times V\in\SLSet$ via the embedding~\eqref{VectEmb}.

First, we need to check that $\bfC=(\SLSet,\otimes)$ is closed relative to $\bfP=\Vect$. We, however, can prove that it is closed (relative to the whole $\SLSet$). For two semi-linear sets $X\times V,Y\times W\in\SLSet$ define an object $\bfhom(X\times V,Y\times W)\in\SLSet$ as the Cartesian product~\eqref{HomSLSet}, where $\Hom\big(X,\Hom(V,W)\big)=\Hom\big(X,\bfhom(V,W)\big)$ is equipped by the structure of a linear space as follows: the linear combination of $f,g\colon X\to\bfhom(V,W)$ with the coefficients $\alpha,\beta\in\KK$ is the map $X\to\bfhom(V,W)$ with the values
\begin{align}
 &(\alpha f+\beta g)(x)=\alpha f(x)+\beta g(x)\in\bfhom(V,W), &&x\in X.
\end{align}
Introduce the map $\eval_{X\times V,Y\times W}\colon\bfhom(X\times V,Y\times W)\otimes(X\times V)\to Y\times W$, which evaluates the map $(\varphi,f_x)\in\bfhom(X\times V,Y\times W)$ on $(x,v)\in X\times V$. It maps an element
\begin{align}
 (\varphi,x,f\otimes v)\in\Hom(X,Y)\times X\times\Big(\Hom\big(X,\bfhom(V,W)\big)\otimes V\Big)
\end{align}
to $\big(\varphi(x),f(x)(v)\big)\in Y\times W$. One can show that $\big(\bfhom(X\times V,Y\times W),\eval_{X\times V,Y\times W}\big)$ is a universal arrow from the functor $-\otimes(X\times V)\colon\SLSet\to\SLSet$ to the object $Y\times W$. Then Theorem~\ref{ThUnit} implies that $(\SLSet,\otimes)$ is closed with the hom-objects defined in this way.

For the case $X=Y=\{0\}$ we have $\{0\}\times\bfhom(V,W)=\bfhom(\{0\}\times V,\{0\}\times W)$ (this identification is the isomorphism $\Phi_{V,W}$ induced by the monoidal functor~\eqref{VectEmb}), so the internal hom-object in $(\Vect,\otimes)$ is the same as the internal hom-object in $(\SLSet,\otimes)$ for the vector spaces embedded by~\eqref{VectEmb}. The internal end is the algebra $\bfend(\{0\}\times V)=\{0\}\times\bfend(V)$ identified with the algebra of operators $\bfend(V)$.

A representation of an algebra $\gA=\{0\}\times\gA$ on $V=\{0\}\times V$ in the sense of Definition~\ref{DefRep} is exactly the algebra morphism $\rho\colon\gA\to\bfend(V)$, so it is a representation of the algebra $\gA$ in the usual sense.

Consider a monoid $\MM\in\Mon(\Set,\times)$ as the semi-linear monoid $\SSS_\MM$ defined by~\eqref{SSSM}. Then a representation of $\MM$ on $V=\{0\}\times V$ is a representation $\wh\rho\colon\SSS_\MM\to\bfend(V)$. This is a morphism $X\times\KK\to\{0\}\times\bfend(V)$ in $\SLSet$ given by the trivial map $0\colon X\to\{0\}$ and a collection of linear maps $f_x=\rho_x\colon\KK\to\bfend(V)$. Define the map $\rho\colon X\to\bfend(V)$ as $\rho(x)=\rho_x(1)$. Any collection $\rho_x$ is given by such map $\rho$. The morphism $\wh\rho=(0,\rho_x)$ is a representation (a morphism in $\Mon(\SLSet,\otimes)$) iff $\rho(x)\rho(y)=\rho(xy)$ for any $x,y\in X$. The latter condition means exactly that $\rho$ is a representation of the monoid $\MM$ on $V$.

Thus a linear representation of a usual monoid or group $\MM$ is its representation in the monoidal category $(\SLSet,\otimes)$ on an object $V$ of the subcategory $\Vect\subset\SLSet$.

\bb{Representations of Lie groups.}
Consider the category of smooth manifolds $\ManInf$. This is a category with finite products. A monoid $\MM\in\Mon(\ManInf,\times)$ is a manifold $X$ with a smooth multiplication $X\times X\to X$ and a unity $e\in X$. If any element $x\in X$ is invertible with respect to this multiplication and the map $X\to X$, $x\mapsto x^{-1}$, is smooth, then $\MM$ is a Lie group, so the category of Lie groups is a full subcategory of $\Mon(\ManInf,\times)$.

Suppose that $\KK$ is an extension of the field of real numbers: $\RR\subset\KK$.
A representation of a Lie group or monoid $\MM\in\Mon(\ManInf,\times)$ on a vector space $V\in\FVect$ is a smooth monoid homomorphism $\rho\colon\MM\to\bfend(V)$, where $\FVect$ is the category of finite-dimensional vector spaces and the algebra $\bfend(V)$ is considered as a smooth manifold $V^*\otimes V$ (this is possible for finite dimensional $V$ only). Such representations can be interpreted in terms of Section~\ref{secGRT} by introducing a smooth version of semi-linear spaces.

Define the category $\SLManInf$. An object of this category is a pair $(X,V)$ of $X\in\ManInf$, $V\in\FVect$, denoted by $X\times V$ and regarded as a Cartesian product of smooth manifolds. The morphisms $F\colon X\times V\to Y\times W$ are smooth maps of the form~\eqref{Fxv}. The objects $X\times V\in\SLManInf$ are called {\it smooth semi-linear spaces} or {\it semi-linear manifolds}.

 We obtain the following faithful (`forgetful') functors. By forgetting semi-linear structure of $X\times V$ we can consider it as an abstract manifold, this gives the forgetful functor $\SLManInf\to\ManInf$. By forgetting the smoothness we obtain $\SLManInf\to\SLSet$. Finally we have the functor $\SLManInf\to\ManInf$ which omits the linear part: $X\times V\mapsto X$.

The category $\SLManInf$ has a monoidal structure with the tensor product defined by~\eqref{TPObSL}, \eqref{TPMorSL}. Monoids in $(\SLManInf,\otimes)$ are called {\it smooth semi-linear monoids}.
The formulae $X\mapsto X\times\KK$ and $V\mapsto \{0\}\times V$ define the embeddings (also faithful and fully faithful respectively)
\begin{align} \label{SLManEmb}
 &(\ManInf,\times)\hookrightarrow(\SLManInf,\otimes), &
 &(\FVect,\otimes)\hookrightarrow(\SLManInf,\otimes).
\end{align}
 The monoidal category $(\SLManInf,\otimes)$ is symmetric, but it is not closed in the usual sense. It is closed relative to the subcategory $\bfP=\FVect$ defined by the right embedding~\eqref{SLManEmb}.
The functors~\eqref{SLManEmb} induce the embeddings
\begin{align} \label{SLMonManEmb}
 &\Mon(\ManInf,\times)\hookrightarrow\Mon(\SLManInf,\otimes), &
 &\FAlg\hookrightarrow\Mon(\SLManInf,\otimes),
\end{align}
which identifies a monoid $\MM\in\Mon(\ManInf,\times)$ (in particular, a Lie group $\MM$) and a finite-dimensional algebra $\gA\in\FAlg$ with the corresponding smooth semi-linear monoids $\SSS_\MM$ and $\{0\}\times\gA$ respectively. A representation of $\MM$ and of $\gA$ on a vector space $V\in\FVect$ can be interpreted as a representation of $\SSS_\MM$ and of $\{0\}\times\gA$ respectively in the same way as for the case of p.~\ref{bbRepUMG}. In particular, the linear representations of Lie groups are representations in the monoidal category $(\SLManInf,\otimes)$.

Similarly, one can interpret the representations of a complex Lie group by using the category of complex manifolds instead of $\ManInf$. In this case we should suppose $\CC\subset\KK$.

\begin{Rem} \normalfont
 Consider the category of vector bundles over non-fixed manifold. A morphisms of two such bundles $E$ and $E'$ with projections $\pi\colon E\to X$ and $\pi'\colon E'\to X'$ is a pair of smooth maps $F\colon E\to E'$ and $\varphi\colon X\to X'$ such that $\varphi\cdot\pi=\pi'\cdot F$ and the induced maps $F_x\colon\pi^{-1}(x)\to(\pi')^{-1}\big(\varphi(x)\big)$ are linear.
 The category $\SLManInf$ coincides with the subcategory of trivial vector bundles $E=X\times V$ (we use the term `semi-linear space' instead of `trivial vector bundle' by the reasons explained in~\cite[\S~4.3.3]{Sqrt}). The role of $(\bfC,\otimes)$ could be played by the category of all the locally trivial vector bundles with the tensor product defined via local triviality, but its subcategory $\SLManInf$ is enough to consider representations of Lie groups from the point of view of the general representation theory.
\end{Rem}

\bbo{Representations of affine algebraic groups.}
Let us fix an infinite subfield $\FF\subset\KK$ and consider the category $\AlgSet_\FF$ of the (affine) algebraic sets over $\FF$.
Since $\AlgSet_\FF$ is a category with finite products, we have the monoidal category $(\AlgSet_\FF,\times)$. Monoids/Hopf monoids in this monoidal category are called {\it (affine) algebraic monoids/groups} (over $\FF$). By using these category instead of $\Set$ or $\ManInf$ we can interpret representations of algebraic groups and monoids on vector spaces in terms of the general representation theory. The case $\FF=\KK$ is described in details in~\cite[\S~4.3]{Sqrt}.

Remind the main notions. By considering the Cartesian products $X\times V$ for $X\in\AlgSet_\FF$ and $V\in\FVect$ as algebraic sets we get the category $\SLAlgSet_\FF$ (the index $\FF$ means that $X$ is over $\FF$, the linear part $V$ is over $\KK$). Its objects are called {\it semi-linear (affine) algebraic sets} or {\it (affine algebraic) semi-linear spaces}. Their morphisms are regular maps of the form~\eqref{Fxv}; their tensor products are defined by the same formula~\eqref{TPObSL}. The objects of $\Mon(\SLAlgSet_\FF,\times)$ are called {\it semi-linear algebraic monoids} (over $\FF$).

The monoidal categories $(\AlgSet_\FF,\times)$ and $(\FVect,\times)$ are identified with the monoidal subcategories of $(\SLAlgSet_\FF,\otimes)$ by $X\mapsto X\times\KK$, $V\mapsto\{0\}\times V$. Then $\bfC=(\SLAlgSet_\FF,\otimes)$ is closed relative to $\bfP=\FVect$. A representation of an algebraic monoid/group or finite-dimensional algebra on $V\in\FVect$ is exactly a representation of the corresponding semi-linear algebraic monoid on $V$.

\subsection{Quantum representations on quadratic algebras}
\label{secQR}

Now let us remind the main points of Quantum Representation Theory developed in~\cite{Sqrt} and apply the general representation theory to this case. Note that the quantum representations can be considered as generalisations of the representations of algebras, algebraic groups and semi-linear algebraic monoids.

\bb{Quantum linear spaces.}
For two $\NN_0$-graded algebras $\A=\bigoplus\limits_{k\ge0}\A_k$ and $\B=\bigoplus\limits_{k\ge0}\B_k$ we denote $\A\circ\B=\bigoplus\limits_{k\ge0}(\A_k\otimes\B_k)$. This operation was introduced by Manin in~\cite{Manin87,Manin88} for quadratic algebras. We call it {\it (white) Manin product}, it defines a symmetric monoidal structure on the category of all the $\NN_0$-graded algebras $\GrAlg$. The unit object is the polynomial algebra $\KK[u]=\bigoplus\limits_{k\ge0}\KK u^k$.

The finitely generated quadratic algebras over $\KK$ (i.e. connected) form a monoidal subcategory $(\FQA,\circ)\subset(\GrAlg,\circ)$. This subcategory is a coclosed monoidal category with the internal cohom-objects $\bfcohom(\A,\B)=\A^!\bullet\B\in\FQA$, $\A,\B\in\FQA$ (see~\cite[\S~4]{Manin88}, \cite[\S~4.4]{Sqrt}). Manin interpreted the opposite category $\FQA^\op$ as a quantum (non-commutative) version of $\FVect$ and called its objects the {\it quantum linear spaces}. The classical finite-dimensional vector spaces are lifted to the quantum level by the contravariant colax monoidal embedding $S^*\colon(\FVect,\otimes)\hookrightarrow(\FQA,\circ)$, $V\mapsto SV^*$. Alternatively one can use the contravariant strong monoidal embedding $T^*\colon(\FVect,\otimes)\hookrightarrow(\FQA,\circ)$, $V\mapsto TV^*$ (see~\cite{Sqrt} for details).

The corepresentations of a comonoid $\OO=(\A,\Delta,\varepsilon)\in\Comon(\FQA,\circ)$ on a quadratic algebra $\B\in\FQA$ or, equivalently, the representations of a monoid in $(\FQA^\op,\circ)$ on a quantum linear space $\B\in\FQA^\op$ are quantum analogues for the finite-dimensional representations of finite-dimensional algebras.

The category $\QA$ consisting of all the quadratic algebras over $\KK$ (including infinitely generated) is a monoidal subcategory of $(\GrAlg,\circ)$. We proved in~\cite[\S~4.4.5]{Sqrt} that for any $\B\in\FQA$ the functor $-\circ\B\colon\QA\to\QA$ has a left adjoint (relative to the whole $\QA$), hence the symmetric monoidal category $(\QA,\circ)$ is coclosed relative to $\bfP=\FQA$ due to the point~(3) of Prop.~\ref{PropRelclosed} applied to the monoidal category $\bfC=(\QA,\circ)$. The corepresentations of comonoids $\OO\in\Comon(\QA,\circ)$ on $\B\in\FQA$ generalise the finite-dimensional representations of infinite-dimensional algebras at least partially.

\bb{Quantum semi-linear spaces.}
The algebraic set can be considered as a particular case the affine scheme. Categorically this fact is interpreted via a fully faithful functor $A\colon\AlgSet_\FF\hookrightarrow\CommAlg_\FF^\op$, where $\CommAlg_\FF$ is the category of the commutative algebras over $\FF$ and $A(X)$ is the algebra of regular $\FF$-valued functions on the algebraic set $X\in\AlgSet_\FF$. Hence $\AlgSet_\FF$ is contravariantly embedded into $\Alg_\FF$, the category of algebras over $\FF$. This is a standard way to interpret the algebras as non-commutative affine algebraic sets/varieties/schemes.

The contravariant functor $S^*\colon\FVect\to\FQA$ is extended to the contravariant faithful functor $S^*\colon\SLAlgSet_\FF\to\GrAlg$ as $S^*(X\times V)=A(X)\otimes_\FF SV^*=A(X)^e\otimes SV^*$, where $\gR^e=\gR\otimes_\FF\KK$ is the extension of scalars. In particular, an object $X\times\KK$ (that is an algebraic set $X\in\AlgSet_\FF$ embedded into $\SLAlgSet_\FF$) is lifted to the quantum level as  $A(X)\otimes_\FF\KK[u]$, the latter is a graded algebra $\A\in\GrAlg$ with components $\A_k\cong A(X)^e$. The contravariant functor $S^*$ has a colax monoidal structure $(\SLAlgSet_\FF,\otimes)\to(\GrAlg,\circ)$. Similarly, the functor $T^*$ is extended to a contravariant faithful strong monoidal functor $T^*\colon(\SLAlgSet_\FF,\otimes)\to(\GrAlg,\circ)$, $T^*(X\times V)=A(X)^e\otimes TV^*$ (if $\FF=\KK$, then $S^*$ and $T^*$ are fully faithful).

The images of $\SLAlgSet_\FF$ under the functors $S^*$ and $T^*$ lie in the category $\QAsc$ of semi-connected quadratic algebras. This is a full subcategory 
$\QAsc\subset\GrAlg$ consisting of the quadratic algebras $\gR\otimes\B$, where $\gR\in\Alg$, $\B\in\QA$. Moreover, these images lie in $\FQAsc$, the full subcategory of the quadratic algebras $\gR\otimes\B$, where $\gR\in\Alg$, $\B\in\FQA$. The objects of $\QAsc^\op$ or at least of $\FQAsc^\op$ can be considered as quantum analogues of the semi-linear algebraic sets over $\KK$ (see details in~\cite{Sqrt}).

It was shown in~\cite[\S~4.4.6]{Sqrt} that for any $\B\in\FQA$ the functor $-\circ\B\colon\QAsc\to\QAsc$ has a left adjoint. Due to the point~(3) of Prop.~\ref{PropRelclosed} the monoidal category $\bfC=(\QAsc,\circ)$ is coclosed relative to $\bfP=\FQA$. Since $\bfcohom(\A,\B)\in\FQA\subset\FQAsc$ for any $\A,\B\in\FQA$, Prop.~\ref{PropRelclosedhom} implies that the monoidal subcategory $\bfC'=(\FQAsc,\circ)$ is also coclosed relative to $\bfP=\FQA$.

Corepresentations of comonoids $\OO\in\Comon(\QAsc,\circ)$ on $\B\in\FQA$ are quantum analogues of representations of semi-linear algebraic monoids. In particular, they give a quantum version of the representations of the affine algebraic groups. A classical representation is considered as a corepresentation on a quadratic algebra via the functor~\eqref{FRepCon}, where $F\colon(\SLAlgSet_\FF,\otimes)\to(\QAsc,\circ)$ is $S^*$ or $T^*$. For the case of an algebraic group or monoid one need to compose these functors with the embedding $\AlgSet_\FF\hookrightarrow\SLAlgSet_\FF$. In the case of the strong monoidal functor $F=T^*$ we obtain a contravariant monoidal embedding of representations~\eqref{FRepConBB}.

One can lift the representations of groups, monoids $\MM\in\Mon(\Set,\times)$ and, more generally, semi-linear monoids $\SSS\in\Mon(\SLSet,\otimes)$ to the quantum level by the same formulae with $A(X)$ replaced by the algebra of all the functions $f\colon X\to\KK$ on $X\in\Set$.

\bb{The case of all graded algebras.}
The category $\QAsc$ is enough to lift the classical representations of algebraic monoids and (finite-dimensional) algebras to the quantum level, however, sometimes we need to consider more general $\NN_0$-graded and $\ZZ$-graded algebras (see e.g.~\cite[\S~6.2.2]{Sqrt}).

Let us extend the monoidal category $(\GrAlg,\circ)$ to the case of all the $\ZZ$-graded algebras. Denote it by $\ZZGrAlg$. The Manin product of the algebras $\A,\B\in\textbf{$\ZZ$-GrAlg}$ is the $\ZZ$-graded vector space with the components $(\A\circ\B)_k=\A_k\otimes\B_k$, $k\in\ZZ$, and multiplication $(a_k\otimes b_k)(a'_l\otimes b'_l)=(a_ka'_l)\otimes(b_kb'_l)\in\A_{k+l}\otimes\B_{k+l}$, where $a_k\in\A_k$, $a'_l\in\A_l$, $b_k\in\B_k$, $b'_l\in\B_l$.

In~\cite{Sqrt} we defined a {\it quantum representation} of a {\it quantum algebra} $\OO\in\Comon(\GrAlg,\circ)$ on $\B\in\FQA$ as a morphism $\omega\colon\bfcoend(\B)\to\OO$ in $\Comon(\GrAlg,\circ)$. More generally, we can define it as a morphism in $\Comon(\ZZGrAlg,\circ)$ from the same $\bfend(\B)$ to a comonoid $\OO\in\Comon(\ZZGrAlg,\circ)$. This notion can be interpreted in terms of Section~\ref{secGRT} as corepresentation (or as a representation in the opposite category) due to the following statement.

\begin{Th} \label{ThGrAlg}
 The monoidal category $\bfC=(\ZZGrAlg,\circ)$ is coclosed relative to $\bfP=\FQA$.
\end{Th}

\noindent{\bf Proof.} Note that the category $\QA$ is a monoidal subcategory of $(\ZZGrAlg,\circ)$ and that $(\QA,\circ)$ is coclosed relative to $\bfP=\FQA$. In order to apply Theorem~\ref{ThCoref} we need to check that $\bfC'=\QA$ is a coreflective subcategory of $\bfC=\ZZGrAlg$. For a $\ZZ$-graded algebra $\A$ denote by $G'\A$ the subalgebra of $\A$ generated by the subspace $\A_1$. This is an $\NN_0$-graded algebra $G'\A$ generated by its first order component $V=\A_1=(G'\A)_1$, so we have the graded epimorphism $a\colon TV\twoheadrightarrow G'\A$, where $TV\in\GrAlg$ is the tensor algebra. The kernel of $a$ is a graded ideal $I=\bigoplus\limits_{k\ge2}I_k$ in $TV$, where $I_k=I\cap V^{\otimes k}$, and we obtain the isomorphism $G'\A\cong TV/I$ in $\GrAlg$. Note that $TV/I$ is a quotient algebra of the quadratic algebra $TV/(I_2)\in\QA$, where $(I_2)$ is the ideal generated by $I_2\subset V^{\otimes2}$. Let $G\A=TV/(I_2)$ and $\varepsilon_\A\colon G\A\to\A$ be the composition
\begin{align}
\varepsilon_\A\colon G\A=TV/(I_2)\twoheadrightarrow TV/I\cong G'\A\hookrightarrow\A.
\end{align}
Its first graded component is the identity morphism $(\varepsilon_\A)_1=\id_V\colon V\isoright\A_1$.
Consider a quadratic algebra $\B=TW/(S)\in\QA$ and a graded homomorphism $f\colon\B\to\A$, where $W\in\Vect$ and $S$ is a subspace of $W^{\otimes 2}$. Since $\B$ is generated by $\B_1=W$, the image of $f$ lies in the subalgebra $G'\A\subset\A$, so it factors as in the commutative diagram
\begin{align} \label{BfGA}
\xymatrix{
\B\ar[rd]^{f}\ar@{-->}[d]_{f'} \\
 G'\A\ar@{^(->}[r]& \A
}
\end{align}
Such morphism $f'\colon\B\to G'\A$ is unique, since the embedding $G'\A\hookrightarrow\A$ is a monomorphism. Moreover, one can show that any $f'\colon\B\to TV/I$ factors as
\begin{align} \label{BTVI}
\xymatrix{
\B\ar[rd]^{f'}\ar@{-->}[d]_{h} \\
 TV/(I_2)\ar@{->>}[r]& TV/I
}
\end{align}
through a unique $h\colon\B\to TV/(I_2)$. Indeed, the first component of $TV/(I_2)\twoheadrightarrow TV/I$ is $\id_V$, so the commutativity of~\eqref{BTVI} is possible only if the first components of $f'$ and $h$ coincide. Since $(f'_1\otimes f'_1)S\subset I_2$, the formula $h_1=f'_1$ correctly defines the graded homomorphism $h\colon\B\to TV/(I_2)$ and the diagram~\eqref{BTVI} is commutative. Thus for any $f\colon\B\to\A$ there is a unique $h\colon\B\to G\A$ such that $\varepsilon_\A\cdot h=f$. By virtue of Theorem~\ref{ThUnit} we get a right adjoint $G\colon\ZZGrAlg\to\QA$ for the embedding $\QA\hookrightarrow\ZZGrAlg$ such that $h=Gf$. By taking the first graded components in the commutative diagrams~\eqref{BfGA}, \eqref{BTVI} we see that the first component $(Gf)_1\colon\B_1\to V$ coincides with $f_1\colon\B_1\to\A_1=V$. Since $\B$ is generated by $\B_1$, this formula uniquely defines the whole graded homomorphism $Gf\colon\B\to G\A=TV/(I_2)$.

Now let us prove that $\phi_{\A,\A'}=G(\varepsilon_\A\circ\varepsilon_{\A'})\colon G\A\circ G\A'\to G(\A\circ\A')$ is an isomorphism for arbitrary $\ZZ$-graded algebras $\A,\A'\in\ZZGrAlg$. Let $V=\A_1$, $V'=\A'_1$, the morphism $\phi_{\A,\A'}$ is uniquely defined by its first order component $(\phi_{\A,\A'})_1=\big(G(\varepsilon_\A\circ\varepsilon_{\A'})\big)_1=(\varepsilon_\A\circ\varepsilon_{\A'})_1=(\varepsilon_\A)_1\otimes(\varepsilon_{\A'})_1=\id_{V\otimes V'}$. Denote by $I=\bigoplus\limits_{k\ge2}I_k$, $I'=\bigoplus\limits_{k\ge2}I'_k$ the kernels of $a\colon TV\twoheadrightarrow G'\A$ and $a'\colon TV'\twoheadrightarrow G'\A'$ respectively, where $I_k=I\cap V^{\otimes k}$, $I'_k=I'\cap (V')^{\otimes k}$. The kernel of the composition $T(V\otimes V')\cong TV\circ TV'\xrightarrow{a\circ a'}G'\A\circ G'\A'\subset\A\circ\A'$ is a graded ideal $J=\bigoplus\limits_{k\ge2}J_k$ with components $J_k=J\cap (V\otimes V')^{\otimes k}$. The subalgebra $G'\A\circ G'\A'\subset\A\circ\A'$ coincides with $G'(\A\circ\A')$, hence $G(\A\circ\A')=T(V\otimes V')/(J_2)$.
By taking into account the isomorphisms $G'\A\cong TV/I$ and $G'\A'\cong TV'/I'$ we see that the ideal $J$ is the kernel of the composition
\begin{align} \label{TVVepi}
 T(V\otimes V')\cong TV\circ TV'\twoheadrightarrow TV/(I_2)\circ TV'/(I'_2)\twoheadrightarrow TV/I\circ TV'/I',
\end{align}
where the epimorphisms $\twoheadrightarrow$ are the Manin products of the corresponding projections. The second order component of the right epimorphism in~\eqref{TVVepi} is the identity map, so the kernel of the second order component of
\begin{align} \label{TVVTVTVII}
 T(V\otimes V')\cong TV\circ TV'\twoheadrightarrow TV/(I_2)\circ TV'/(I'_2)
\end{align}
equals $J_2$. Since~\eqref{TVVTVTVII} is a morphism of quadratic algebras, its kernel is the ideal $(J_2)$. Hence we obtain the isomorphism $TV/(I_2)\circ TV'/(I'_2)\cong T(V\otimes V')/(J_2)$, whose first order component is the identity map $\id_{V\otimes V'}$. This isomorphism is exactly the morphism $\phi_{\A,\A'}\colon G\A\circ G\A'\to G(\A\circ\A')$. In particular, $\phi_{\A,\A'}$ is an isomorphism for any $\A\in\ZZGrAlg$, $\A'\in\FQA$, so the conditions of Theorem~\ref{ThCoref} are fulfilled. \qed

Since $\bfcohom(\A,\B)\in\FQA\subset\GrAlg$ for any $\A,\B\in\FQA$, Prop.~\ref{PropRelclosedhom} implies that the category $(\GrAlg,\circ)$ is also coclosed relative to $\FQA$. The internal cohom-functor in $\bfC=(\GrAlg,\circ)$ or $(\ZZGrAlg,\circ)$ with parametrising subcategory $\bfP=\FQA$ equals to the composition of the functor $\bfcohom\colon\FQA^\op\times\FQA\to\FQA$ and the embedding $\FQA\hookrightarrow\GrAlg$ or $\FQA\hookrightarrow\ZZGrAlg$.

\section*{Conclusion}
\addcontentsline{toc}{section}{Conclusion}

\noindent{\bf General representation theory.}
Relative adjunctions introduced by Ulmer in~\cite{Ulm} can be used to generalise the notion of closed monoidal category in such way that the approach of~\cite[\S~3]{Sqrt} works literally for this more general case. Thus we obtain a general representation theory for relatively closed monoidal categories and for monoidal functors between such categories. In particular, we described the tensor product of representations and how it is translated by a monoidal functor. Any representation can be considered as a left action, however we proved all the statements about representations independently from the theory of actions, we used the natural transformations $\pi$ and $\Phi$.

\vspace*{2mm}\noindent{\bf Classical representations and semi-linear spaces.}
This general approach covers the classical types of representations on vector spaces: representations of algebras, groups, monoids, Lie groups and algebraic groups. To interpret representations of groups and monoids as representations in the general sense we introduced semi-linear spaces with some monoidal product. The monoids in the monoidal category of semi-linear spaces (the semi-linear monoids) generalise algebras, monoids and groups, so the representation theory of such monoids includes the classical representations. The same can be done for the groups and monoids with geometric structures.

\vspace*{2mm}\noindent{\bf Quantum Representation Theory.}
We proved that the categories of the $\NN_0$-graded and, more generally, of the $\ZZ$-graded algebras with the Manin product `$\circ$' are relatively closed, where the parametrising subcategory consists of the connected finitely generated quadratic algebras: $\bfP=\FQA$. This means that the general representation theory described here is applicable for quantum representations introduced in~\cite[\S~5.1]{Sqrt}. We expect that it will be also applied to generalise quantum representations for other cases such as super-case, topological case, corepresentations on the connected algebras, a quantum version of the projective representations.

\end{document}